\documentclass[11pt]{article}

\usepackage{amsfonts,amssymb,amstext}
\usepackage{amsmath}
\usepackage{amstext}
\usepackage{amsopn}
\usepackage{amsthm}

\makeatletter

\newcommand{\tr}{\hbox{tr}}
\newcommand{\n}{\noindent}
\newcommand{\q}{\quad}

\newcommand{\E}{\mathbb{E}}
\newcommand{\g}{\gamma}
\newcommand{\I}{\varphi}
\newcommand{\G}{\Gamma}
\newcommand{\de}{\delta}
\newcommand{\De}{\Delta}

\newcommand{\al}{\alpha}
\newcommand{\la}{\lambda}
\newcommand{\f}{\infty}
\newcommand{\vs}{\varepsilon}
\newcommand{\cd}{\cdot}
\newcommand{\si}{\sigma}

\newcommand{\be}{\beta}

\newcommand{\Om}{\Omega}
\newcommand{\om}{\omega}

\newcommand{\sm}{\setminus}
\newcommand{\inl}{\int_{-\pi}^{\pi}}
\newcommand {\ol} {\overline}

\newcommand {\uu}{{\bf u}}
\newcommand {\vv}{{\bf v}}

\newcommand{\ConvD}{\overset{d}{\rightarrow}}

\newcommand{\ConvFDD}{\overset{f.d.d.}{\longrightarrow}}

\newcommand{\EqFDD}{\overset{f.d.d.}{=}}
\newcommand{\Var}{\mathrm{Var}}
\newcommand{\Cum}{\mathrm{Cum}}
\newcommand{\Cov}{\mathrm{Cov}}

\newcommand {\sn} {\subsection}
\newcommand {\ssn} {\subsubsection}

\newtheorem{thm}{Theorem}[section]

\newtheorem{pp}{Proposition}[section]
\newtheorem{cor}{Corollary}[section]
\newtheorem{exa}{Example}[section]
\newtheorem{rem}{Remark}[section]

\newtheorem{asn}{Assumption}[section]

\usepackage{enumerate}

\numberwithin{equation}{section}

\newcommand{\beq}{\begin{equation}}
\newcommand{\eeq}{\end{equation}}
\newcommand{\bea}{\begin{eqnarray}}
\newcommand{\eea}{\end{eqnarray}}
\newcommand{\beaa}{\begin{eqnarray*}}
\newcommand{\eeaa}{\end{eqnarray*}}
\newcommand{\wh}{\widehat}
\usepackage[active]{srcltx}

\begin{document}

\title{A Survey on Limit Theorems for Toeplitz Type Quadratic Functionals
of Stationary Processes and Applications}

\author{Mamikon S. Ginovyan and Murad S. Taqqu}


\date{\today}

\maketitle

\begin{abstract}
\noindent
This is a survey of recent results on central and non-central
limit theorems for quadratic functionals of stationary processes.
The underlying processes are  Gaussian, linear or L\'evy-driven
linear processes with memory, and are defined either in
discrete or continuous time.
We focus on limit theorems for Toeplitz and tapered Toeplitz type quadratic
functionals of stationary processes with applications in parametric and
nonparametric statistical estimation theory.
We discuss questions concerning Toeplitz matrices and operators,
Fej\'er-type singular integrals, and L\'evy-It\^o-type and Stratonovich-type
multiple stochastic integrals. These are the main tools for obtaining
limit theorems.

\end{abstract}

\tableofcontents

\section{Introduction}
\label{int}
A significant part of large-sample statistical inference relies on limit
theorems of probability theory, which involves sums and quadratic functionals
of stationary observations. Depending on the memory (or dependence) structure
of the underlying processes, these functionals, once normalized, can have different
limits, and the proofs of such limit theorems generally use different methods.
In this paper, we focus on the quadratic functionals. The term 'central
limit theorem' (CLT) refers to a statement that a suitably standardized
quadratic functional converges in distribution to a Gaussian random variable.
Limit theorems where a suitably standardized quadratic functional converges
in distribution to a non-Gaussian random variable are termed 'non-central
limit theorems' (NCLT).

We present results on central and non-central
limit theorems for Toeplitz and tapered Toeplitz type quadratic
functionals of stationary processes with applications in parametric and
nonparametric statistical estimation theory.
The underlying processes are  Gaussian, linear or L\'evy-driven
linear processes with memory, and are defined either in
discrete or continuous time.
We also discuss some questions concerning Toeplitz matrices and operators,
Fej\'er-type singular integrals, L\'evy-It\^o-type and Stratonovich-type
multiple stochastic integrals, and power counting theorems.
These are the main tools for obtaining
limit theorems, but they are also of interest in themselves.

\sn{Notation and conventions.}
The following notation and conventions are used throughout the paper.

The symbol '$:=$' stands for 'by definition'. c.t.: = continuous-time;
d.t.: = discrete-time; s.d.:= spectral density; c.f.:= covariance function;
CLT:= central limit theorem; NCLT:= non-central limit theorem.
The symbol '$\stackrel{d}=$' stands for equality of the
finite-dimensional distributions.
The symbol '$\ConvD$' stands for convergence in distribution.
The symbol '$\ConvFDD$' stands for convergence of finite-dimensional distributions.
The symbol '$\Longrightarrow$' stands for weak convergence.
The notation $X_T\ConvD \eta \sim N(0,\sigma^2)$ as $T\to\f$
will mean that the distribution of the random variable
$X_T$ tends (as $T\to\f$) to the centered normal
distribution with variance $\sigma^2$.
$\E[\cdot]$: = expectation operator.
$\tr[A]$: = trace of an operator (matrix) $A$.
${\mathbb I}_A(\cdot)$: = indicator of a set $A\subset \Lambda$.
WN$(0,1)$: = standard white-noise.

The standard symbols $\mathbb{N}$, $\mathbb{Z}$ and $\mathbb{R}$
denote the sets of natural, integer and real numbers, respectively.
By $\Lambda$ we denote the frequency domain, that is, $\Lambda:=\mathbb{R}$
in the c.t.\ case, and $\Lambda:= [-\pi.\pi]$ in the d.t.\ case.
By $L^p(\mu):=L^p(\Lambda,\mu)$ ($p\geq $1) we denote the weighted
Lebesgue space with respect to the measure $\mu$, and by $||\cdot||_{p,\mu}$
we denote the norm in $L^p(\mu)$.
In the special case where $d\mu(\la)=d\la$, we will use 
$L^p$ and $||\cdot||_{p}$, respectively.
By $l^p$ ($p\geq $1) we denote the the space of $p$-summable sequences.
The letters $C$ and $c$ with or without indices
are used to denote positive constants, the values of which can vary from
line to line.

\sn{The functionals under consideration}
Let $\{X(u), \ u\in \mathbb{U}\}$ be a centered real-valued stationary process
with spectral density (s.d.)\ $f(\la)$, $\la\in \Lambda$ and covariance function (c.f.)\
$r(t)$, $t\in \mathbb{U}$. We consider simultaneously the continuous-time (c.t.)\
case, where $\mathbb{U}=\mathbb{R}:=(-\f,\f)$, and the discrete-time (d.t.)\ 
case, where $\mathbb{U}=\mathbb{Z}:=\{0,\pm 1,\pm 2,\ldots \}$.
The domain $\Lambda$ of the frequency variable $\la$ is $\Lambda=\mathbb{R}$
in the c.t.\ case, and $\Lambda:= [-\pi.\pi]$ in the d.t.\ case.

We first survey the recent results concerning the asymptotic distribution
(as $T\to\f$) of the following Toeplitz type quadratic functionals of the process $X(u)$:
\beq
\label{MTc-1}
Q_T:=
 \left \{
 \begin{array}{ll}
\sum_{t=1}^T\sum_{s=1}^T\widehat g(t-s)X(t)X(s) & \mbox{in the d.t.\  case},\\
\\
\int_0^T\int_0^T\widehat g(t-s)X(t)X(s)\,dt\,ds&  \mbox{in the c.t.\ case},
\end{array}
\right. \eeq
where
\beq
\label{c-2}
\widehat g(t):= \int_\Lambda e^{i\la t}\,g(\la)\,d\la, \quad t\in \mathbb{U}
\eeq
is the Fourier transform of some real, even, integrable function
$g(\la),$ $\la\in\Lambda$. We will refer to $g(\la)$ and to its Fourier transform
$\widehat g(t)$ as a {\it generating function} and {\it generating kernel} for the functional
$Q_T$, respectively.
In the d.t.\  case the functions $f(\la)$ and $g(\la)$ are assumed to be $2\pi$-periodic
and periodically extended to $\mathbb{R}$.
In the c.t.\  case the process $X(u)$ is assumed measurable and mean square continuous.

The limit distributions of the functionals in (\ref{MTc-1}) are completely determined
by the spectral density $f(\la)$ (or covariance function $r(t)$) and the generating
function $g(\la)$  (or generating kernel $\widehat g(t)$), and depending
on their properties, the limits can be either Gaussian (that is, $Q_T$ with an appropriate
normalization obeys the central limit theorem), or non-Gaussian.

The following two questions arise naturally:

\begin{itemize}
\item[(a)]
Under what conditions on $f(\la)$ (resp. $r(t)$) and $g(\la)$ (resp. $\widehat g(t)$)
will the limits be Gaussian? (CLT-problem).\\
\item[(b)] Describe the limit distributions, if they are non-Gaussian (NCLT-problem).
\end{itemize}

We discuss here these questions both for d.t.\  and c.t.\  stationary processes,
and survey the recent results.

We also survey recent results concerning functional central and non-central
limit theorems for the following processes, generated by quadratic functionals in \eqref{MTc-1}:

\beq
\label{MTc-1p}
Q_T(\tau):=
 \left \{
 \begin{array}{ll}
\sum_{t=1}^{[T\tau]}\sum_{s=1}^{[T\tau]}\widehat g(t-s)X(t)X(s) &
\mbox{in the d.t.\  case},\\
\\
\int_0^{T\tau}\int_0^{T\tau}\widehat g(t-s)X(t)X(s)\,dt\,ds&
\mbox{in the c.t.\  case},
\end{array}
\right.
\eeq
where $\tau\in[0,1]$ and $[\,\cdot\,]$ stands for the greatest integer.

We say that a functional central limit theorem (FCLT) for $Q_T(\tau)$ holds if
the process $Q_T(\tau)$ with an appropriate normalization converges weakly in
$C[0,1]$ in the c.t.\  case (and in $D[0,1]$ in the d.t.\  case) to Brownian motion.
We say that we have a functional non-central limit theorem (FNCLT) for the process
$Q_T(\tau)$ if the limit is non-Gaussian or, if Gaussian, it is not Brownian motion.

\vskip2mm
\n
{\em The tapered case.}
In the spectral analysis of stationary processes the data are frequently
tapered before calculating the statistics of interest. Instead of the original data
$\{X(t)$, $0\le t\le T\}$ the tapered data $\{h(t)X(t)$, $0\le t\le T\}$
with the data taper $h(t)$ are used for all further calculations.
The benefits of tapering the data have been widely reported in the literature
(see, e.g., Bloomfield \cite{Bl}, Brillinger \cite{Bri2},  Dahlhaus \cite{D1,D3, D2},
Dahlhaus and K\"unsch \cite{DK}, Guyon \cite{Gu}, and references therein).
For example, data-tapers are introduced to reduce the so-called 'leakage effects',
that is, to obtain better estimation of the spectrum of the model in the case
where it contains high peaks.
Other application of data-tapers is in situations in which some of the
data values are missing. Also, the use of tapers leads to bias reduction,
which is especially important when dealing with spatial data. In this case,
the tapers can be used to fight the so-called 'edge effects'.

In this case, to establish asymptotic properties of the corresponding estimators,
we have to study the asymptotic distribution (as $T\to\f$) of
the following Toeplitz type tapered quadratic functionals of the process $X(u)$:
\beq
\label{MTc-1T}
Q_T^h:=
 \left \{
 \begin{array}{ll}
\sum_{t=1}^T\sum_{s=1}^T\widehat g(t-s)h_T(t)h_T(s)X(t)X(s) & \mbox{in the d.t.\  case},\\
\\
\int_0^T\int_0^T \widehat g(t-s)h_T(t)h_T(s)X(t)X(s)\,dt\,ds&  \mbox{in the c.t.\  case},
\end{array}
\right. \eeq
where $\widehat g(t)$, $t\in\mathbb{U}$
is the Fourier transform of some integrable even function $g(\la)$,
$\la\in\Lambda$, and $h_T(t):= h(t/T)$ with a taper function $h(t)$, $t\in[0,1]$.

Quadratic functionals of the form \eqref{MTc-1} and \eqref{MTc-1T}
arise naturally in the context of nonparametric and parametric statistical
estimation of the spectrum of stationary processes based on the non-tapered and
tapered data, respectively. And their limiting distributions are necessary
to establish asymptotic properties of the corresponding estimators.
For instance, when we are interested in nonparametric estimation of a linear
integral functional $J(f)$ in $L^p(\Lambda)$, $p > 1$, then a natural
statistical estimator for $J(f)$ is the linear integral functional of the
empirical spectral density (periodogram) of the process $X(t)$, which is
a quadratic functional of the form \eqref{MTc-1} in the non-tapered case,
and of the form \eqref{MTc-1T} in the tapered case.
In the case of parametric estimation, for example, the Whittle estimation procedure
is based on the minimization of quadratic functionals of the form \eqref{MTc-1} and \eqref{MTc-1T}.

\sn{A brief history}
The problem of describing the asymptotic distribution of Toeplitz type
quadratic forms and functionals of stationary processes has a long history,
and goes back to the classical monograph by Grenander and Szeg\H{o} \cite{GS},
where the CLT-problem for Gaussian processes was considered as an application
of the authors' theory of the asymptotic behavior of the trace of products of
truncated Toeplitz matrices and operators.

Later the CLT-problem has been studied by a number of authors. Here we mention only
some significant contributions. For the d.t.\  short memory processes,
the problem was studied by Bentkus \cite{B}, Ibragimov \cite{I1963} and
Rosenblatt \cite{R2}, in connection with statistical estimation of the spectral
and covariance functions.
Since 1986, there has been a renewed interest in this problem, related to the
statistical inferences for long memory (long-range dependence) and intermediate
memory (anti-persistent)  processes (see, e.g., Avram \cite{A},
Fox and Taqqu \cite{FT1}, Giraitis and Surgailis \cite{GSu}, Giraitis and Taqqu
\cite{GT1998}, Has'minskii and Ibragimov \cite{IH1}, Ginovian and Sahakian
\cite{GS2005}, Terrin and Taqqu \cite{TT2}, and references therein).
In particular, Avram \cite{A}, Fox and Taqqu \cite{FT1}, Ginovian and Sahakian
\cite{GS2005}, Giraitis and Surgailis \cite{GSu}, Giraitis and Taqqu \cite{GT1998}
have obtained sufficient conditions for quadratic form $Q_T$ to obey the central
limit theorem.

In the case of c.t.\ stationary Gaussian processes the CLT-problem for Toeplitz type
quad\-ratic functionals was studied in a number of papers. We cite merely the papers
Avram et al. \cite{AvLS}, Bai et al. \cite{BGT1, BGT2}, Bryc and Dembo \cite{BD},
Ginovyan \cite{G1988b, G1994, G2011a}, Ginovyan and Sahakyan \cite{GS2007},
Ibragimov \cite{I1963}, Leonenko and Sakhno \cite{Leonenko:2006},
where additional references can be found.
The NCLT-problem have been studied in Bai et al. \cite{BGT2},
Giraitis and Taqqu \cite{GT2001}, and Terrin and Taqqu \cite{TT1}.

Central and non-central limit theorems for tapered quadratic forms of a
d.t.\  long memory
Gaussian stationary fields have been proved in Doukhan et al. \cite{DLS}.
A central limit theorem for tapered quadratic functionals $Q_T^h$,
in the case where the underlying model $X(t)$ is a L\'evy-driven c.t.\
stationary linear process has been proved in Ginovyan and Sahakyan \cite{GS2019}.

The problem of nonparametric and parametric estimation of the spectrum of the process $X(t)$
both for d.t.\  and c.t.\  cases based on the non-tapered data has been considered by many authors
(see, e.g., Avram et al. \cite{AvLS}, Bentkus \cite{B}, Dahlhaus \cite{Dah1}, Dzhaparidze
\cite{Dah1}, Fox and Taqqu \cite{FT2}, Gao et al. \cite{GAHT}, Ginovyan \cite{G1988a, G1988b, G1995, G2003, G2011a,
G2011b}, Giraitis et al. \cite{GKSu}, Giraitis and Surgailis \cite{GSu}, Giraitis and Taqqu \cite{GT1999b}, Guyon \cite{Gu},
Hasminskii and Ibragimov \cite{IH1}, Ibragimov \cite{I1963, I1967}, Leonenko and Sakhno \cite{Leonenko:2006},
Taniguchi \cite{Tan87}, Taniguchi and Kakizawa \cite{TK}, Taqqu \cite{Tq2}, and references therein).
The problem in the tapered case was studied in Alomari et al. \cite{ALRST},
Dahlhaus \cite{D1,D2}, Dahlhaus and K\"unsch \cite{DK}, Ginovyan \cite{G2020e},
Ginovyan and Sahakyan \cite{GS2019}, and Lude\~na and Lavielle \cite{LL}.

\sn{Frequency-domain conditions}
Conditions that are expressed in terms of the spectral density $f(\la)$ and the generating function
$g(\la)$ will be termed frequency-domain conditions, while conditions that are in terms
of the covariance function $r(t)$ and the generating kernel $\widehat g(t)$ will be termed
time-domain conditions.

There are three sets of frequency-domain conditions for functionals
of the form (\ref{MTc-1}) to obey the CLT, and these conditions separate the ranges
of CLT and NCLT:
\begin{itemize}
\item[(a)] the $(L^p, L^q)$ condition,
\item[(b)] the $(\al,  \be)$ condition,
\item[(c)] the trace condition.
\end{itemize}
All three are compensation conditions, meaning that the good behavior of one function,
say $g$ can compensate for the bad behavior of $f$ and vice versa.

\vskip2mm
\n
{\bf(a) The $(L^p, L^q)$ condition.} Let $f \in L^p$ $(p\ge2)$ and
$g\in L^q$ $(q\ge2)$. If $1/p+1/q\le1/2,$ then the functional $Q_T$
with an appropriate normalization obeys the CLT (see Theorem \ref{cth1}(C)),
while when $1/p+1/q>1/2$, then, in general, $Q_T$ does not obey the CLT.
This condition goes back to the classical works of Rosenblatt \cite{R1,R2},
where when estimating the covariance function $r(t)$ of a d.t.\  Gaussian process,
it was shown that for a sample covariance $\widehat r(t)$  (which is a functional of the
form (\ref{MTc-1}) with $g(\la)=\cos(t\la)$) to obey the CLT, the spectral density $f$
should satisfy the condition $f \in L^p$ $(p\ge2)$.

\medskip
\n
{\bf (b) The $(\al,  \be)$ condition.} If both the spectral density $f$ and the generating
function $g$ are regularly varying at the origin of orders $\alpha$ and $\beta$, respectively,
then it is the sum $\alpha+\beta$ that determines the limiting distribution of $Q_T$.
When $\alpha+\beta\le 1/2,$ then the limiting distribution of $Q_T$ is Gaussian, that is,
$Q_T$ with an appropriate normalization obeys the CLT (see Theorem \ref{cth2}),
while when $\alpha+\beta>1/2,$ then the limiting distribution of $Q_T$ is non-Gaussian
(see Theorem \ref{Thm:fNCLT}).

\medskip
\n
{\bf (c) The trace condition.} This condition, which is more general and implies both conditions
(a) and (b), is an implicit condition. It is expressed in terms of traces of products of
truncated Toeplitz matrices (in the d.t.\  case) and operators (in the c.t.\  case).
The idea here is to approximate the traces of products of Toeplitz matrices and operators
(which are no longer Toeplitz) by the traces of a Toeplitz matrix and a Toeplitz operator,
respectively.
Let $A_T(f)$ denote either the $T\times T$ Toeplitz matrix $B_T(f)$ or the $T$-truncated
Toeplitz operator $W_T(f)$ generated by the spectral density $f$,
and let $A_T(g)$ denote either the $T\times T$ Toeplitz matrix $B_T(g)$, or the
$T$-truncated Toeplitz operator $W_T(g)$ generated by the function $g$
(for definitions see Section \ref{fdc}). The trace condition is: if
$$fg\in L^2$$
and
$$T^{-1}\bigl[\tr\bigl(A_T(f)A_T(g)\bigr)^2 -\tr \bigl(A_T(f^2g^2)\bigr)\bigr]\to 0
\q {\rm as} \q T\to\f,$$
then the quadratic functional $Q_T$ in (\ref{MTc-1}) with an appropriate
normalization obeys the CLT (see Theorems \ref{cth1}(A) and \ref{GSth1}).

We will also discuss the time-domain counterparts of $(L^p, L^q)$ and $(\al,  \be)$ conditions.

\sn{Methods and tools.}
The most commonly used methods to prove central limit theorems for Toeplitz type
quadratic forms and functionals are:
\begin{itemize}
\item[(a)] the method of cumulants or moments,
\item[(b)] the approximation method,
\item[(c)] the method of characteristic functions.
\end{itemize}

To prove the non-central limit theorems for Toeplitz type
quadratic functionals, was used:
\begin{itemize}
\item[(a)] the spectral representation of the underlying process,
\item[(b)] the properties of L\'evy-It\^o-type and Stratonovich-type multiple
stochastic integrals,
\item[(c)] the power counting theorem.
\end{itemize}

Some details of the above methods are described in Section \ref{methods}.

\sn{The structure of the paper}
The paper is structured as follows.
In Section \ref{model} we describe the model of interest - a stationary process,
and recall some key notions and results from the theory of stationary processes.
In Section \ref {G-CLT} we present sufficient conditions for Toeplitz type quadratic
forms and functionals of the form \eqref{MTc-1} to obey the CLT in the case where
the model is either a Gaussian or a linear process.
Section \ref{FLTGL} contains functional central and noncentral limit theorems for
processes generated by Toeplitz type quadratic forms and functionals for Gaussian
and linear models.
Section \ref{FLTL} is devoted to the functional central and noncentral limit
theorems for L\'evy-driven linear models.
In Section \ref{Tap1} we discuss the case of tapered Toeplitz quadratic functionals,
and state central limit theorems.
Section \ref{App} contains some applications, involving nonparametric estimation
of spectral functionals and Whittle parametric estimation procedure.
In Section \ref{methods} we briefly discuss the methods and tools, used to prove
central and noncentral limit theorems for Toeplitz type quadratic forms and functionals.


\section{The model: second-order stationary process}
\label{model}
In this section we introduce the model of interest - a second-order stationary process,
and recall some key notions and results from the theory of stationary processes.

\sn{Key notions and some basic results}

\ssn{Second-order (wide-sense) stationary process}
Let $\{X(u), \ u\in \mathbb{U}\}$ be a centered real-valued second-order (wide-sense)
stationary process defined on a probability space $(\Omega, \mathcal{F}, P)$ with covariance
function $r(t)$, that is,
\[
\E[X(u)]=0, \q r(u)=\E[X(t+u)X(t)], \q u, t\in\mathbb{U},
\]
where $\E[\cdot]$ stands for the expectation operator with respect to measure $P$.
We consider simultaneously the c.t.\  case, where
$\mathbb{U}=\mathbb{R}:=(-\f,\f)$, and the d.t.\  case, where
$\mathbb{U}=\mathbb{Z}:=\{0,\pm 1,\pm 2,\ldots \}$.
We assume that $X(u)$ is a {\it non-degenerate process}, that is,
${\rm Var}[X(u)]=\E|X(u)|^2=r(0)>0$.
(Without loss of generality, we assume that $r(0)=1$).
In the c.t.\  case the process $X(u)$ is also assumed
mean-square continuous, that is, $\mathbb{E}[X(t)-X(s)]^2\to0$ as
$t\to s$.
This assumption is equivalent to that of the covariance function $r(u)$
be continuous at $u=0$ (see, e.g., Cram\'er and Leadbetter \cite{CL}, Section 5.2).

\ssn{Spectral representations}
By the Herglotz theorem in the d.t.\  case, and the Bochner-Khintchine
theorem in the c.t.\  case (see, e.g., Cram\'er and Leadbetter \cite{CL},
Doob \cite{Doob}, Ibragimov and Linnik \cite{IL}),
there is a finite measure $\mu$ on $(\Lambda, \mathfrak{B}(\Lambda))$,
where $\Lambda=\mathbb{R}$ in the c.t.\  case, and $\Lambda= [-\pi.\pi]$
in the d.t.\  case, and $\mathfrak{B}(\Lambda)$ is the Borel $\si$-algebra on $\Lambda$,
such that for any $t\in \mathbb{U}$ the covariance function $r(t)$
admits the following {\sl spectral representation}:
\beq
\label{i1}
r(u)=\int_\Lambda \exp\{i\la u\}d\mu(\la), \q u\in\mathbb{U}.
\eeq
\n

The measure $\mu$ in (\ref{i1}) is called the {\sl spectral measure} of the
process $X(u)$. The function  $F(\la): =\mu[-\pi,\la]$ in the d.t.\  case
and $F(\la): =\mu[-\f,\la]$ in the c.t.\  case, is called the {\sl spectral
function} of the process $X(t)$.
If $F(\la)$ is absolutely continuous (with respect to Lebesgue measure),
then the function $f(\la):=dF(\la)/d\la$ is called the {\sl spectral density}
of the process $X(t)$. Notice that if the spectral density $f(\la)$ exists, then
$f(\la)\geq 0$, $f(\la)\in L^1(\Lambda)$, and \eqref{i1} becomes
\begin{equation}
\label{mo1}
r(u)=\int_\Lambda\exp\{i\la u\}f(\la)d\la, \q  u\in\mathbb{U}.
\end{equation}
Thus, the covariance function $r(u)$ and the spectral
function $F(\la)$ (resp. the spectral density function $f(\la)$) are equivalent
specifications of the second order properties for a stationary process
$\{X(u), \ u\in\mathbb{U}\}$.

By the well-known Cram\'er theorem (see, e.g., Cram\'er and Leadbetter \cite{CL})
for any stationary process $\{X(u), \, u\in \mathbb{U}\}$ with spectral
measure $\mu$ there exists an orthogonal stochastic measure $Z=Z(B)$,
$B\in\mathfrak{B}(\Lambda)$, such that for every $u\in \mathbb{U}$
the process $X(u)$ admits the following {\sl spectral representation}:
\begin{equation}
\label{i18}
X(u)=\int_\Lambda \exp\{i\la u\}dZ(\la), \q  u\in\mathbb{U}.
\end{equation}
Moreover, $E\left[|Z(B)|^2\right]=\mu(B)$ for every $B\in\mathfrak{B}(\Lambda)$.
For definition and properties of orthogonal stochastic measures and
stochastic integral in (\ref{i18}) we refer, e.g., Cram\'er and Leadbetter \cite{CL}.

\ssn{Kolmogorov's isometric isomorphism theorem}
Given a probability space $(\Om, \mathcal{F}, P)$, define the $L^2$-space of
random variables $\xi= \xi(\om)$,
$\E[\xi]=0$:
\begin{equation}
\label{i19}
L^2(\Om):=\left\{\xi: \, ||\xi||^2:=
\int_\Om |\xi(\om)|^2dP(\om)<\f\right\}.
\end{equation}
Then $L^2(\Om)$ becomes a Hilbert space with the following
inner product: for $\xi, \eta\in L^2(\Om)$
\begin{equation}
\label{i20}
(\xi,\eta)= \E[\xi \eta]=\int_\Om \xi(\om){\eta(\om)}dP(\om).
\end{equation}
\n
For $a, b\in \mathbb{U}$,  $-\f\le a\le b\le \f,$
we define  the space $H_a^b(X)$ to be the closed linear
subspace of the space $L^2(\Om)$ spanned by the random
variables $X(u,\om)$, $u\in [a,b]$:
\begin{equation}
\label{i21}
 H_{a}^b(X): = \ol{sp}\{X(u), \,\, a \le u\le b\}_{L^2(\Om)}.
\end{equation}

\n Observe that the subspace $H_{a}^b(X)$ consists of
all finite linear combinations $\sum_{k=1}^n c_k X(u_k)$
($a \le u_k\le b$), as well as, their $L^2(\Om)$-limits.\\
The space $H(X):=H_{-\f}^\f(X)$ is called the {\it Hilbert
space generated by the process} $X(u)$, or the {\it time-domain} of $X(u)$.

Let $\mu$ be the spectral measure of the process $\{X(u), \, u\in \mathbb{U}\}$.
Consider the weighted $L^2$-space $L^2(\mu):=L^2(\mu,\Lambda)$ of complex-valued functions
$\I(\la), \, \la\in\Lambda$, defined by
\beq
\label{i22}
L^2(\mu):=\left\{\I(\la): \,||\I||^2_\mu:= \int_\Lambda|\I(\la)|^2d\mu(\la)<\f\right\}.
\eeq
Then $L^2(\mu)$ becomes a Hilbert space with the following
inner product: for $\I, \psi\in L^2(\mu)$
\begin{equation}
\label{i23}
(\I, \psi)_\mu= \int_\Lambda\I(\la)\ol{\psi}(\la)d\mu(\la).
\end{equation}
The Hilbert space $L^2(\mu,\Lambda)$ is called the {\it frequency-domain} of the process $X(u)$.


{\sl Kolmogorov's isometric isomorphism theorem} states that for any stationary process
$X(u)$, $u\in \mathbb{U}$, with spectral measure $\mu$ there exists a unique
isometric isomorphism $V$ between the time- and frequency-domains $H(X)$ and $L^2(\mu)$,
such that $V[X(u)]=e^{iu\la}$ for any $u\in \mathbb{U}.$

Thus, any linear problem in the time-domain $H(X)$
can be translated into one in the frequency-domain
$L^2(\mu)$, and vice versa. This fact allows to
study stationary processes using analytic methods.

\sn{Linear processes. Existence of spectral density functions}
We will consider here stationary processes possessing spectral density
functions. For the following results we refer to
Cram\'er and Leadbetter \cite{CL}, Doob \cite{Doob}, Ibragimov and Linnik \cite{IL}.
\begin{itemize}
\item[(a)]
The spectral function $F(\la)$ of a d.t.\  stationary process $\{X(u), \,u \in \mathbb{Z}\}$
is absolutely continuous (with respect to the Lebesgue measure), $F(\la)=\int_{-\pi}^\la f(x)dx$,
if and only if it can be represented as an infinite moving average:
\begin{equation}
\label{dlp}
X(u) = \sum_{k=-\f}^{\f}a(u-k)\xi(k), \quad
         \sum_{k=-\f}^{\f}|a(k)|^2 < \f,
\end{equation}
where $\{\xi(k), k\in \mathbb{Z}\}\sim$ WN(0,1) is a standard white-noise, that is,
a sequence of orthonormal random variables.
\item[(b)]
The covariance function $r(u)$ and the spectral density $f(\la)$ of $X(u)$ are given by formulas:
\begin{equation}
\label{dcv}
r(u)= \E X(u)X(0)=\sum_{k=-\f}^{\f}a(u+k) a(k),
\end{equation}
and
\beq
\label{dsd}
f(\la) = \frac{1}{2\pi} |\wh a(\la)|^2 =
\frac{1}{2\pi}\left|\sum_{k=-\f}^{\f}a(k)e^{-ik\la}\right|^2,
\quad \la\in\Lambda.
\eeq
\item[(c)]
In the case where $\{\xi(k), k\in \mathbb{Z}\}$ is a sequence of Gaussian random variables,
the process $\{X(u), \,u \in \mathbb{Z}\}$ is Gaussian.
\end{itemize}
Similar results hold for c.t.\  processes. Indeed, the following results hold.
\begin{itemize}
\item[(a)]
The spectral function $F(\la)$ of a c.t.\  stationary process $\{X(u), \,u \in \mathbb{R}\}$
is absolutely continuous (with respect to Lebesgue measure), $F(\la)=\int_{-\f}^\la f(x)dx$,
if and only if it can be represented as an infinite continuous moving average:
\begin{equation}
\label{clp}
X(u)=\int_{\mathbb{R}} a(u-t)d\xi(t), \q
         \int_{\mathbb{R}}|a(t)|^2dt < \f,
\end{equation}
where $\{\xi(t), t\in \mathbb{R}\}$ is a process with orthogonal increments and
$\E|d\,\xi(t)|^2 = dt$.
\item[(b)]
The covariance function $r(u)$ and the spectral density
$f(\la)$ of $X(u)$ are given by formulas:
\begin{equation}
\label{ccv}
r(u)= \E X(u)X(0)=\int_{\mathbb{R}} a(u+x)a(x)dx,
\end{equation}
and
\beq
\label{csd}
f(\la) = \frac{1}{2\pi} |\wh a(\la)|^2 =
\frac{1}{2\pi}\left|\int_{\mathbb{R}} e^{-i\la t}a(t)dt\right|^2,
\quad \la\in\mathbb{R}.
\eeq
\item[(c)]
In the case where $\{\xi(t), t\in \mathbb{R}\}$ is a Gaussian process,
the process $\{X(u), \,u \in \mathbb{Z}\}$ is Gaussian.
\end{itemize}

\sn{L\'evy-driven linear process}
We first recall that a L\'evy process, $\{\xi(t),\ t\in\mathbb{R}\}$ is a process with
independent and stationary increments, continuous in probability, with sample-paths which
are right-continuous with left limits (c$\grave{a}$dl$\grave{a}$g) and
$\xi(0) = \xi(0-) = 0$.
The Wiener process $\{B(t), \ t\ge 0\}$ and the centered Poisson process
$\{N(t)-\E N(t), \ t\ge 0\}$ are typical examples of centered L\'evy processes.
A L\'evy-driven linear process $\{X(t),\ t\in\mathbb{R}\}$ is a real-valued c.t.\
stationary process defined by \eqref {clp},
where $\xi(t)$ is a L\'evy process satisfying the conditions:
\[
\text{$\E \xi(t)=0$, $\E \xi^2(1)=1$ and $\E\xi^4(1)<\infty$.}
\]

In the case where $\xi(t)=B(t)$, $X(t)$ is a Gaussian process.

The function $a(\cdot)$ in representations \eqref{dlp} and \eqref{clp}
plays the role of a {\it time-invariant filter}, and the linear processes
defined by \eqref{dlp} and \eqref{clp} can be viewed as the output of
a linear filter $a(\cdot)$ applied to the process $\{\xi(u),\ t\in\mathbb{U}\}$,
called the innovation or driving process of $X(u)$.

Processes of the form (\ref{clp}) appear in many fields of science
(economics, finance, physics, etc.), and cover a large class of popular models
in c.t.\  time series modeling.
For instance, the so-called c.t.\  autoregressive moving average (CARMA) models,
which are the c.t.\  analogs of the classical autoregressive moving average (ARMA)
models in d.t.\  case, are of the form (\ref{clp}) and play a
central role in the representations of c.t.\  stationary time series
(see, e.g., Brockwell \cite{Br1}).

\sn{Dependence (memory) structure of the model}
In the frequency domain setting, the statistical and spectral analysis
of stationary processes requires
{\em two types of conditions\/} on the spectral density $f(\la).$
The first type controls the {\em singularities} of $f(\la)$,
and involves the {\em dependence (or memory) structure } of the
process, while the second type -- controls the {\em smoothness} of
$f(\la).$

We will distinguish the following types of stationary models:

(a) short memory (or short-range dependent),

(b) long memory (or long-range dependent),

(c) intermediate memory (or anti-persistent).

\noindent
The memory structure of a stationary process is essentially a
measure of the dependence between all the variables in the process,
considering the effect of all correlations simultaneously.
Traditionally memory structure has been defined in the time domain
in terms of decay rates of the autocorrelations, or in the
frequency domain in terms of rates of explosion of low frequency spectra
(see, e.g., Beran \cite{Ber}, Beran et al. \cite{BFGK}, Giraitis et al. \cite{GKSu},
Gu\'egan \cite{Gu1}, Robinson \cite{Ro2}, and references therein).

It is convenient to characterize the memory structure in terms
of the spectral density function.

\ssn{Short memory models}
Much of statistical inference
is concerned with {\it short memory} stationary models,
where the spectral density $f(\lambda)$ of the
model is bounded away from zero and infinity,
that is, there are constants $C_1$ and $C_2$ such that
\begin{equation*}
0< C_1 \le f(\la) \le C_2 <\f.
\end{equation*}

A typical d.t.\  short memory model example is the stationary
Autoregressive Moving Average (ARMA)$(p,q)$ process $X(t)$
defined to be a stationary solution of the difference equation:
\begin{equation*}
\psi_p(B)X(t)=\theta_q(B)\varepsilon(t), \q t\in\mathbb{Z},
\end{equation*}
where $\psi_p$ and $\theta_q$ are polynomials of degrees $p$ and $q$,
respectively, $B$ is the backshift operator defined by $BX(t)=X(t-1)$,
and $\{\vs(t), t\in\mathbb{Z}\}$ is a d.t.\  white noise,
that is, a sequence of zero-mean, uncorrelated random variables
with variance $\si^2$.
The covariance $r(k)$ of (ARMA)$(p,q)$ process is
exponentially bounded:
\[
|r(k)|\le Cr^{-k}, \q k=1,2,\ldots;\q0<C<\f; \,\, 0<r<1,
\]
and the spectral density $f(\la)$ is a rational function
(see, e.g., Brockwell and Davis \cite{BD}, Section 3.1):
%
\begin{equation}
\label{arma}
f(\la) = \frac{\si^2}{2\pi}\cd\frac
{|\theta_q(e^{-i\la})|^2}{|\psi_p(e^{-i\la})|^2}.
\end{equation}

A typical c.t.\  short-memory model example is the stationary
c.t.\  ARMA$(p,q)$ processes, denoted by CARMA$(p,q)$,
which is defined to be the solution of a $p$th order
stochastic differential equation with a suitable initial condition and driven by a
standard Brownian  motion and its derivatives
up to and including order $0\leq q < p$.
The spectral density function $f(\la)$ of the CARMA$(p,q)$ process is given by the following formula (see, e.g., Brockwell \cite{Br1}):
\begin{equation}
\label{CAS}
f(\la)=\frac{\si^2}{2\pi}\cd\frac{|\be(i\la)|^2}{|\al(i\la)|^2},
\end{equation}
where $\al(z)$ and $\be(z)$ are
polynomials of degrees $p$ and $q$, respectively.

Another important c.t.\  short memory model is the
{\it Ornstein-Uhlenbeck} process, which is a Gaussian stationary
process with covariance function $r(t)=\si^2e^{-\al |t|}$ ($t\in\mathbb{R}$),
and spectral density
\begin{equation}
\label{OU}
f(\la) = \frac{\si^2}{\pi}\cd\frac {\al^2}{\la^2+\al^2},\q \al>0,  \, \la\in\mathbb{R}.
\end{equation}

\ssn{Discrete-time long-memory and anti-persistent models}
Data in many fields of science (economics, finance, hydrology, etc.),
however, is well modeled by stationary processes whose spectral densities are
{\it unbounded} or {\it vanishing} at some fixed points  (see, e.g.,
Beran \cite{Ber}, Gu\'egan \cite{Gu1}, Palma \cite{Pal}, Taqqu \cite{Tq1}
and references therein).

A {\it long-memory} model is defined to be a
stationary process with {\it unbounded} spectral density,
and an {\it anti-persistent} model -- a stationary
process with {\it vanishing} (at some fixed points) spectral density.

In the discrete context, a basic long-memory model
is the Autoregressive Fractionally Integrated Moving Average
(ARFIMA)$(0,d,0))$ process $X(t)$ defined to be a stationary
solution of the difference equation
(see, e.g., Brockwell and Davis \cite{BD}, Section 13.2):
\begin{equation*}
(1-B)^dX(t)=\varepsilon(t), \q0<d<1/2,
\end{equation*}
where $B$ is the backshift operator and
$\varepsilon(t)$ is a d.t.\  white noise defined above.
The spectral density $f(\la)$ of $X(t)$ is given by
%
\begin{equation}
\label{MT0}
f(\la)=|1-e^{-i\la}|^{-2d}=(2\sin(\la/2))^{-2d},
\q0<\la\le\pi, \q0<d<1/2.
\end{equation}
Notice that $f(\la)\thicksim c\, |\la|^{-2d}$ as $\la\to0$, that is, $f(\la)$ blows up
at $\la=0$ like a power function, which is the typical behavior of
a long memory model.

A typical example of an {\it anti-persistent} model is the
ARFIMA$(0,d,0)$ process $X(t)$ with spectral density 
$f(\la)=|1-e^{-i\la}|^{-2d}$
with $d<0$, which vanishes at
$\la=0$.
Note that the condition $d<1/2$ ensures that
$\int_{-\pi}^\pi f(\la)d\la<\f$,
implying that the process $X(t)$ is well defined because
$E|X(t)|^2=\int_{-\pi}^\pi f(\la)d\la.$

Data can also occur in the form of a realization of a 'mixed'
short-long-intermediate-memory stationary process $X(t)$.
A well-known example of such a process, which appears in many applied
problems, is an ARFIMA$(p,d,q)$ process $X(t)$ defined to be a
stationary solution of the difference equation:
\begin{equation*}
\psi_p(B)(1-B)^dX(t)=\theta_q(B)\varepsilon(t), \q d<1/2,
\end{equation*}
where $B$ is the backshift operator, $\varepsilon(t)$ is a
d.t.\  white noise, and $\psi_p$ and $\theta_q$ are
polynomials of
degrees $p$ and $q$, respectively.
The spectral density $f_X(\la)$ of $X(t)$ is given by
%
\begin{equation}
\label{AA}
f_X(\la)=|1-e^{-i\la}|^{-2d}f(\la), \q d<1/2,
\end{equation}
where $f(\la)$ is the spectral density of an ARMA$(p,q)$ process,
given by (\ref{arma}).
\noindent
Observe that for $ 0<d<1/2$ the model $X(t)$ specified by (\ref{AA})
displays long-memory, for $d<0$ -- intermediate-memory,
and for $d=0$ -- short-memory.
For $d \ge1/2$ the function $f_X(\la)$ in (\ref{AA}) is not integrable,
and thus it cannot represent a spectral density of a stationary
process. Also, if $d \le-1$, then the series $X(t)$ is not
invertible in the sense that it cannot be used to recover a white noise
$\varepsilon(t)$ by passing $X(t)$ through a linear filter
(see, e.g., Bondon and Palma \cite{BP1}, and Brockweel and Davis \cite{BD}).

The ARFIMA$(p,d,q)$ processes, first introduced by Granger and Joyeux
\cite{GJ}, and Hosking \cite{Hos}, became very popular due to their
ability in
providing a good characterization of the long-run properties of many
economic and financial time series. They are also very useful for
modeling multivariate time series, since they are able to capture a
larger number of long term equilibrium relations among economic
variables than the traditional multivariate ARIMA models
(see, e.g., Gu\'egan \cite{Gu1}, and Henry and Zaffaroni \cite{HZ}
for a survey on this topic).

Another important long-memory model is the fractional Gaussian noise (fGn).
To define fGn we first consider the {\it fractional Brownian motion\/} (fBm)
\{$B_H(t), t\in\mathbb{R}\}$ with Hurst index $H$, $0<H<1$,
defined to be a centered Gaussian $H$-self-similar process having stationary
increments, that is, $B_H(t)$ satisfies the following conditions:

(a) $B_H(0)=0$, $\mathbb{E}[B_H(t)]=0$, $t\in\mathbb{R}$;

(b) $\{B_H(at), t\in\mathbb{R}\} \stackrel{d}
= \{a^HB_H(t), t\in\mathbb{R}\}$ for any $a>0$;

(c) $\{B_H(t+u)-B_H(u), t\in\mathbb{R}\} \stackrel{d}
=\{B_H(t), t\in\mathbb{R}\}$ for each fixed $u\in\mathbb{R}$;

(d) the covariance function is given by
\[
{\rm Cov}(B_H(s),B_H(t))
=\frac{\sigma_0^2}2\left[|t|^{2H}-|s|^{2H}-|t-s|^{2H}\right],
\]
where the symbol $\stackrel{d}=$ stands for equality of the
finite-dimensional distributions, and $\sigma_0^2={\rm Var} B_H(1)$.
Then the increment process
\[
\{X(k):= B_H(k+1)-B_H(k), k\in\mathbb{Z}\},
\]
called {\it fractional Gaussian noise\/} (fGn), is a d.t.\
centered Gaussian stationary process with covariance function
%
\begin{equation}
\label{MT1}
r(k)=\frac{\sigma_0^2}2\left[|k+1|^{2H}-|k|^{2H}-|k-1|^{2H}\right],
\q k\in\mathbb{Z},
\end{equation}
and spectral density function
%
\begin{equation}
\label{MT2}
f(\la)=c\, |1-e^{-i\la}|^{2}\sum_{k=-\f}^\f|\la+2\pi k|^{-(2H+1)},
\q-\pi\le\la\le\pi,
\end{equation}
where $c$ is a positive constant.

It follows from (\ref{MT2}) that $f(\la)\thicksim c\, |\la|^{1-2H}$
as $\la\to0$, that is, $f(\la)$ blows up if $H>1/2$
and tends to zero if $H<1/2$.
Also, comparing (\ref{MT0}) and (\ref{MT2}), we observe that,
up to a constant, the spectral density of fGn has the same
behavior at the origin as ARFIMA$(0,d,0)$ with $d=H-1/2.$

Thus, the fGn $\{X(k), k\in\mathbb{Z}\}$ has long-memory
if $1/2<H<1$ and is anti-percipient if $0<H<1/2$.
The variables $X(k)$, $k\in\mathbb{Z}$, are independent if $H=1/2$.
For more details we refer to Samorodnisky and Taqqu \cite{ST}
and Taqqu \cite{Tq1,Tq3}.

\ssn{Continuous-time long-memory and anti-persistent models}
In the continuous context, a basic process which has commonly been used
to model long-range dependence is the fractional Brownian motion (fBm)
$B_H$ with Hurst index $H$, defined above, which can be regarded as a
Gaussian process having a "spectral density":
%
\begin{equation}
\label{lr3}
f(\la)=c|\la|^{-(2H+1)}, \q c>0, \ \,\,
0<H<1, \  \,\, \la\in\mathbb{R}.
\end{equation}
The form (\ref{lr3}) can be understood in a generalized sense
(see Yaglom \cite{Ya}, Section 24, Flandrin \cite{Fl}, Solo \cite{So}),
since the fBm $B_H$ is a nonstationary process.

A proper stationary model in lieu of fBm is the
{\it fractional Riesz-Bessel motion} (fRBm),
introduced in Anh et al. \cite{AAR}, and defined as a
c.t.\  Gaussian process $X(t)$ with
spectral density
\begin{equation}
\label{rb1}
f(\la)=c\,|\la|^{-2\al}(1+\la^2)^{-\be},\q\la\in\mathbb{R},
\, 0<c<\f, \, 0<\al<1, \,\be>0.
\end{equation}
The exponent $\al$ determines the long-range dependence,
while the exponent $\be$
indicates the second-order intermittency of the process
(see, e.g., Anh et al. \cite{ALM} and Gao et al. \cite{GAHT}).

Notice that the process $X(t)$, specified by (\ref{rb1}),
is stationary if $0<\al<1/2$
and is non-stationary with stationary increments if $1/2\le\al<1.$
\noindent Observe also that the spectral density (\ref{rb1}) behaves
as $O(|\la|^{-2\al})$ as $|\la|\to0$ and as $O(|\la|^{-2(\al+\be)})$
as $|\la|\to\f$.
Thus, under the conditions $0<\al<1/2$, $\be>0$ and $\al+\be>1/2,$
the function $f(\la)$ in (\ref{rb1}) is well-defined for both
$|\la|\to0$ and $|\la|\to\f$
due to the presence of the component $(1+\la^2)^{-\be}$, $\be>0$,
which is the Fourier transform of the Bessel potential.

Comparing (\ref{lr3}) and (\ref{rb1}), we observe that the
spectral density of fBm is the limiting case as $\be\to0$
that of fRBm with Hurst index $H=\al-1/2.$

Another important c.t.\  long-memory model is the CARFIMA$(p,H,q)$ processes, which is defined to be the solution of a $p$th order
stochastic differential equation with a suitable initial condition and driven by a
fractional Brownian motion with Hurst parameter $H$ and its derivatives
up to and including order $0\leq q < p$.
The spectral density function $f(\la)$
of the CARFIMA$(p,H,q)$ process is given by the following formula (see, e.g., Tsai and Chan \cite{TC}):
\begin{equation}
\label{CAA}
f(\la)=\frac{\si^2}{2\pi}\G(2H+1)\sin(\pi H)|\la|^{1-2H}\frac{|\be(i\la)|^2}{|\al(i\la)|^2},
\end{equation}
where $\al(z)$ and $\be(z)$ are
polynomials of degrees $p$ and $q$, respectively.
Notice that for $H=1/2$, the spectral density function given by \eqref{CAA} becomes
that of the short-memory CARMA$(p,q)$ process, given by \eqref{CAS}.

\section{CLT for Toeplitz type quadratic functionals for Gaussian and linear processes}
\label{G-CLT}

In this section we present sufficient conditions for quadratic forms and functionals
of the form \eqref{MTc-1} to obey the CLT.
The processes considered will be d.t.\  and c.t.\  Gaussian
or linear processes with memory.
The matrix and the operator that characterize the quadratic form and
functional are Toeplitz.

As it was mentioned in the introduction, the limit distributions of the functionals in (\ref{MTc-1})
are completely determined by the spectral density $f(\la)$ (or covariance function $r(t)$)
and the generating function $g(\la)$  (or generating kernel $\widehat g(t)$).
Conditions that are in terms of the spectral density $f(\la)$ and the generating function
$g(\la)$ will be called frequency-domain conditions, while conditions that are in terms
of the covariance function $r(t)$ and the generating kernel $\widehat g(t)$ will be called
time-domain conditions.

\sn{Frequency domain conditions}
\label{fdc}

Let $\{X(u), \ u\in \mathbb{U}\}$ be a centered real-valued
Gaussian stationary process with spectral density
$f(\la)$, $\la\in \Lambda$ and covariance function
$r(t):=\widehat f(t)$, $t\in \mathbb{U}$, where $\mathbb{U}$
and $\Lambda$ are as in Section \ref{model}.
\n We are interested in the asymptotic distribution
(as $T\to\f$) of the following Toeplitz type quadratic
functionals of the process $X(u)$:
\beq
\label{c-1}
Q_T:= Q_T(f,g)=
 \left \{
 \begin{array}{ll}
\sum_{t=1}^T\sum_{s=1}^T\widehat g(t-s)X(t)X(s) &
\mbox{in the d.t.\  case},\\
\\
\int_0^T\int_0^T\widehat g(t-s)X(t)X(s)\,dt\,ds&
\mbox{in the c.t.\  case},
\end{array}
\right.
\eeq
where $\widehat g(t)$
is the Fourier transform of some real, even, integrable function $g(\la),$ $\la\in\Lambda$.
In the d.t.\  case the functions $f(\la)$ and $g(\la)$ are assumed to be $2\pi$-periodic
and periodically extended to $\mathbb{R}$.

\n
{\it Note.} We include the function $f$ in the notation $Q_T(f,g)$
to emphasize that the distribution of the quadratic form depends also on the spectral density $f$.
Let $Q_T$ be as in (\ref{c-1}). By $\widetilde Q_T$ we denote the standard
normalized quadratic functional:
\begin{equation}
\label{c-3}
\widetilde Q_T:= T^{-1/2}\,\left(Q_T-\E[Q_T]\right).
\end{equation}
As before, the notation
\begin{equation}
\label{c-4}
\widetilde Q_T\ConvD \eta \sim N(0,\sigma^2) \q {\rm as}\q T\to\f
\end{equation}
will mean that the distribution of the random variable
$\widetilde Q_T$ tends (as $T\to\f$) to the centered normal
distribution with variance $\sigma^2$.

Toeplitz matrices and operators arise naturally
in the theory of stationary processes, and serve as tools, to study many
topics of the spectral and statistical analysis of d.t.\  and c.t.\
stationary processes.

We first define the truncated Toeplitz matrices and operators, generated
by integrable real symmetric functions.

Let $\psi(\la)$ be an integrable real symmetric function defined on
$\Lambda=[-\pi, \pi]$. For  $T=1, 2,\ldots$,
the {\it $(T\times T)$-truncated Toeplitz matrix\/} generated by
$\psi(\la)$, denoted by $B_T(\psi)$, is defined by the following
equation (see, e.g., Ginovyan and Sahakyan \cite{GS2005}, and
Grenander and Szeg\H{o} \cite{GS}):
\beq
\label{IMT2-1}
B_T(\psi):=\|\widehat\psi(t-s)\|_{t,s=1,2\ldots,T}
=\left(
\begin{array}{llll}
\widehat \psi(0)&\widehat \psi(-1)&\cdots&\widehat \psi(1-T)\\
\widehat \psi(1)&\widehat \psi(0)&\cdots&\widehat \psi(2-T)\\
\cdots&\cdots&\cdots&\cdots\\
\widehat \psi(T-1)&\widehat \psi(T-2)&\cdots&\widehat\psi(0)
\end{array}
\right),
\eeq
where $\widehat\psi(t)=\int_\Lambda e^{i\la t}\,\psi(\la)\,d\la$
$(t\in \mathbb{Z})$ are the Fourier coefficients of $\psi$.

Given a real number $T>0$ and an integrable real symmetric function $\psi(\la)$
defined on $\mathbb{R}:=(-\f, \f)$, the {\it $T$-truncated Toeplitz operator\/}
(also called {\it Wiener-Hopf operator}) generated by $\psi(\la)$, denoted by $W_T(\psi)$,
is defined by the following equation (see, e.g., Ginovyan \cite{G1994},
 Ginovyan and Sahakyan \cite{GS2007},
Grenander and Szeg\H{o} \cite{GS}, Ibragimov \cite{I1963}, and Kac \cite{Kc}):
\beq
\label{IMT3-1}
[W_T(\psi)u](t)=\int_0^T\hat\psi(t-s)u(s)ds,
\q u(s)\in L^2[0,T],
\eeq
where $\widehat\psi(t)=\int_\mathbb{R} e^{i\la t}\,\psi(\la)\,d\la$
$(t\in \mathbb{R})$ is the Fourier transform of $\psi(\la)$.

Let $A_T(f)$ and $A_T(g)$ denote either the $T$-truncated Toeplitz
operators (in the c.t.\  case), or the $T\times T$ Toeplitz matrices
(in the d.t.\  case) generated by the functions $f$ and $g$, respectively.
Observe that $A_T(f)$ is the covariance matrix (or operator) of the process
$\{X(u), \ u\in \mathbb{U}\}$.

We assume below that $f, g \in L^1(\Lambda)$, and with no loss of generality,
that $g\ge 0$. Also, we set
\beq
\label{c-6}
\si^2_0: = 16\pi^3\int_\Lambda f^2(\la)g^2(\la)\,d\la.
\eeq
As usual $\Lambda=[-\pi, \pi]$ in the d.t.\  case and
$\Lambda=\mathbb{R}$ in the c.t.\  case.

The theorems that follow contain sufficient conditions expressed in terms
of $f(\la)$ and $g(\la)$ to ensure central limit theorems
for standard normalized quadratic functionals $\widetilde Q_T$ both
for d.t.\  and c.t.\  Gaussian processes.
Some of the assumptions imposed on $f$ allow for long-range dependence ($f (0) = \f$),
others for discontinuities at other frequencies. Sometimes the good behavior
of one function, say $g$ can compensate for the bad behavior of $f$ and vice versa.

\begin{thm}
\label{cth1}
Let $f$, $g$,  $A_T(f)$, $A_T(g)$, and $\widetilde Q_T$ be as above.
Each of the following conditions is sufficient for functional $Q_T$ to obey the CLT, that is,
\beq
\label{c-7}
\widetilde Q_T\ConvD \eta \sim N(0,\sigma_0^2) \q {\rm as}\q T\to\f,
\eeq
with $\si^2_0$ given by (\ref{c-6}).
\begin{itemize}
\item[{\bf (A)}]
$f\cdot g\in L^1(\Lambda)\cap L^2(\Lambda)$ and
\beq
\label{c-8}
\chi_2(\widetilde Q_T):=
\frac2T\tr\bigl[A_T(f)A_T(g)\bigr]^2 \longrightarrow
\sigma_0^2,
\eeq
where $\tr[A]$ stands for the trace of the operator (or the matrix) $A$.
\item[{\bf (B)}]
The function
\beq
\label{c-9}
\varphi({\bf u}):=\varphi(u_1, u_2,u_3)=\int_\Lambda
f(\la)g(\la-u_1)f(\la-u_2)g(\la-u_3)\,d\la
\eeq
belongs to $L^2(\Lambda^3)$ and is continuous at ${\bf 0}=(0,0,0).$

\item[{\bf (C)}]
$f \in L^1(\Lambda)\cap L^p(\Lambda)$ $(p\ge2)$
and  $g\in L^1(\Lambda)\cap L^q(\Lambda)$ $(q\ge2)$
with $1/p+1/q\le1/2.$

\item[{\bf (D)}]
$f\in L^1(\Lambda)\cap L^2(\Lambda)$,
\,$g\in L^1(\Lambda)\cap L^2(\Lambda)$,
$fg\in L^2(\Lambda)$ and
\beq
\label{c-10}
\nonumber
\int_\Lambda  f^2(\la)g^2(\la-\mu)\,d\la \longrightarrow
 \int_\Lambda f^2(\la)g^2(\la)\,d\la \quad {\rm as} \quad \mu\to0.
\eeq
\end{itemize}
\end{thm}

\begin{rem}
{\rm
Observe that assertion (A) implies assertions (B) -- (D), and
assertion (B) implies assertions (C) and (D) (see Giraitis and Surgailis \cite{GSu},
Ginovyan and Sahakyan \cite{GS2005, GS2007}).
For the d.t.\  case:
assertions (A) and (D) were proved in Giraitis and Surgailis \cite{GSu}
(see also Giraitis et al. \cite{GKSu});
assertions (B)  was proved in Ginovyan and Sahakyan \cite{GS2005};
assertion (C) for $p=q=\f$ was first established by Grenander and Szeg\H{o}
(\cite{GS}, Section 11.7), while the case $p=2$, $q=\f$ was proved by Ibragimov
\cite{I1963} and Rosenblatt \cite{R2}, in the general d.t.\  case
assertion (C) was proved by Avram \cite{A}.
For the c.t.\  case, assertions (A) -- (D) were proved in Ginovyan \cite{G1994}
and in Ginovyan and Sahakyan \cite{GS2007}.}
\end{rem}

The following theorem provides conditions for the CLT to hold when either
$f$ or $g$ is a sufficiently smooth function.
Given a number $\al$ ($0<\al<1$). Denote by $\rm{Lip}_\al(\Lambda)$ the class of Lipschitz
functions on $\Lambda$. By definition, $\phi\in \rm{Lip}_\al(\Lambda)$ if there exists
a constant $C<\f$ such that
\[
|\phi(x)-\phi(y)|\leq C|x-y|^\al \q \text{for all} \q x,y\in \Lambda.
\]
\begin{thm}\label{GT-G}
Let either $f\in \rm{Lip}_\al(\Lambda)$ or $g\in \rm{Lip}_\al(\Lambda)$
with $\al>1/2$, and let $fg\in L^2(\Lambda)$.
Then $\widetilde Q_T\ConvD \eta \sim N(0,\sigma^2_0)$ as $T\to\f$ with $\si^2_0$ given by (\ref{c-6}).
\end{thm}
Theorem \ref{GT-G} for the d.t.\  case was proved in Giraitis and Taqqu \cite{GT1997}.
For the c.t.\  case it can be proved similarly.

To state the next theorem, we need to introduce a class of slowly varying functions
at zero. Recall that a function $u(\la)$, $\la\in\mathbb{R}$,
is called slowly varying at zero if it is non-negative and for any $t>0$
\[
\lim_{\la\rightarrow 0} \frac{u(t\la)}{u(\la)}\rightarrow 1.
\]
\n
Denote by $SV_0(\Lambda)$ the class of slowly varying functions at zero
$u(\la)$, $\la\in\Lambda$, satisfying the following conditions:
for some $a>0$, $u(\la)$ is bounded on $[-a,a]$, $\lim_{\la\to0}u(\la)=0,$ \
$u(\la)=u(-\la)$ and $0<u(\la)<u(\mu)$\ for\ $0<\la<\mu<a$.
An example of a function belonging to $SV_0(\Lambda)$
is $u(\la)= \left|\ln |\la|\right|^{-\gamma}$ with $\gamma>0$ and $a=1$.

\begin{thm}\label{cth2}
Assume that the functions $f$ and $g$ are integrable on $\mathbb{R}$
and bounded outside any neighborhood of the origin, and satisfy for some $a>0$
the following conditions:
\beq \label{cm-0}
f(\lambda)\le |\lambda|^{-\alpha}L_1(\lambda),
\q
|g(\lambda)|\le |\lambda|^{-\beta}L_2(\lambda),\q \lambda\in [-a,a],
\eeq
for some $\alpha<1, \ \beta<1$ with $\alpha+\beta\le1/2$, where $L_1(x)$ and
$L_2(x)$ are slowly varying functions at zero satisfying
\bea
\label{cm-00}
L_i\in SV_0(\mathbb{R}), \ \ \lambda^{-(\alpha+\beta)}L_i(\lambda)\in L^2[-a,a], \ i=1,2.
\eea
Then $\widetilde Q_T\ConvD \eta \sim N(0,\sigma^2_0)$ as $T\to\f$ with $\si^2_0$ given by (\ref{c-6}).
\end{thm}
\begin{rem}
{\rm
For the d.t.\  case the result of Theorem \ref{cth2} under the condition $\al+\be<1/2$ was first
obtained by Fox and Taqqu \cite{FT1}. For the general case, including the critical
value $\al+\be=1/2$ it was proved in Ginovyan and Sahakyan \cite{GS2005}.
For the c.t.\  case the result was proved in Ginovyan and Sahakyan \cite{GS2007}.}
\end{rem}

\begin{rem}\label{Rem:B follow A}
{\rm   The conditions $\al<1$ and $\be<1$ in Theorem \ref{cth2} ensure that the
Fourier transforms of $f$ and $g$ are well defined. Observe that when $\al>0$
the process $X(t)$ may exhibit long-range dependence.
We also allow here  $\alpha+\beta$ to assume the critical value 1/2.
The assumptions $f\cdot g\in L^1(\mathbb{R})$, $f,g\in L^\infty(\mathbb{R}\setminus [-a,a])$
and (\ref{cm-00}) imply that $f\cdot g \in L^2(\mathbb{R})$, so that the variance
$\sigma^2_0$ in (\ref{c-6}) is finite.}
\end{rem}
\begin{rem}
{\rm In Theorem \ref{cth2}, the assumption that $L_1(x)$ and $L_2(x)$  belong to $SV_0(\mathbb{R})$
instead of merely being slowly varying at zero is done in order to deal with the critical case $\alpha+\beta=1/2$.
Suppose that we are away from this critical case, namely, $f(x)=|x|^{-\alpha}l_1(x)$ and $g(x)=|x|^{-\beta}l_2(x)$, where $\alpha+\beta< 1/2$, and $l_1(x)$ and $l_2(x)$ are slowly
varying at zero functions. Assume also that $f(x)$ and $g(x)$ are integrable and bounded on
$(-\infty,-a)\cup(a,+\infty)$ for any $a>0$. We claim that Theorem \ref{cth2} applies.
Indeed, choose $\alpha'>\alpha$, $\beta'>\beta$ with $\alpha'+\beta'<1/2$.
Write
$f(x)=|x|^{-\alpha'} |x|^{\delta}l_1(x)$, where $\delta=\alpha'-\alpha>0$.
Since $l_1(x)$ is slowly varying, when $|x|$ is small enough,  for some
$\epsilon\in (0,\delta)$ we have $
|x|^{\delta}l_1(x)\le |x|^{\delta-\epsilon}$. Then  one can  bound
$|x|^{\delta-\epsilon}$ by $c\left|\ln |x|\right|^{-1}\in SV_0(\mathbb{R})$
for small $|x|<1$. Hence one  has  when  $|x|<1$ is small enough,
$$
f(x)\le |x|^{-\alpha'} \left(c  \left|\ln |x|\right|^{-1} \right).
$$
Similarly, when $|x|<1$ is small enough, one has
$$
g(x)\le |x|^{-\beta'} \left(c  \left|\ln |x|\right|^{-1} \right).
$$
All the assumptions in Theorem \ref{cth2} are now readily checked with
$\alpha$ and $\beta$ replaced by $\alpha'$ and $\beta'$, respectively.}
\end{rem}

\begin{rem}
{\rm The functions
\beq
\label{c-13}
f(\la)= |\la|^{-\al}|\ln|\la||^{-\g} \q \mbox{and} \q
g(\la)= |\la|^{-\be}|\ln|\la||^{-\g},
\eeq
where $\al<1,$\ $\be<1,$\ $\al+\be\le1/2$ and $\g>1/2$, provide examples of
spectral density $f(\la)$ and generating function $g(\la)$
satisfying the conditions of Theorem \ref{cth2}
(see Ginovyan and Sahakyan \cite{GS2005}).}
\end{rem}
\begin{rem}
{\rm The slowly varying functions $L_1$ and $L_2$ in (\ref{cm-0})
are important because they provide a great flexibility
in the choice of spectral density $f$ and generating function $g$.
Theorem \ref{cth2} shows that in the critical case $\al+\be=1/2$
the limit distribution of the standard normalized quadratic form
$\widetilde Q_T$ given by (\ref{c-3}) is Gaussian and depends on the
slowly varying factors $L_1$ and $L_2$ through $f$ and $g$.}
\end{rem}

\begin{rem}
{\rm The functions $f(\la)$ and $g(\la)$ in Theorem \ref{cth2}
have singularities at the point $\la=0$, and are bounded in
any neighborhood of this point. It can be shown that the choice
of the point $\la=0$ is not essential, and instead
it can be taken to be any point $\la_0\in[-a, a]$.
Using the properties of the products of Toeplitz matrices
and operators it can be shown that Theorem \ref{cth2}
remains valid when $f(\la)$ and $g(\la)$ have
singularities of the form (\ref{c-13}) at the
same finite number of points of the segment $[-a, a]$
(see Ginovyan and Sahakyan \cite{GS2005}).}
\end{rem}

\begin{rem}
{\rm Assertion (A) of Theorem \ref{cth1} implies Theorem \ref{cth2}.
On the other hand, for functions $f(\la)=\la^{-3/4}$ and $g(\la)=\la^{3/4}$
satisfying the conditions of Theorem \ref{cth2}, the function
$\I(u_1, u_2,u_3)$ in \eqref{c-9} is not defined for \,$u_2=0$, $u_1\neq 0$, $u_3\neq 0$,
showing that assertion (B) of Theorem \ref{cth1} generally does not
imply Theorem \ref{cth2} (see Ginovyan and Sahakyan \cite{GS2005}).}
\end{rem}

Giraitis and Surgailis \cite{GSu} proved that assertion (A) of
Theorem \ref{cth1} remains valid for d.t.\  linear processes.
More precisely, they proved the following theorem.
\begin{thm}
\label{GSth1}
Let $\{X(u), \,u \in \mathbb{Z}\}$ be a d.t.\  stationary linear process of the
form \eqref{dlp} with spectral density $f$. Let $Q_T$ be a quadratic form generated by a
function $g$ given by \eqref{c-1}, and let $B_T(f)$ and $B_T(g)$ be the $T\times T$
Toeplitz matrices generated by the functions $f$ and $g$, respectively (see \eqref{IMT2-1}).
Assume that
\beq
\label{gsc-8}
\frac1T\tr\bigl[B_T(f)B_T(g)\bigr]^2 \longrightarrow 8\pi^3 \inl
f^2(\lambda) g^2(\lambda) d\lambda<\f,
\eeq
Then the CLT holds for $Q_T$, that is,
\beq
\label{c-77}
\widetilde Q_T:= T^{-1/2}\,\left(Q_T-\E[Q_T]\right)\ConvD \eta \sim N(0,\sigma^2)
\q {\rm as}\q T\to\f,
\eeq
where
\begin{equation}\label{var}
\sigma^2=16\pi^3 \inl f^2(\lambda) g^2(\lambda) d\lambda
+\kappa_4 \left[2\pi\inl f(\lambda) g(\lambda) d\lambda\right]^2,
\end{equation}
and $\kappa_4$ is the fourth cumulant of $\xi(0)$.
\end{thm}
The next result, which was proved in Giraitis and Taqqu \cite{GT1997},
extends Theorem \ref{GT-G} to the case of d.t.\  linear processes.
\begin{thm}\label{GT-L}
Let $\{X(u), \,u \in \mathbb{Z}\}$ be a d.t.\  stationary linear process of the
form \eqref{dlp} with spectral density $f$ and $\E[\xi^2(0)]=1$, and let $Q_T$ be a quadratic
form generated by a function $g$ given by \eqref{c-1}.
Suppose either $f\in \rm{Lip}_\al(\Lambda)$ or $g\in \rm{Lip}_\al(\Lambda)$
with $\al>1/2$, and $fg\in L^2(\Lambda)$.
Then $\widetilde Q_T\ConvD \eta \sim N(0,\sigma^2)$ as $T\to\f$ with $\si^2$
given by (\ref{var}).
\end{thm}

Giraitis and Surgailis \cite{GSu} conjectured that the relation \eqref{gsc-8},
and hence, the CLT for $Q_T$, holds under the single condition that the integral
on the right hand side of \eqref{gsc-8} is finite.  Ginovyan \cite{G1993} showed
that the finiteness of this integral does not guarantee convergence in \eqref{gsc-8},
and conjectured that positiveness and finiteness of the integral
in \eqref{gsc-8} must be sufficient for $Q_T$ to obey the CLT.

The next proposition shows that the condition of positiveness
and finiteness of the integral in \eqref{gsc-8} also
is not sufficient for $Q_T$ to obey the CLT (see Ginovyan and Sahakyan \cite{GS2005}).

\begin{pp}
\label{m1993}
There exist a spectral density $f(\la)$ and a generating function
$g(\la)$ such that
\beq
\label{c-14}
0<\inl f^2(\la)\,g^2(\la)\,d\la<\f
\eeq
and
\beq
\label{c-15}
\lim_{T\to \infty}\sup\chi_2(\widetilde Q_T)
=\lim_{n\to \infty}\sup\frac2T\tr\left(B_T(f)B_T(g)\right)^2=\infty,
\eeq
that is, the condition (\ref{c-14}) does not guarantee convergence
in (\ref{gsc-8}), and hence is not sufficient for $Q_T$ to obey the CLT.
\end{pp}
To construct functions $f(\la)$ and $g(\la)$ satisfying
(\ref{c-14}) and (\ref{c-15}), for a fixed $p\ge2$
we choose a number $q>1$ to satisfy $1/p+1/q>1$,
and for such $p$ and $q$ we consider the functions
$f_0(\la)$ and $g_0(\la)$ defined by
\beq
\label{c-16}
f_0(\la)=
 \left \{
 \begin{array}{ll}
\left({2^s}/{s^2}\right)^{1/p}, & \mbox{if $2^{-s-1}\le
\la \le 2^{-s},\,s=2m$}\\
0,&  \mbox{if $2^{-s-1}\le\la \le 2^{-s},\,s=2m+1$},
\end{array}
\right.
\eeq
\beq
\label{c-17}
g_0(\la)=
 \left \{
\begin{array}{ll}
 \left({2^s}/{s^2}\right)^{1/q}, & \mbox{if $2^{-s-1}\le
\la \le 2^{-s},\,s=2m+1$}\\
0,& \mbox{if $2^{-s-1}\le\la \le
2^{-s},\,s=2m$},\end{array} \right.
\eeq
where $m$ is a positive integer. For an
arbitrary finite positive constant $C$ we set
$g_\pm(\la)=g_0(\la)\pm C$.
Then the functions $f=f_0$ and $g=g_+$ or $g=g_-$
satisfy (\ref{c-14}) and (\ref{c-15})
(for details we refer to Ginovyan \cite{G1993}, and Ginovyan and Sahakyan \cite{GS2005}).
Consequently, for these functions the quadratic form $Q_T$ does not obey the CLT,
and it is of interest to describe the limiting non-Gaussian
distribution of $Q_T$ in this special case.

\sn{Time domain conditions}
\label{time-c}

In this subsection we present time-domain sufficient conditions for a quadratic form
$Q_T$ of the form \eqref{c-1} to obey the CLT. That is, the conditions are in terms
of the covariance function $r(t)$ and the generating kernel $\widehat g(t)$.
The processes considered here will be d.t.\  Gaussian or linear
processes with memory.
Observe that the time-domain conditions stated below are, in general, not equivalent
to the frequency-domain conditions stated in Theorems \ref{cth1} and \ref{cth2}.
The method of proof for establishing the CLT is also different under the time-domain
conditions. In this case, diagrams are used and the method of moments is applied.
This is because, in the frequency-domain, one can use an approximation technique,
essentially replacing the possibly unbounded spectral density by a bounded one,
which allowed us to approximate the bivariate quadratic forms by univariate sums
of $m$-dependent random variables. Such an approximation technique, however,
does not work with time-domain conditions, and one has to deal directly with the
covariances, which decrease slowly.

The next results contain sufficient conditions in terms
of the covariance function $r(t)$ and the generating kernel $\widehat g(t)$
ensuring central limit theorems for standard normalized quadratic form
$\widetilde Q_T$ for the d.t.\  Gaussian processes
(see Fox and Taqqu \cite{FT}, and Giraitis and Taqqu \cite{GT1998}).
\begin{thm}
\label{GT-1}
If the covariance function $r(t)$ and the generating kernel $\widehat g(t)$
satisfy the following condition:
$$r \in l^p \, (p\ge1) \q\text{and}\q
\widehat g\in l^q \, (q\ge1) \q\text{with}\q 1/p+1/q\ge3/2,$$
then the CLT holds for the quadratic form $Q_T$ with limiting variance
$\si^2_0$ given by (\ref{c-6}).
\end{thm}

\begin{rem}
{\rm In fact, in Giraitis and Taqqu \cite{GT1998} (see also Giraitis and Taqqu \cite{GT1999})
was proved more general result
stating that Theorem \ref{GT-1} is true for quadratic forms of the form:
\begin{equation}
\label{nc-1}
Q_{T,m,n}=\sum_{t=1}^T\sum_{s=1}^T\widehat g(t-s)P_{m,n}(X(t),X(s)), \q T\in\mathbb{N},
\end{equation}
where $P_{m,n}(X(t),X(s))$ is a bivariate Appell polynomial (Wick power) of the linear
variables $X(t)$ and $X(s)$, $m, n \ge 0$, $m + n \ge 1$, and $X(t)$ is a d.t.\  linear process
of the form \eqref{dlp}. Also, observe that $Q_{T,1,1}=Q_{T}$.}
\end{rem}

The next theorem, which was proved in in Giraitis and Taqqu \cite{GT1998},
shows that under some rather restrictive conditions on the
covariance function $r(t)$ and the generating kernel $\widehat g(t)$,
the long-range dependence of the process $X(t)$ can be compensated by
the fast decay of the generating kernel $\widehat g(t)$ in such a way that the CLT for $Q_T$ holds.
These conditions ensure, in fact, that the sufficient conditions in the
frequency domain, provided in Theorem \ref {cth2}, are satisfied.
The theorem involves quasi-monotone sequences: a sequence $\{b(t), t\in\mathbb{Z}\}$
is said to be {\it quasi-monotonically convergent} to $0$ if $b(t)\to 0$ and
$b(t + 1) \leq b(t)(1 + c/t)$ as $t \to 0$ for some $c > 0$.
The sequence $b(t)$ has bounded variation if $\sum_{t=1}^\f|b(t+1)-b(t)| <\f$.

\begin{thm}
\label{GT-2}
Assume that
$r(t)= |t|^{-\g_1}L_1(|t|),$ $\widehat g(t)=|t|^{-\g_2}L_2(|t|)$ with $0<\g_1, \g_2<3$
and $\g_1 + \g_2>3/2$.
Suppose in addition that both sequences $\{r(t)\}$ and $\{{\widehat g(t)}\}$
have bounded variation and are quasi-monotonically convergent to $0$;
if $1 < \g_1 < 3$, $r(t)$ has the same sign for large $t$ and satisfies
$\sum_{t\in\mathbb{Z}}r(t) = 0$; if $1 < \g_2 < 3$, $\widehat g(t)$ has
the same sign for large $t$ and satisfies $\sum_{t\in\mathbb{Z}}\widehat g(t) = 0$.
Then the CLT holds for the quadratic form $Q_T$. The limiting variance
is expressed by (\ref{c-6}) provided that $0<\g_1, \g_2<1$.
\end{thm}

\sn{Operator conditions}
In this subsection we assume that the model $X(t)$ is a d.t.\
Gaussian process defined on a probability space $(\Omega,\mathcal{F}, P)$,
and we are interested in asymptotic normality of the quadratic form
$Q_T:=Q_T(f, g)$ given by \eqref{c-1} (d.t.\  case).
Recall the notation (see formula \eqref{c-3}
\beq
\label{c-5s}
\widetilde Q_T(f,g):= T^{-1/2}\,\left(Q_T(f,g)-\E[Q_T(f,g)]\right).
\eeq
and (see formula \eqref{c-6}
\beq
\label{c-6s}
\si^2(f,g): = 16\pi^3\inl f^2(\la)g^2(\la)\,d\la.
\eeq
We denote by $L^2(dP)$ the $L^2$-space constructed from a probability
measure $P$ such that $X(t)\in L^2(dP)$, $t\in\mathbb{Z}$.

Solev and Gerville-Reache \cite{SG2006} observed that for a fixed spectral
density $f$ and a number $T\in \mathbb{N}$ the quadratic form $\widetilde Q_T(f, g)$
in \eqref{c-5s} can be regarded as the value of a linear operator
$\mathbf{\widetilde Q}_T: g \mapsto \widetilde Q_T(f, g)$.
It turns out that in order to study the asymptotic normality of the quadratic form
$\widetilde Q_T(f, g)$, it is enough to understand for what sets $B_\vs$
(possibly, depending on $f$) we have
\begin{equation}
\label{sol-1}
\limsup_{T\to\f}\sup_{g\in B_\vs} ||\widetilde Q_T(f,g)||_{L^2(dP)}\leq \vs.
\end{equation}
It can be shown that the function $||\widetilde Q_T(f,g)||_{L^2(dP)}$ 
possesses the following symmetry property regarded as a function of $(f, g)$:
if $g$ is a nonnegative function, then
\begin{equation}
\label{sol-2}
||\widetilde Q_T(f,g)||^2_{L^2(dP)}=||\widetilde Q_T(g,f)||^2_{L^2(dP)}.
\end{equation}
The next result plays a key role in the proof of the asymptotic normality
of $\widetilde Q_T(f,g)$ (see Solev and Gerville-Reache \cite{SG2006}).

\begin{thm}
\label{st1}
Let $\mathcal{G}_0$ and $\mathcal{G}$ ($\mathcal{G}_0\subset\mathcal{G}$)
be linear subsets in the space $L^1[-\pi,\pi]$. Assume that the following
conditions hold:
\begin{itemize}
\item[(a)] for every $g_0\in \mathcal{G}_0$ the limit
$$
\lim_{T\to\f}||\widetilde Q_T(f,g_0)||^2_{L^2(dP)}=\si^2(f,g_0)<\f.
$$
exists, and
$$
\widetilde Q_T(f,g_0)\ConvD \eta_0 \sim N(0,\sigma^2(f,g_0)) \q {\rm as}\q T\to\f;
$$
\item[(b)] for every $g\in \mathcal{G}$ and a number $\vs>0$ there exists
a function $g_0\in \mathcal{G}_0$ such that
$$
\limsup_{T\to\f} ||\widetilde Q_T(f,g-g_0||_{L^2(dP)}\leq \vs,
\q ||g-g_0)||_{f^2}\leq \vs,
$$
where $||\cdot||_{L^2(dP)}$ and $||\cdot||_{f^2}$ are the norms in the spaces
$L^2(dP)$ and $L^2(f^2, [-\pi,\pi])$, respectively. Then $\si^2(f,g)<\f$, the limit
$$
\lim_{T\to\f}||\widetilde Q_T(f,g)||^2_{L^2(dP)}=\si^2(f,g)
$$
exists, and
\begin{equation}
\label{sol-3}
\widetilde Q_T(f,g)\ConvD \eta \sim N(0,\sigma^2(f,g)) \q {\rm as}\q T\to\f.
\end{equation}
\end{itemize}
\end{thm}
Thus, if the results holds for every $g_0\in \mathcal{G}_0$ (part (a)),
then it holds for every $g\in \mathcal{G}$ close to $g_0$ (part (b)).

Given a function $\phi\in L^1[-\pi,\pi]$ and an interval $I\subset[-\pi,\pi]$
of length $|I|$, we set
$$
\{\phi\}_I:=\frac1{|I|}\int_I\phi(x)dx,
$$
and define
\begin{equation}
\label{sol-4}
\la(\phi):=\sup_I\{\phi\}_I\times\{1/\phi\}_I,
\end{equation}
where the supremum is taken over all intervals $I\subset[-\pi,\pi]$.
The condition $\la(\phi)<\f$ implies a type of smoothness
of the function $\psi = \log|\phi|$ in the neighborhoods of points where
it turns to infinity (for details, see the paper Hunt et al. \cite{HMW}).

The following result was obtained by Solev and Gerville-Reache \cite{SG2006}.
\begin{thm}
\label{st2} Assume that the following conditions hold:
\begin{itemize}
\item[(a)] $f\leq f_*$, $|g|\leq g_*$ and $\la(f_*), \la(g_*)<\f$,
\item[(b)] $f_*g_*\in L^2[-\pi,\pi]$.
\end{itemize}
Then the quadratic form $Q_T(f,g)$ obeys the CLT, that is, the relation \eqref{sol-3} is satisfied.
\end{thm}

\section{Functional limit theorems for Gaussian and linear models}
\label{FLTGL}

In this section, we establish weak convergence in $C[0,1]$ for c.t.\  case
and in $D[0,1]$ for d.t.\  case of
normalized stochastic processes, generated by Toeplitz type quadratic
functionals of a Gaussian stationary process,
exhibiting long-range dependence.
Specifically, we are interested in describing the limiting behavior (as $T\to\f$)
of the following process $Q_T(\tau)$, generated by Toeplitz type quadratic functionals
of a Gaussian stationary process $\{X(u), \ u\in\mathbb{U}\}$ with spectral density $f$:

\beq
\label{eq:quad proc}
Q_T(\tau)=
 \left \{
 \begin{array}{ll}
\sum_{t=1}^{[T\tau]}\sum_{s=1}^{[T\tau]}\widehat g(t-s)X(t)X(s) &
\mbox{in the d.t.\  case},\\
\\
\int_0^{T\tau}\int_0^{T\tau}\widehat g(t-s)X(t)X(s)\,dt\,ds &
\mbox{in the c.t.\  case},
\end{array}
\right.
\eeq
where $\tau\in[0,1]$, $\widehat g(t)$, $t\in\mathbb{U}$ is the Fourier transform
of some integrable even function $g(\la)$, $\la\in\Lambda$, and $[\,\cdot\,]$
stands for the greatest integer.

The limit of the process (\ref{eq:quad proc}) again is completely
determined by the spectral density $f(\la)$ (or covariance function $r(t)$)
and the generating function $g(\la)$ (or generating kernel $\widehat g(t)$),
and depending on their properties, the limiting process can be Gaussian or not.
We say that we have a {\it functional non-central limit theorem} (FNCLT) for the process
$Q_T(\tau)$ if the limit is non-Gaussian or, if Gaussian, it is not Brownian motion.
We say that a {\it functional central limit theorem} (FCLT) for $Q_T(\tau)$ holds if we have
the following weak convergence in $C[0,1]$ in the c.t.\  case
(and in $D[0,1]$ in the d.t.\  case):
\beq
\label{BM}
 \widetilde Q_T(\tau) \Longrightarrow \sigma B(\tau),\quad \tau\in[0,1],
\eeq
where $\sigma>0$, $B(\tau)$ is a standard Brownian motion, and $\widetilde Q_T(\tau)$
is the standard normalized version of the process $Q_T(\tau)$ in (\ref{eq:quad proc}):
\begin{equation}\label{eq:tilde Q_T(t)}
 \widetilde Q_T(\tau) :=T^{-1/2}\left(Q_T(\tau)-\E[Q_T(\tau)]\right),
 \quad \tau\in[0,1].
\end{equation}

\subsection{Functional CLT for $Q_T(\tau)$}
We first state functional central limit theorems, when the limit of the
normalized process \eqref{eq:tilde Q_T(t)} is Brownian motion.

In the d.t.\  case the next result was obtained in Giraitis and Taqqu \cite{GT2001}
(see also Fox and Taqqu \cite{FT2} and Giraitis and Taqqu \cite{GT1997}).
\begin{thm}
\label{GT-1F}
If the covariance function $r(t)$ and the generating kernel $\widehat g(t)$
satisfy the following condition: $r \in l^p$ $(p\ge1)$ and
$\widehat g\in l^q$ $(q\ge1)$ with $1/p+1/q\ge3/2,$
then for the process $Q_T(\tau)$ the FCLT holds, that is, the convergence \eqref{BM}
holds in $D[0,1]$ with limiting variance $\si^2$ given by
\beq
\label{BMsi}
\si^2:=\sum_{t,s,v\in\mathbb{Z}}\widehat g(t)\widehat g(s)
\Cov \left(P(X(v),X(v+t)),P(X(0),X(s)\right),
\eeq
where $P(X(t),X(s)):=X(t)X(s)-\E[X(t)X(s)]$.
\end{thm}

Now we consider the c.t.\  case.
The next result, which is an extension of Theorem \ref{cth1}(A),
involves the convergence of finite-dimensional distributions of the process
$\widetilde Q_T(\tau)$ to those of a standard Brownian motion
(see Bai et al. \cite{BGT1}).
\begin{thm}\label{Thm:CLT-G}
Assume that the spectral density $f(\la)$ and the generating function $g(\la)$
satisfy the following conditions:
\begin{equation}
\label{eq:CLT condition 1}
f\cdot g\in L^1(\mathbb{R})\cap L^2(\mathbb{R})
\end{equation}
and
\begin{equation}\label{eq:CLT condition 2}
\E [\widetilde{Q}_T^2(1)] \rightarrow 16\pi^3\int_\mathbb{R} f^2(x) g^2(x)dx
~\text{ as }T\rightarrow\infty.
\end{equation}
Then we have the following convergence of the finite-dimensional distributions
\[
 \widetilde Q_T(\tau) \ConvFDD \sigma_0 B(\tau),
\]
where $\widetilde Q_T(\tau)$ is as in (\ref{eq:tilde Q_T(t)}),
$B(\tau)$ is a standard Brownian motion, and (see \eqref{c-6}):
\beq\label{eq:sigma}
\si^2_0: = 16\pi^3\int_\mathbb{R} f^2(x) g^2(x)dx.
\eeq
\end{thm}

To extend the convergence of finite-dimensional distributions in Theorem
\ref{Thm:CLT-G} to the weak convergence in the space $C[0,1]$,
we impose an additional condition on the underlying Gaussian process $X(t)$
and on the generating function $g$.
It is convenient to impose this condition in the time domain, that is,
on the covariance function $r:=\hat{f}$ and the generating kernel $a:=\hat{g}$.
The following condition is an analog of the assumption in Theorem \ref{GT-1F}
\begin{equation}\label{eq:CLT tightness cond}
r(\cdot) \in L^p(\mathbb{R}), \quad a(\cdot) \in L^q(\mathbb{R})
\quad\text{for some} \quad p,q\ge 1, ~ \frac{1}{p}+\frac{1}{q}\ge \frac{3}{2}.
\end{equation}
\begin{rem}
{\rm In fact under (\ref{eq:CLT condition 1}), the condition (\ref{eq:CLT tightness cond})
is  sufficient for the convergence in (\ref{eq:CLT condition 2}).
Indeed, let $\bar{p}=p/(p-1)$ and $\bar{q}=q/(q-1)$ be the H\"{o}lder
conjugates of $p$ and $q$, respectively. Since $1\le p,q\le 2$,
one has by the Hausdorff-Young inequality  that
$\|f\|_{\bar{p}}\le c_p\|r\|_p,~ \|g\|_{\bar{q}}\le c_q \|a\|_q, $
and hence, in view if (\ref{eq:CLT tightness cond}) we have
\beq \label{young}
f(\cdot) \in L^{\bar{p}},\quad g(\cdot) \in L^{\bar{q}}, \quad \frac{1}{\bar{p}}+\frac{1}{\bar{q}}=2-\frac{1}{p}-\frac{1}{q}\le 1/2.
\eeq
Then the convergence in (\ref{eq:CLT condition 2}) follows from Theorem \ref{cth1}(C)
(see the proof of Theorem 3 of Ginovyan and Sahakyan \cite{GS2007} in the
c.t.\  case, and Giraitis and Surgailis \cite{GSu} in the d.t.\  case)}.
\end{rem}
\begin{rem}
{\rm Observe that condition (\ref{eq:CLT tightness cond}) is satisfied if
the functions $r(t)$ and $a(t)$ satisfy the following:
there exist constants $C>0$, $\alpha^*$ and $\beta^*$,  such that
\begin{equation}\label{eq:CLT tightness cond power}
|r(t)|\le C  (1\wedge|t|^{\alpha^*-1}), \qquad |a(t)|\le C (1\wedge |t|^{\beta^*-1}),
\end{equation}
where $0<\alpha^*,\beta^*<1/2$ and $\alpha^*+\beta^*< 1/2$.
Indeed, to see this, note first that $r(\cdot), \ a(\cdot)\in L^\infty(\mathbb{R})$. Then one can choose $p,q\ge 1$ such that $p(\alpha^*-1)<-1$ and $q(\beta^*-1)<-1$, which entails that $r(\cdot)\in L^{p}(\mathbb{R})$ and $a(\cdot)\in L^{q}(\mathbb{R})$. Since $1/p+1/q<2-\alpha^*-\beta^*$ and  $2-\alpha^*-\beta^*>3/2$, one can further choose $p,q$ to satisfy $1/p+1/q \ge 3/2$.}
\end{rem}

The next results, two functional central limit theorems, extend Theorems
\ref{cth1}(A) and \ref{cth2} to weak convergence in the space $C[0,1]$ of the stochastic
process $\widetilde Q_T(\tau)$ to a standard Brownian motion.

\begin{thm}\label{Thm:fCLT C[0,1]}
Let the spectral density $f$ and the generating function $g$ satisfy
condition (\ref{eq:CLT condition 1}). Let the covariance function $r(t)$ and
the generating kernel $a(t)$ satisfy condition (\ref{eq:CLT tightness cond}).
Then we have the following weak convergence in $C[0,1]$:
\[
 \widetilde Q_T(\tau) \Longrightarrow \sigma_0 B(\tau),
\]
where $\widetilde Q_T(\tau)$ is as in (\ref{eq:tilde Q_T(t)}),
$\sigma_0$ is as in (\ref{eq:sigma}), and $B(\tau)$ is a standard Brownian motion.
\end{thm}
\vskip2mm
Recall the class $SV_0(\mathbb{R})$ of slowly varying at zero functions
$u(x)$, $x\in\mathbb{R}$, satisfying the following conditions:
for some $a>0$, $u(x)$ is bounded on $[-a,a]$, $\lim_{x\to0}u(x)=0,$ \
$u(x)=u(-x)$ and $0<u(x)<u(y)$\ for\ $0<x<y<a$.

\begin{thm}\label{Thm:fCLT C[0,1] B}
Assume that the functions $f$ and $g$ are integrable on $\mathbb{R}$
and bounded outside any neighborhood of the origin, and satisfy for some $a>0$
\beq
\label{c-11}
f(x)\le |x|^{-\alpha}L_1(x), \q
|g(x)|\le |x|^{-\beta}L_2(x),\q x\in [-a,a]
\eeq
for some $\alpha<1, \ \beta<1$ with $\alpha+\beta\le1/2$,
where $L_1(x)$ and $L_2(x)$ are slowly varying at zero functions
satisfying
\begin{equation}\label{eq:L1 L2 cond}
L_i\in SV_0(\mathbb{R}),
\quad x^{-(\alpha+\beta)}L_i(x)\in L^2[-a,a], \quad i=1,2.
\end{equation}
 Let, in addition, the covariance function $r(t)$ and
the generating kernel $a(t)$ satisfy condition (\ref{eq:CLT tightness cond}). Then we have the following weak convergence in $C[0,1]$:
\[
 \widetilde Q_T(\tau) \Longrightarrow \sigma_0 B(\tau),
\]
where $\widetilde Q_T(\tau)$ is as in (\ref{eq:tilde Q_T(t)}),
$\sigma_0$ is as in (\ref{eq:sigma}), and $B(\tau)$ is a standard Brownian motion.
\end{thm}
\begin{rem}
{\rm The proofs of Theorems \ref{Thm:CLT-G} - \ref{Thm:fCLT C[0,1] B},
given in Bai et al. \cite{BGT1} uses the method developed in
Ginovyan and Sahakyan \cite{GS2007}, which itself is based on the
approximations of traces of the products of truncated Toeplitz operators
(see Section \ref{cm}).}
\end{rem}

\subsection{Functional NCLT for $Q_T(\tau)$}

\subsubsection{Non-central limit theorems (discrete-time)}
\label{NCLT-d}

For d.t.\  Gaussian processes the problem of description of the limit
distribution of the quadratic form:
\begin{equation}
\label{nc-1m}
Q_T:=\sum_{t=1}^T\sum_{s=1}^T\widehat g(t-s)X(t)X(s),
\q T\in\mathbb{N},
\end{equation}
if it is non-Gaussian, goes back to the classical papers of Rosenblatt \cite{R1}--\cite{R3}.
Later this problem was studied in a series of papers by Fox and Taqqu,
Taqqu, and Terrin and Taqqu (see, e.g., \cite{FT}, \cite{Tq1},
\cite{Tq4}, \cite{TT1}, \cite{TT2}, and references therein).
Suppose that the spectral density $f(\la)$ and
the generating function $g(\la)$ are regularly varying functions
at the origin:
\begin{equation}
\label{nc-2}
f(\la)= |\la|^{-\al}L_1(\la) \q{\rm and} \q
g(\la)= |\la|^{-\be}L_2(\la), \q\al<1, \be<1,
\end{equation}
where $L_1(\la)$ and $L_2(\la)$ are slowly varying functions
at zero, which are bounded on bounded intervals.
The conditions $\al<1$ and $\be<1$ ensure that the Fourier
coefficients of $f$ and $g$ are well defined. When $\al>0$
the model $\{X(t), t\in\mathbb{Z}\}$ exhibits long memory.

It is the sum $\al+\be$ that determines the asymptotic behavior
of the quadratic form $Q_T$. If $\al+\be\le1/2$, then by
Theorem \ref{cth2} the standard normalized quadratic form
\[
T^{-1/2}\left(Q_T-\E[Q_T]\right)
\]
converges in distribution to a Gaussian random variable.
If $\al+\be>1/2$, the convergence to a Gaussian distribution fails.

Consider the embedding of the discrete sequence
$\{Q_T, \,T\in\mathbb{N}\}$ into a conti\-nuous-time process
$\{Q_T(\tau), \, T\in\mathbb{N}, \, \tau\in[0,1]\}$ defined by
\begin{equation}
\label{nc-3}
Q_T(\tau):=\sum_{t=1}^{[T\tau]}\sum_{s=1}^{[T\tau]}\widehat g(t-s)X(t)X(s),
\end{equation}
where $[\,\cdot\,]$ stands for the greatest integer.
The problem of interest is to describe the limit distribution of the
following normalized process:
\beq
\label{NDP}
\widetilde Q_T(\tau):=d_T^{-1}\left(Q_T(\tau)-\E[Q_T(\tau)]\right)
\eeq
where $d_T$ is a suitably chosen normalization factor.

In \cite{R1} (see also \cite{R3}) Rosenblatt showed that if a
d.t.\  centered Gaussian process $X(t)$ has covariance function
$r(t)=(1+t^2)^{\al/2-1/2}$ with $1/2<\al<1$, then the random variable
\[
\widetilde Q_T:=T^{-\al}\sum_{t=1}^T\left[X^2(t)-1\right]
\]
has a non-Gaussian limiting distribution, and described this distribution
in terms of characteristic functions. This is a special case of
\eqref{NDP} with $\tau=1$, $1/2<\al<1$, $\be=0$ and $d_T=T^{\al}$.
In \cite{Tq1} (see also \cite{Tq4}) Taqqu extended Rosenblatt's result
by showing that the stochastic process
\[
\widetilde Q_T(\tau):=T^{-\al}\sum_{t=1}^{[T\tau]}\left[X^2(t)-1\right]
\]
converges (as $T\to\f$) weakly in $D[0,1]$ to a process (called the Rosenblatt process)
which has the double Wiener-It\^o integral representation
%
\begin{equation}
\label{nc-6}
Q(\tau):=C_\al\int_{\mathbb{R}^2}^{''}\frac{e^{i\tau(x+y)}-1}{i(x+y)}
|x|^{-\al/2}|y|^{-\al/2}dZ(x)dZ(y),
\end{equation}
where $Z(\cd)$ is a complex-valued Gaussian random measure defined
on the Borel $\si$-algebra $\mathcal{B}(\mathbb{R})$, and satisfying
$EZ(B)=0$, $E|Z(B)|^2=|B|$, and $\ol{Z(-B)}=Z(B)$ for any
$B\in\mathcal{B}(\mathbb{R})$.
The double prime in the integral (\ref{nc-6}) indicates that the
integration excludes the diagonals $x=\pm y$.

Notice that the distribution of the random variable $Q(\tau)$ in
(\ref{nc-6}) for $\tau =1$ is described in Veillette and Taqqu \cite{VT}.

The next result, proved in Terrin and Taqqu \cite{TT1} using power
counting theorems (see Section \ref{pct}),
describes the non-Gaussian limit distribution
of the suitable normalized process $Q_T(\tau)$.
\begin{thm}
\label{NCT}
Let $f(\la)$ and $g(\la)$ be as in (\ref{nc-2}) with $\al<1$, $\be<1$
and slowly varying at zero and bounded on bounded intervals factors
$L_1(\la)$ and $L_2(\la)$. Let the process $Q_T(\tau)$ be as in (\ref{nc-3}).
Then for $\al+\be>1/2$ the process:
\begin{equation}
\label{1-16}
Z_T(\tau):=\frac1{T^{\al+\be}L_1(1/T)L_2(1/T)}\,
\left(Q_T(\tau)-\E[Q_T(\tau)]\right)
\end{equation}
converges (as $T\to\f$) weakly in $D[0,1]$ to
\begin{equation}
\label{nc-4}
Z(\tau):=\int_{\mathbb{R}^2}^{''}K_\tau(x,y)dZ(x)dZ(y),
\end{equation}
where
\begin{equation}
\label{nc-5}
K_\tau(x,y)=|xy|^{-\al/2}\int_\mathbb{R}\frac{e^{i\tau(x+u)}-1}{i(x+u)}
\cd\frac{e^{i\tau(y-u)}-1}{i(y-u)}|u|^{-\be}du,
\end{equation}
The double prime in the integral (\ref{nc-4}) indicates that the
integration excludes the diagonals $x=\pm y$.
\end{thm}

\begin{rem}
{\rm The limiting process in (\ref{nc-4}) is real-valued, non-Gaussian,
and satisfies $\E[Z(\tau)]=0$ and $\E[Z^2(\tau)]=\int_{\mathbb{R}^2}|K_\tau(x,y)|^2dxdy$.
It is self-similar with parameter $H=\al+\be\in(1/2,2)$, that is, the
processes $\{Z(a\tau), \tau\in[0,1]\}$ and $\{a^HZ(\tau), \tau\in[0,1]\}$
have the same finite dimensional distributions for all $a>0$.}
\end{rem}

\subsubsection{Non-central limit theorems (continuous-time)}
\label{NCLT-c}

Now we state a \emph{non-central limit theorem} for the process $Q_T(\tau)$
in the c.t.\  case, that is, when $Q_T(\tau)$ in \eqref{eq:quad proc}
is defined through integrals.
Let the spectral density $f$ and the generating function $g$ satisfy
\beq
\label{fnc-2}
f(x)= |x|^{-\al}L_1(x) \q {\rm and} \q
g(x)= |x|^{-\be}L_2(x), \q x\in \mathbb{R},\q \al<1, \, \be<1,
\eeq
with   slowly varying at zero  functions $L_1(x)$ and $L_2(x)$
such that
\beq
\nonumber
\int_{\mathbb{R}}|x|^{-\alpha}L_1(x)dx<\infty \q {\rm and} \q
\int_{\mathbb{R}}|x|^{-\beta}L_2(x)dx<\infty.
\eeq
We assume in addition that the functions  $L_1(x)$ and $L_2(x)$ satisfy
the following condition, called Potter's bound (see Giraitis et al.
\cite{GKSu}, formula (2.3.5)):
{\it for any $\epsilon>0$ there exists a constant $C>0$ so that if $T$ is large enough, then}
\begin{equation}\label{eq:potter}
\frac{L_i(u/T)}{L_i(1/T)}\le C (|u|^{\epsilon}+|u|^{-\epsilon}), \quad i=1,2.
\end{equation}
Note that a  sufficient condition for (\ref{eq:potter}) to hold is that
$L_1(x)$ and $L_2(x)$ are bounded on intervals $[a,\infty)$ for any $a>0$,
which is the case for the slowly varying functions in Theorem \ref{Thm:fCLT C[0,1] B}.

Now we are interested in the limit process of the following normalized version
of the process $Q_T(\tau)$ given by (\ref{eq:quad proc}), with $f$ and $g$ as in (\ref{fnc-2}):
\begin{equation}\label{eq:Z_T}
Z_T(\tau):= \frac{1}{T^{\alpha+\beta} L_1(1/T)L_2(1/T)}\left( Q_T(\tau) -\E [Q_T(\tau)] \right).
\end{equation}

\begin{thm}\label{Thm:fNCLT}
Let $f$ and $g$ be as in (\ref{fnc-2}) with $\alpha<1$, $\beta<1$
and slowly varying at zero functions $L_1(x)$ and $L_2(x)$ satisfying (\ref{eq:potter}),
and let $Z_T(\tau)$ be as in (\ref{eq:Z_T}).
Then if $\alpha+\beta>1/2$, we have the following weak convergence
in the space $C[0,1]$:
\begin{equation*}
Z_T(\tau)\Longrightarrow Z(\tau),
\end{equation*}
where the limit process $Z(\tau)$ is given by
\begin{equation}\label{eq:limit proc}
Z(\tau)=\int_{\mathbb{R}^2}'' H_\tau(x_1,x_2) W(dx_1) W(dx_2),
\end{equation}
with
\begin{equation}\label{eq:H_t}
H_\tau(x_1,x_2)=|x_1x_2|^{-\alpha/2}\int_{\mathbb{R}} \left[\frac{e^{i\tau(x_1+u)}-1}{i(x_1+u)}\right]\cdot\left[\frac{e^{i\tau(x_2-u)}-1}{i(x_2-u)}\right] |u|^{-\beta}du~,
\end{equation}
where $W(\cdot)$ is a complex Gaussian random measure with Lebesgue control measure,
and the double prime in the integral (\ref{eq:limit proc}) indicates that the
integration excludes the diagonals $x_1=\pm x_2$.
\end{thm}

\begin{rem}
{\rm Let $P_T$ and $P$ denote the measures generated in $C[0,1]$
by the processes $Z_T(\tau)$ and $Z(\tau)$ given by (\ref{eq:Z_T}) and
(\ref{eq:limit proc}), respectively. Then Theorem \ref{Thm:fNCLT} can be
restated as follows: under the conditions of Theorem \ref{Thm:fNCLT},
the measure $P_T$ converges weakly in $C[0,1]$ to the measure $P$ as
$T\rightarrow\infty$.
Similar assertions can be stated for Theorems \ref{Thm:fCLT C[0,1]}
and \ref{Thm:fCLT C[0,1] B}.}
\end{rem}

\begin{rem}
{\rm Comparing Theorems \ref{NCT} and \ref{Thm:fNCLT},
we see that the limit process $Z(\tau)$ is the same both for continuous
and discrete time models.
Also, it is worth noting that although the statements of Theorems \ref{NCT}
and \ref{Thm:fNCLT} are similar, the proofs are different.
The proof of Theorem \ref{NCT}, given in Terrin and Taqqu \cite{TT1},
uses hard technical analysis based on power counting theorems,
while the proof of Theorem \ref{Thm:fNCLT}, given in Bai et al. \cite{BGT2},
is simple and uses the spectral representation of the underlying process
and properties of Wiener-It\^o integrals.}
\end{rem}

\section{Functional limit theorems for L\'evy-driven linear models}
\label{FLTL}

\label{FLTG}
In this section, we survey results involving functional central and
non-central limit theorems for a suitably normalized stochastic process:
\begin{equation}\label{eq:Q_T(t)}
Q_T(\tau):=\int_{0}^{T\tau}\int_0^{T\tau} \widehat g(t-s) X(t)X(s) dtds,
\q \tau \in [0,1]
\end{equation}
in the general case where the underlying model $\{X(t), \ t\in\mathbb{R}\}$
is a c.t.\  linear process driven from {\it L\'evy noise} $\xi(t)$ with
time invariant filter $a(\cdot)$, that is,
\begin{equation}
\label{clpL}
X(t)=\int_{\mathbb{R}} a(t-s)d\xi(s), \q
         \int_{\mathbb{R}}|a(s)|^2ds < \f.
\end{equation}

These theorems, which were stated and proved in Bai et al. \cite{BGT2},
extend the results stated in Section \ref{FLTGL} for
Wiener-driven processes 
and show that under some $(L^p,L^q)$-type conditions imposed on the
filter $a(\cdot)$ and the generating kernel $\widehat g(\cdot)$ of
the quadratic functional, the process $Q_T(\tau)$ obeys a central
limit theorem, that is, the finite-dimensional distributions of the
standard $\sqrt{T}$ normalized process $Q_T(\tau)$ in \eqref{eq:Q_T(t)}
tend to those of a normalized {\it standard Brownian motion}.
In contrast, when the functions $a(\cdot)$ and $\widehat g(\cdot)$
have slow power decay, then we have a non-central limit theorem for
$Q_T(\tau)$, that is, the finite-dimensional distributions of the
process $Q_T(\tau)$, normalized by $T^{\gamma}$ for some $\gamma>1/2$,
tend to those of a non-Gaussian non-stationary-increment self-similar
process which can be represented by a double stochastic Wiener-It\^o
integral on $\mathbb{R}^2$.

We point out that the proofs of the central limit theorems given in
Bai et al. \cite{BGT2} were based on an approximation approach,
which reduces the quadratic integral form to a single integral form,
while the proofs of the non-central limit theorems, used the spectral representation of the
underlying process, the properties of Wiener-It\^o integrals, and a continuous analog
of a method to establish convergence in distribution of quadratic functionals to
double Wiener-It\^o integrals, developed by Surgailis \cite{Su1982}
(see also Giraitis et al. \cite{GKSu}).

It is worth noting that if the underlying process $X(t)$ is not necessarily Gaussian,
additional complications arise due to the contribution of the random diagonal
term in the double stochastic integral with respect to L\'evy noise,
which is not present in the case of Gaussian noise (see Remark \ref{cont-diagonal} below).
For this model, in Avram et al. \cite{AvLS}, a central limit theorem for the quadratic functional
$Q_T(1)$ was stated (without proof) under some $(L^p,L^q)$-type
conditions imposed on the spectral density $f(\la)$ and the
generating function $g(\la)$ (see Remark \ref{Rem:avram} below).
For a related study of the sample covariances of L\'evy-driven moving
average processes we refer to the papers by Cohen and Lindner
\cite{CL2013}, and Spangenberg \cite{Sp2015}.

In this section, we follow the paper Bai et al. \cite{BGT2},
and will use the following notation.\\
The symbol $\ast$ will stand for the convolution:
\[(\phi\ast \psi)(u)=\int_{\mathbb{R}}\phi(u-x)\psi(x)dx,
\]
while the symbol  $\bar{\ast}$ is used to denote the reversed convolution:
\[(\phi^{\bar{\ast}2})(u)=(\phi\bar{\ast} \phi)(u) =\int_{\mathbb{R}}\phi(u+x)\phi(x)dx.
\]
By $\mathcal{F}$ and $\mathcal{F}^{-1}$ we denote the Fourier  and the
inverse Fourier transforms:
\[
(\mathcal{F}\phi)(u)=\wh{\phi}(u)=\int_{\mathbb{R}} e^{ixu} \phi(x)dx, \quad
(\mathcal{F}^{-1}\phi)(u)=\frac{1}{2\pi}\int_{\mathbb{R}} e^{-ixu} \phi(x)dx.
\]
We will use the following well-known identities:
\begin{equation}\label{eq:fourier 1}
\mathcal{F}(\phi\ast \psi)= \mathcal{F}(\phi)\cdot \mathcal{F}(\psi)
\q \text{and} \q
\mathcal{F}(\phi\bar{\ast} \psi)= \mathcal{F}(\phi)\cdot \overline{\mathcal{F}(\psi)}.
\end{equation}

\subsection{Central limit theorems}

Let $\{X(t), \ t\in\mathbb{R}\}$ be a centered real-valued linear process
given by (\ref{clpL}) with  filter $a(\cdot)\in L^2(\mathbb{R})$.
Recall that the covariance function $r(t)$ and the spectral density
$f(\la)$ of $X(t)$ are given by formulas (\ref{ccv}) and (\ref{csd}),
respectively.

The theorem that follows contains $(L^p,L^q)$-type time-domain sufficient
conditions for the process $Q_T(\tau)$ to obey the central limit theorem
(see Bai et al. \cite{BGT2}).
\begin{thm}\label{Thm:CLT}
Let $X(t)$ be as in (\ref{clpL}), and let $Q_T(\tau)$ be as in (\ref{eq:Q_T(t)}).
Assume that
\begin{equation}
\label{eq:CLT condition-L}
a(\cdot)\in L^p(\mathbb{R})\cap L^2(\mathbb{R}),
\quad \widehat g(\cdot)\in L^q(\mathbb{R})
\end{equation}
with
\begin{equation}\label{eq:CLT condition p,q-L}
1\le p,q\le 2,\quad \frac{2}{p}+\frac{1}{q}\ge \frac{5}{2}.
\end{equation}
Then
\begin{equation}\label{eq:conv clt}
\widetilde{Q}_T(\tau):= T^{-1/2}\left(Q_T(\tau)- \E [Q_T(\tau)]\right)
\ConvFDD \sigma B(\tau),
\end{equation}
where the symbol $\ConvFDD$ stands for convergence of finite-dimensional distributions, $B(\tau)$ is a standard Brownian motion, and
\begin{equation}\label{eq:sigma^2}
\sigma^2=\int_{\mathbb{R}} [2K_A(v) +\kappa_4 K_B(v)]~dv,
\end{equation}
where $\kappa_4$ is the fourth cumulant of $\xi(1)$, and
\begin{equation}\label{eq:R_A R_B}
K_A(v)=\Big((a\ast b)^{\bar{\ast}2}\cdot a^{\bar{\ast}2}\Big)(v),\quad K_B(v)
= \Big((a{\ast}b)\cdot a\Big)^{\bar{\ast}2}(v).
\end{equation}
\end{thm}

\begin{rem}
{\rm Young's convolution inequality 
(see, e.g., Bogachev \cite{Bog2007}, Theorem 3.9.4)
states that for any numbers $p,p_1, q$ satisfying
$1\le p,q\le p_1\le\f$ and $\frac{1}{p_1}=\frac{1}{p}+\frac{1}{q}-1$,
and for any functions $f\in L^p(\mathbb{R})$,
$g\in L^q(\mathbb{R})$, the function $f\ast g$ is defined almost everywhere,
$f\ast g\in L^{p_1}(\mathbb{R})$, and 
\begin{equation}\label{eq:young}
\|f\ast g\|_{p_1}\le \|f\|_p \|g\|_q.
\end{equation}
Applying this inequality to the convolution in (\ref{ccv}), we get
$\|r\|_{p_1}\le \|a\|_p^2<\infty$, where $1+1/p_1=2/p$. Hence the relations
(\ref{eq:CLT condition-L}) and (\ref{eq:CLT condition p,q-L}) imply that
\begin{equation}\label{eq:usual time domain}
r(\cdot)\in L^{p_1}(\mathbb{R}),\quad b(\cdot)\in L^q(\mathbb{R}),\quad
\frac{1}{p_1}+\frac{1}{q}=\frac{2}{p}-1+\frac{1}{q}\ge \frac{5}{2}-1= \frac{3}{2}.
\end{equation}
The condition (\ref{eq:usual time domain})  is sufficient for the convergence in
Theorem \ref{Thm:CLT} to hold in the case where $\xi(t)$ is Brownian motion
(see Theorem \ref{Thm:fCLT C[0,1]}).
In fact, in this case, the convergence in Theorem \ref{Thm:CLT} holds under
even a weaker condition imposed  on the generating function $g(\lambda)$ and
the spectral density $f(\lambda)$ of $X(t)$
(see Theorem \ref{Thm:CLT-G})}.
\end{rem}
\begin{rem}
\label{cont-diagonal}
{\rm In contrast to the cases where the model is either a d.t.\
linear process (see Giraitis and Surgailis \cite{GSu}), or a c.t.\
Gaussian process (see Bai et al. \cite{BGT1}), it is convenient to impose
the time-domain conditions (\ref{eq:CLT condition-L}) and (\ref{eq:CLT condition p,q-L}) on the functions $a(\cdot)$ and $\widehat g(\cdot)$,
instead  of on the
spectral density $f(\lambda)$ and the generating function $g(\lambda)$.
This allows the analysis of the random diagonal term which arises from the
double stochastic integral with respect to a non-Gaussian L\'evy process.
In the d.t.\  case the random diagonal term is estimated by the full
double sum (see Giraitis and Surgailis \cite{GSu}, relation (2.3)),
while in the c.t.\  Gaussian case, there is no such  random diagonal term
(see Ginovyan and Sahakyan \cite{GS2007}).
In the c.t.\  non-Gaussian case, there is a random diagonal term
in the form of a single stochastic integral that cannot be controlled by
the  double integral, and hence it should be treated separately
(see the proof of Theorem 2.1 of Bai et al. \cite{BGT2})}.
\end{rem}

\begin{rem}
{\rm Observe that the long-run variance $\sigma^2$ given by
(\ref{eq:sigma^2}) can be expressed in terms of the spectral density
$f(\lambda)$ and the generating function $g(\lambda)$, provided that
these functions satisfy some regularity conditions.
Indeed, using the equalities in (\ref{eq:fourier 1}) and the Parseval-Plancherel
theorem, under suitable integrability conditions on $a(\cdot)$ and $b(\cdot)$,
we can write
\begin{align*}
\int_{\mathbb{R}}K_A(v)dv&= \int_{\mathbb{R}}(a\ast b)^{\bar{\ast}2}(v) a^{\bar{\ast}2}(v) dv
=\frac{1}{2\pi} \int_{\mathbb{R}} \mathcal{F}\left((a\ast b)^{\bar{\ast}2}\right)(\lambda)
\overline{\mathcal{F}\left(a^{\bar{\ast}2}\right)(\lambda)} d\lambda=\\
&=\frac{1}{2\pi} \int_{\mathbb{R}} |\mathcal{F}(a\ast b)(\lambda)|^2 |\mathcal{F}(a)(\lambda)|^2 d\lambda
= \frac{1}{2\pi}\int_{\mathbb{R}} |\wh{a}(\lambda) \wh{b}(\lambda)|^2 |\wh{a}(\lambda)|^2 d\lambda\\
&= 8\pi^3 \int_{\mathbb{R}} f(\lambda)^2 g(\lambda)^2 d\lambda,
\end{align*}
where in the last equality  we used the fact $|\wh{a}|^2=2\pi f$ and  $\wh{b}=2\pi g$
(because $b(\cdot)$ is an even function).
 Similarly, we have
\begin{align*}
\int_{\mathbb{R}}K_B(v)dv&= \int_{\mathbb{R}}dv\int_{\mathbb{R}}dx\Big((a{\ast}b)
\cdot a\Big)(x)\Big((a{\ast}b)\cdot a\Big)(x+v)
= \left( \int_{\mathbb{R}} (a{\ast}b)(x) a(x) dx\right)^2\\
&=\frac{1}{4\pi^2}\left( \int_{\mathbb{R}} \wh{a}(x) \wh{b}(x)
\overline{\wh{a}(x)} d\lambda\right)^2
= 4\pi^2 \left[\int_{\mathbb{R}} f(\lambda) g(\lambda) d\lambda\right]^2.
\end{align*}
So an alternative expression for $\sigma^2$ in (\ref{eq:sigma^2}) is
\begin{equation}\label{eq:long run var}
\sigma^2=16\pi^3 \int_{\mathbb{R}} f(\lambda)^2 g(\lambda)^2 d\lambda
+\kappa_4 \left[2\pi\int_{\mathbb{R}} f(\lambda) g(\lambda) d\lambda\right]^2,
\end{equation}
which should be compared with Avram et al. \cite{AvLS} (Theorem 4.1),
and formula \eqref{var} for an analogous expression in the d.t.\  case.} 
\end{rem}

\begin{rem}
{\rm The d.t.\  analog of Theorem \ref{Thm:CLT} with $t=1$ and $\xi$
being Gaussian was established in Giraitis and Surgailis \cite{GSu}.
A special case of Theorem \ref{Thm:CLT} with $t=1$ and $\xi$  being Gaussian
was established in Ginovian \cite{G1994}, and in Ginovian and Sahakyan
\cite{GS2007}. Theorem \ref{Thm:CLT}  for Wiener-driven model ($\kappa_4=0$)
was proved in Bai et al. \cite{BGT1}.}
\end{rem}
\begin{rem}\label{Rem:avram}
{\rm For L\'evy-driven model with $t=1$ and $\sigma^2$ given by
(\ref{eq:long run var}), a version of Theorem \ref{Thm:CLT} was stated in
Avram et al. \cite{AvLS} (Theorem 4.1). They impose $(L^p,L^q)$-type conditions on the
spectral density $f(\cdot)$ and the generating function $g(\cdot)$, and
assume the existence of all moments of the driving L\'evy process $\xi(t)$.
The details of the proof of Theorem 4.1 in Avram et al. \cite{AvLS}
is unfortunately omitted, and it is not clear how the omitted
details of the method-of-moment proof can be carried out given the complexity
of computing the moments of multiple integrals with respect to non-Gaussian
L\'evy noise (see Peccati and Taqqu \cite{PT2011}, Chapter 7).}
\end{rem}

The following corollary contains sufficient conditions for the assumptions
in Theorem \ref{Thm:CLT} to hold.
These conditions involve  bounds on the tails of functions  $a(\cdot)$ and
$b(\cdot)$ by suitable power functions (see Bai et al. \cite{BGT2}).
\begin{cor}\label{Cor:clt sufficient}
The convergence in (\ref{eq:conv clt}) holds if the functions
$a(\cdot)$ and $b(\cdot)$ satisfy the following conditions:
\begin{equation}\label{eq:a b time domain}
a(\cdot),\,\, b(\cdot)\in L^{\infty}(\mathbb{R}),\quad
|a(x)|\le c |x|^{\alpha/2-1}, \quad |b(x)|\le  c|x|^{\beta-1}
\end{equation}
with
\[
0<\alpha,\beta<1, \quad\alpha+\beta<1/2.
\]
\end{cor}

\subsection{Non-central limit theorems}
\label{NCLT-cL}

We now state the non-central limit theorems.
We make the following assumptions on the functions $a(\cdot)$ and $b(\cdot)$,
and on their Fourier transforms $\wh{a}(\cdot)$ and $\wh{b}(\cdot)$.
\begin{asn}
\label{(As1)}
{\rm The Fourier transform $\wh{a}(\cdot)$ of $a(\cdot)\in L^2(\mathbb{R})$ satisfies
\[
\wh{a}(x)= A(x)|x|^{-\alpha/2}L_1^{1/2}(x),
\]
where $L_1(x)$ is an even non-negative function slowly varying at zero and bounded
on intervals $[c,\infty)$ for any $c>0$, and   $A(x)$ is a complex-valued function
satisfying $|A(x)|=1$, and  $\lim_{x\rightarrow 0^+} A(x)=A_0$ for some $A_0$ on
the complex unit circle (since $\wh{a}(-x)=\overline{\wh{a}(x)}$, we also have
$\lim_{x\rightarrow 0^-} A(x)=\overline{A_0}$).}
\end{asn}
\begin{asn}
\label{(As2)}
{\rm The   generating function $\wh{b}(\cdot)\in L^1(\mathbb{R})$ and satisfies
\[
\wh{b}(x)=|x|^{-\beta}L_2(x),
\]
where $L_2(x)$ is an even non-negative function slowly varying at zero and
bounded on intervals $[c,\infty)$ for any $c>0$.}
\end{asn}
\begin{asn}
\label{(As3)}
{\rm The parameters $\alpha$ and $\beta$ above satisfy
\begin{equation}\label{eq:alpha beta range}
-1/2<\alpha<1, \quad -1/2<\beta<1,\quad \alpha+\beta>1/2.
\end{equation}}
\end{asn}
\begin{asn}
\label{(As4)}
{\rm There exist numbers $\alpha^*$ and $\beta^*$ satisfying
\[
0<\alpha^*,\beta^*<1\quad 1<\alpha^*+\beta^*<\alpha+\beta+1/2,
\]
such that
\[
|a(x)|\le C |x|^{\alpha^*/2-1},\quad |b(x)|\le C|x|^{\beta^*-1}.
\]}
\end{asn}

The proof of the following theorem can be found in Bai et al. \cite{BGT2}.
\begin{thm}\label{Thm:NCLT}
Suppose that Assumptions \ref{(As1)} - \ref{(As4)} hold. Then as $T\rightarrow\infty$ we have
\begin{equation}\label{eq:nclt}
Z_T(\tau):= \frac{1}{T^{\alpha+\beta}L_1(1/T)L_2(1/T)}
\left(Q_T(\tau)- \E Q_T(\tau)\right)\ConvFDD  Z_{\alpha,\beta}(\tau),
\end{equation}
where
\begin{equation}\label{eq:Z_alpha,beta}
Z_{\alpha,\beta}(\tau)=\frac{1}{2\pi}\int_{\mathbb{R}^2}'' |x_1x_2|^{-\alpha/2}
\int_{\mathbb{R}} \frac{e^{i\tau(x_1+u)}-1}{i(x_1+u)}
\frac{e^{i\tau(x_2-u)}-1}{i(x_2-u)}  |u|^{-\beta}du~  W(dx_1)W(dx_2),
\end{equation}
where the double prime $''$ indicates the exclusion of the hyper-diagonals
$u_p=\pm u_q$, $p\neq q$ and $W(\cdot)$ is a complex-valued Brownian motion
(see Section \ref{sec:prelim}). 
\end{thm}

\begin{rem}\label{Rem:reg fourier}
{\rm The regular variation conditions on $\wh{a}(\cdot)$ and $\wh{b}(\cdot)$
in Assumptions \ref{(As1)} - \ref{(As3)} generally do not follow from the corresponding regular
variation conditions imposed on the inverse Fourier transforms $a(\cdot)$
and $b(\cdot)$. This implication only holds under some additional assumptions
on the slowly varying factors $L_1(\cdot)$ and $L_2(\cdot)$ of $a(\cdot)$
and $b(\cdot)$. For instance, it will hold if we have (see Bingham et al.
\cite{BGT}, formula (4.3.7))
\begin{equation}\label{eq:a,b power}
a(x)=x^{\alpha/2-1}\ell_1(x)1_{[0,\infty)}(x),\quad
b(x)= |x|^{\beta-1}\ell_2(x),
\end{equation}
where $0<\alpha<1$, $0<\beta<1$, $\alpha+\beta>1/2$, and  $\ell_1(x)$ and
$\ell_2(x)$ are even non-negative functions which are locally bounded,
slowly varying at infinity  and \emph{quasi-monotone}.
Recall that a slowly varying function $l(\cdot)$ is said to be quasi-monotone
if it has locally bounded variation, and for all $\delta>0$, one has
(see Bingham et al. \cite{BGT}, Section 2.7):
\[
\int_0^x t^\delta |d\ell(t)|=O(x^\delta l(x))\quad {\rm as} \quad x\to\infty.
\]
A sufficient condition for a slowly varying $\ell(x)$ with locally bounded
variation to be quasi-monotone is that  $x^{\delta}\ell(x)$ is increasing
and $x^{-\delta}\ell(x)$ is decreasing when $x$ is large enough, for any
$\delta>0$ (see Theorem 1.5.5 and Corollary 2.7.4 in Bingham et al. \cite{BGT}).

Notice also that Assumption \ref{(As4)} will be satisfied if (\ref{eq:a,b power}) holds.}
\end{rem}

\begin{rem}
{\rm Let the functions $a(\cdot)$ and $b(\cdot)$ be as in (\ref{eq:a,b power})
with $\alpha<0$ or $\beta<0$ (by (\ref{eq:alpha beta range}) only one of
$\alpha$ and $\beta$ can be negative).
Assume that $\alpha<0$ and $\beta>0$. Then for the corresponding regular
variation of $\wh{a}(\cdot)$ to hold, one needs to impose in addition
that $\int_0^\infty a(x)dx=0$. In this case, one does not need to assume
quasi-monotonicity for $\ell_1$ (see Corollary 1.40 of Soulier \cite{Sou}).
Similar considerations hold if $\beta<0$ and $\alpha>0$ instead.}
\end{rem}

\begin{rem}
{\rm Note that Assumption \ref{(As1)} holds with $\alpha=0$ if $a(\cdot)\in L^1(\mathbb{R})$
and $\int_0^\infty a(x)\neq 0$, and Assumption \ref{(As2)} holds with $\beta=0$ if
$b(\cdot)\in L^1(\mathbb{R})$ and $\int_0^\infty b(x)\neq 0$.}
\end{rem}

The next theorem contains  time-domain representations for the limiting process
$Z_{\alpha,\beta}(\tau)$ in (\ref{eq:Z_alpha,beta}) in the case $\alpha,\beta\ge 0$ (see Bai et al. \cite{BGT2}, Theorem 2.3).
\begin{thm}\label{Pro:time domain rep}
The limiting process $Z_{\alpha,\beta}(\tau)$ in (\ref{eq:Z_alpha,beta})
admits the following time-domain representations:
\begin{enumerate}[(a)]
\item when $\alpha>0$, $\beta>0$:
\begin{equation}\label{eq:alpha>0 beta>0}
Z_{\alpha,\beta}(\tau)\EqFDD c_{\alpha,\beta}
\int_{\mathbb{R}^2}'
\int_0^\tau\int_0^\tau |u-v|^{\beta-1} (u-x_1)_+^{\alpha/2-1} (v-x_2)_+^{\alpha/2-1} du dv
~ B(dx_1)B(dx_2),
\end{equation}
where
$$c_{\alpha,\beta}=\frac{\Gamma(1-\beta)\sin(\beta\pi/2)}{\pi\Gamma(\alpha/2)^2};$$
\item when $\alpha>1/2$, $\beta=0$:
\begin{equation}\label{eq:alpha>0 beta=0}
Z_{\alpha,\beta}(\tau)\EqFDD c_\alpha
\int_{\mathbb{R}^2}' \int_0^\tau  (u-x_1)_+^{\alpha/2-1} (u-x_2)_+^{\alpha/2-1} du  ~ B(dx_1)B(dx_2),
\end{equation}
where
$$c_{\alpha}=\frac{\sin(\alpha\pi/2)\Gamma(1-\alpha/2)}{\pi \Gamma(\alpha/2)};$$
\item when $\alpha=0$, $\beta>1/2$:
\begin{equation}\label{eq:alpha=0 beta>0}
Z_{\alpha,\beta}(\tau)\EqFDD c_\beta\int_{[0,\tau]^2}'
|x_1-x_2|^{\beta-1} ~ B(dx_1)B(dx_2),
\end{equation}
where
$$c_\beta=\frac{\Gamma(1-\beta)\sin(\beta\pi/2)}{\pi}.$$
Here $B(\cdot)$ is the real Brownian random measure and the prime $'$ in the integrals
indicates the exclusion of the diagonals.
\end{enumerate}
\end{thm}
\begin{rem}
{\rm In view of (\ref{eq:Q_T(t)}) and (\ref{eq:a,b power}),
the representation (\ref{eq:alpha>0 beta>0})  gives an explicit insight
of the convergence in Theorem \ref{Thm:NCLT} (see Theorem \ref{Thm:nclt time domain} below).
The process in (\ref{eq:alpha>0 beta=0}) is known as  Rosenblatt process
(see Taqqu \cite{Tq1}), and the corresponding convergence in Theorem
\ref{Thm:NCLT} is the c.t.\  analog of the d.t.\  case
considered in Fox and Taqqu \cite{FT}.
The representation  (\ref{eq:alpha=0 beta>0}) is obtained because for
$\alpha=0$, the underlying process $X(t)$ has short memory and in this case,
one expects that in the limit  $X(t)dt$ in (\ref{eq:Q_T(t)}) can be replaced
by the white noise $B(dt)$.}
\end{rem}
\begin{rem}
{\rm It is of interest to obtain appropriate elementary expressions for the
time-domain representation of the limiting process $Z_{\alpha,\beta}(\tau)$,
similar to \eqref{eq:alpha>0 beta>0} - \eqref{eq:alpha=0 beta>0},
in the cases where either $\alpha$ or $\beta$ satisfying (\ref{eq:alpha beta range})
is negative.}
\end{rem}

Using the time-domain representation (\ref{eq:alpha>0 beta>0}),
one can state a non-central limit theorem  in the case where
$\alpha,\beta>0$ without going to the spectral domain.
This simplifies the assumptions imposed on the functions $a(\cdot)$ and $b(\cdot)$ (see Bai et al. \cite{BGT2}).
\begin{thm}\label{Thm:nclt time domain}
Suppose that the functions  $a(\cdot)$ and $b(\cdot)$ are given
by (\ref{eq:a,b power}), where
 $0<\alpha<1$, $0<\beta<1$, $\alpha+\beta>1/2$, and  $\ell_1(x)$ and $\ell_2(x)$
 are even functions slowly varying at infinity and bounded on bounded intervals.
 Then as $T\rightarrow\infty$ we have
\begin{align*}
&\frac{1}{T^{\alpha+\beta}\ell_1(T)\ell_2(T)}\left(Q_T(\tau)-
\E [Q_T(\tau)]\right)
\ConvFDD \\
&\int_{\mathbb{R}^2}' \int_0^\tau\int_0^\tau |u-v|^{\beta-1}
(u-x_1)_+^{\alpha/2-1} (v-x_2)_+^{\alpha/2-1} du dv ~ B(dx_1)B(dx_2).
\end{align*}
\end{thm}

\section{CLT for Tapered Toeplitz Quadratic Functionals}
\label{Tap1}
\sn{The problem}
\label{H-6.1}

In this section we consider a question concerning asymptotic distribution
(as $T\to\f$) of the following {\it tapered} Toeplitz type quadratic
functionals of the centered stationary process $X(u)$, $u\in \mathbb{U}$,
with spectral density $f(\la)$, $\la\in \Lambda$:
\beq
\label{tq-1}
Q_T^h:=
 \left \{
 \begin{array}{ll}
\sum_{t=1}^T\sum_{s=1}^T\widehat g(t-s)h_T(t)h_T(s)X(t)X(s)
& \mbox{in the d.t.\  case},\\
\\
\int_0^T\int_0^T \widehat g(t-s)h_T(t)h_T(s)X(t)X(s)\,dt\,ds
& \mbox{in the c.t.\  case},
\end{array}
\right. \eeq
where $\widehat g(t)$ is the Fourier transform of some integrable even function $g(\la)$ and
\beq
\label{Tap}
h_T(t):= h(t/T)
\eeq
with a taper function $h(\cdot)$ to be specified below.

The limit distribution of the functional (\ref{tq-1}) is completely
determined by the functions $f$, $g$ and $h$, and depending on their
properties it can be either Gaussian (that is, $Q^h_T$ with an appropriate
normalization obey central limit theorem), or non-Gaussian.

We discuss here the case where the limit distribution is Gaussian,
and present sufficient conditions in terms of functions $f$, $g$ (and $h$)
ensuring central limit theorems for a standard normalized tapered quadratic
functional $Q_T^h$.

We will assume that the taper function  $h(\cdot)$
satisfies the following assumption.

\begin{asn}
\label{(H)}
{\rm The taper $h:\mathbb{R}\to\mathbb{R}$ is a continuous nonnegative
function of bounded variation and of bounded support $[0, 1]$, such that
$H_k\neq 0$, where}
\beq
\label{t4}
H_{k}: = \int_0^1h^k(t)dt, \q k\in \mathbb{N}:=\{1,2,\ldots\}.
\eeq
\end{asn}

\n
{\it Note.} The case where $h(t)={\mathbb I}_{[0,1]}(t)$, where
${\mathbb I}_{[0,1]}(\cdot)$ denotes the indicator of the segment $[0,1]$,
will be referred to as the {\it non-tapered case}.

\begin{rem}
{\rm In the d.t.\  case, an example of a taper function $h(t)$ satisfying
Assumption \ref{(H)} is the Tukey-Hanning taper function
$h(t)=0.5(1-\cos(\pi t))$ for $t\in [0, 1]$.
In the c.t.\  case, a simple example of a taper function $h(t)$ satisfying
Assumption \ref{(H)} is the function $h(t)=1-t$ for $t\in [0, 1]$ (see, e.g., Anh et al. \cite{ALS}).}
\end{rem}

\sn{Statistical motivation}
Much of the statistical inferences (parametric and nonparametric estimation, hypotheses testing)
about the spectrum or the covariance of a stationary process $\{X(u), \ u\in \mathbb{U}\}$
is based on the sample:
\beq \label{1-d}
\mathbf{X}_T: =
\left \{
\begin{array}{ll}
 \{X(t),  \, t=1,\ldots, T\} & \mbox{in the d.t.\  case},\\
\{X(t),  \, 0\leq t\leq T\} & \mbox{in the c.t.\  case}
\end{array}
\right.
\eeq
In the statistical analysis of stationary processes, however, the data are frequently
tapered before calculating the statistics of interest, and the statistical inference
procedure, instead of the original data $\mathbf{X}_T$ given by \eqref{1-d},
is based on the tapered data $\mathbf{X}_T^h$:
\beq \label{t2}
\mathbf{X}_T^h: =
\left \{
\begin{array}{ll}
 \{h_T(t)X(t),  \, t=1,\ldots, T\} & \mbox{in the d.t.\  case},\\
\{h_T(t)X(t),  \, 0\leq t\leq T\} & \mbox{in the c.t.\  case},
\end{array}
\right.
\eeq
with a taper function  $h(t)$, $t\in\mathbb{R}$.

The benefits of tapering the data have been widely reported in the literature
(see, e.g., Bloomfield \cite{Bl}, Brillinger \cite{Bri2},  Dahlhaus \cite{D1,D3, D2},
Dahlhaus and K\"unsch \cite{DK}, Guyon \cite{Gu}, and references therein).
For example, data-tapers are introduced to reduce the so-called 'leakage effects',
that is, to obtain better estimation of the spectrum of the model in the case
where it contains high peaks.
Other application of data-tapers is in situations in which some of the
data values are missing. Also, the use of tapers leads to bias reduction,
which is especially important when dealing with spatial data. In this case,
the tapers can be used to fight the so-called 'edge effects'.

Quadratic functionals of the form \eqref{tq-1} appear both in nonparametric and
parametric estimation of the spectrum of the process $X(t)$ based on the tapered
data \eqref{t2}.
For instance, when we are interested in nonparametric estimation of a linear
integral functional in $L^p(\Lambda)$, $p > 1$ of the form:
\beq
\label{t1}
J=J(f): = \int_\Lambda f(\lambda)g(\lambda)d\lambda,
\eeq
where $g(\lambda) \in L^q(\Lambda)$, \,$1/p + 1/q = 1$, then a natural
statistical estimator for $J(f)$ is the linear integral functional of the
{\it empirical spectral density (periodogram)} of the process $X(t)$.
To define this estimator, we first introduce some notation.

Denote by $H_{k,T}(\la)$ the tapered Dirichlet type kernel, defined by
\beq \label{t3}
H_{k,T}(\la): =
\left \{
\begin{array}{ll}
\sum_{t=1}^T h_{T}^k(t)e^{-i\la t} & \mbox{in the d.t.\  case},\\
\\[-1mm]
\int_0^T h_{T}^k(t)e^{-i\la t}dt & \mbox{in the c.t.\  case}.
\end{array}
\right.
\eeq

Define the finite Fourier transform of the tapered data (\ref{t2}):
\beq
\label{t5}
d^h_T(\la): =
\left \{
\begin{array}{ll}
\sum_{t=0}^T h_{T}(t)X(t)e^{-i\la t} & \mbox{in the d.t.\  case},\\
\\[-1mm]
 \int_0^Th_T(t)X(t)e^{-i\la t}dt & \mbox{in the c.t.\  case}.
\end{array}
\right.
\eeq
and the tapered periodogram $I^h_{T}(\lambda)$ of the process $X(t)$:
\bea
\label{t6}
I^h_{T}(\lambda):= \frac1{C_T}\,d^h_T(\la)d^h_T(-\la)=
\left \{
\begin{array}{ll}
\frac 1{C_T}\left |\sum_{t=0}^T h_{T}(t)X(t)e^{-i\la t}\right|^2 &
\mbox{in the d.t.\  case},\\
\\[-1mm]
\frac1{C_T}\left|\int_0^Th_T(t)X(t)e^{-i\la t}dt\right|^2 &
\mbox{in the c.t.\  case}.
\end{array}
\right.
\eea
where
\beq
\label{t55}
C_T:= 2\pi H_{2,T}(0)=2\pi\int_0^Th_T^2(t)dt=2\pi H_2\,T \neq 0.
\eeq
Notice that for non-tapered case ($h(t)={\mathbb I}_{[0,1]}(t)$),
we have $C_T= 2\pi T$.

As an estimator $J_T^h$ for functional $J(f)$, given by \eqref{t1},
based on the tapered data \eqref{t2}, we consider the averaged tapered
periodogram (or a simple "plug-in" statistic), defined by
\bea
\label{t7}
J_T^h &=& J(I^h_{T}):= \int_\Lambda I^h_{T}(\la)g(\lambda)d\lambda.
\eea
In view of \eqref{tq-1}, (\ref{t6}) and (\ref{t7}) we have
\bea
\label{t77}
J_T^h =C_T^{-1}Q_T^h=\left \{
\begin{array}{ll}
\frac 1{C_T}\sum_{t=1}^T \sum_{s=1}^T\widehat g(t-s)h_T(t)h_T(s)X(t)X(s) &
\mbox{in the d.t.\  case},\\
\\[-1mm]
\frac 1{C_T}\int_0^T\int_0^T \widehat g(t-s)h_T(t)h_T(s)X(t)X(s)\,dt\,ds &
\mbox{in the c.t.\  case},
\end{array}
\right.
\eea
where $C_T$ is as in (\ref{t55}), and  $\widehat g(t)$ is the
Fourier transform of function $g(\la)$.

Thus, to study the asymptotic properties of the estimator $J_T^h$, we have
to study the asymptotic distribution (as $T\to\f$) of the tapered Toeplitz type
quadratic functional $Q_T^h$ given by (\ref{tq-1}) (for details see Section \ref{App}).

\sn{Central limit theorems for tapered quadratic functional $Q_T^h$}
\label{CLT}

We will use the following notation.
By $\widetilde Q^h_T$ we denote the standard normalized
quadratic functional:
\begin{equation}
\label{1-4}
 \widetilde Q^h_T= T^{-1/2}\,\left(Q^h_T-\E [Q^h_T]\right).
\end{equation}

Also, we set
\beq\label{1-8}
\si^2_h: = 16\pi^3H_4 \int_{\Lambda} f^2(\la)g^2(\la)\,d\la,
\eeq
where $H_4$ is as in \eqref{t4}.
The notation
\begin{equation}
\label{1-55}
\widetilde Q^h_T\ConvD \eta \sim N(0,\sigma^2_h) \q {\rm as}\q T\to\f
\end{equation}
will mean that the distribution of the random variable
$\widetilde Q^h_T$ tends (as $T\to\f$) to the centered normal
distribution with variance $\sigma^2_h$.

Let $\psi(\la)$ be an integrable real symmetric function defined on
$[-\pi, \pi]$, and let $h(t)$, $t\in[0,1]$ be a taper function.
For  $T=1, 2,\ldots$, the {\it $(T\times T)$-truncated tapered Toeplitz matrix\/}
generated by $\psi$ and $h$, denoted by $B_T^h(\psi)$, is defined by the following
equation (see \eqref{IMT2-1} for non-tapered case):
\beq
\label{1-7M}
B_T^h(\psi):=\|\widehat\psi(t-s)h_T(t)h_T(s)\|_{t,s=1,2\ldots,T},
\eeq
where $\widehat\psi(t)$ $(t\in \mathbb{Z})$ are the Fourier coefficients of $\psi$.

Given a real number $T>0$ and an integrable real symmetric function $\psi(\la)$
defined on $\mathbb{R}$, the {\it $T$-truncated tapered Toeplitz operator\/}
(also called {\it tapered Wiener-Hopf operator}) generated by $\psi$ and a taper function $h$,
denoted by $W_T^h(\psi)$ is defined as follows (see \eqref{IMT3-1} for non-tapered case):
\begin{equation}
\label{1-7}
[{W}^h_T(\psi)u](t)=\int_0^T\hat\psi(t-s)u(s)h_T(s)ds,
\q u(s)\in L^2([0,T]; h_T),
\end{equation}
where $\hat\psi(\cdot)$ is the Fourier transform of $\psi(\cdot)$,
and $L^2([0,T]; h_T)$ denotes the weighted $L^2$-space with respect to
the measure $h_T(t)dt$.

Let $A_T^h(f)$ be either the $T\times T$ tapered Toeplitz matrix $B_T^h(f)$,
or the $T$-truncated tapered Toeplitz operator $W_T^h(f)$
generated by the spectral density $f$ and taper $h$, and let $A_T^h(g)$
denote either the $T\times T$ tapered Toeplitz matrix, or the $T$-truncated tapered
Toeplitz operator generated by the functions $g$ and $h$.

\ssn{CLT for Gaussian models}
The theorems that follow extend the results of Theorems \ref{cth1} and \ref{cth2}
to the tapered case. We assume that the model process $X(t)$ is Gaussian, and
with no loss of generality, that $g\ge 0$. The following theorems were  proved in
Ginovyan and Sahakyan \cite{GS2019a}.

\begin{thm}\label{th2}
Assume that $f\cdot g\in L^1(\Lambda)\cap L^2(\Lambda)$,
the taper function $h$ satisfies Assumption \ref{(H)}, 
and for $T\to \infty$
\beq\label{1-11}
\chi_2(\widetilde Q^h_T)=\frac2T\tr\bigl[A^h_T(f)A^h_T(g)\bigr]^2 \longrightarrow \si^2_h,
\eeq
where $\si^2_h$ is as in \eqref{1-8}.
Then $\widetilde Q^h_T\ConvD \eta \sim N(0,\sigma^2_h)$ as $T\to\f$.
\end{thm}

\begin{thm}\label{th3}
Assume that the function
\beq\label{1-10}
\varphi(x_1, x_2,x_3)=\int_{\Lambda}
f(u)g(u-x_1)f(u-x_2)g(u-x_3)\,du
\eeq
belongs to $L^2(\Lambda^3)$ and is continuous at $(0,0,0)$,
and the taper function $h$ satisfies Assumption \ref{(H)}.
Then $\widetilde Q^h_T\ConvD \eta \sim N(0,\sigma^2_h)$ as $T\to\f$.
\end{thm}

\begin{thm}\label{th1}
Assume that $f(\la)\in L^p(\Lambda)$ $(p\ge1)$
and  $g(\la)\in L^q(\Lambda)$ $(q\ge1)$
with $1/p+1/q\le1/2$, and the taper function $h$ satisfies Assumption \ref{(H)}.
Then $\widetilde Q^h_T\ConvD \eta \sim N(0,\sigma^2_h)$ as $T\to\f$.
\end{thm}

\begin{thm}\label{th4}
Let $f\in L^2(\Lambda)$, \,$g\in L^2(\Lambda)$, $fg\in L^2(\Lambda)$,
\begin{equation}
\label{1-12} \int_{\Lambda}
 f^2(\la)g^2(\la-\mu)\,d\la \longrightarrow
 \int_{\Lambda} f^2(\la)g^2(\la)\,d\la \quad {\rm as} \quad \mu\to0,
\end{equation}
and let the taper function $h$ satisfy Assumption \ref{(H)}.
Then $\widetilde Q^h_T\ConvD \eta \sim N(0,\sigma^2_h)$ as $T\to\f$.
\end{thm}

To state the next theorem, we recall the class $SV_0(\mathbb{R})$ of slowly
varying functions at zero
$u(\la)$, $\la\in\mathbb{R}$, satisfying the following conditions:
for some $a>0$, $u(\la)$ is bounded on $[-a,a]$, $\lim_{\la\to0}u(\la)=0,$ \
$u(\la)=u(-\la)$ and $0<u(\la)<u(\mu)$\ for\ $0<\la<\mu<a$.

\begin{thm}\label{th5}
Assume that the functions $f$ and $g$ are integrable on $\mathbb{R}$
and bounded outside any neighborhood of the origin, and satisfy for some $a>0$
\beq \label{m-0}
f(\lambda)\le |\lambda|^{-\alpha}L_1(\lambda),
\q
|g(\lambda)|\le |\lambda|^{-\beta}L_2(\lambda),\q \lambda\in [-a,a],
\eeq
for some $\alpha<1, \ \beta<1$ with $\alpha+\beta\le1/2$, where $L_1(x)$ and
$L_2(x)$ are slowly varying functions at zero satisfying
\bea
\label{m-00}
L_i\in SV_0(\mathbb{R}), \ \ \lambda^{-(\alpha+\beta)}L_i(\lambda)\in L^2[-a,a], \ i=1,2.
\eea
Also, let the taper function $h$ satisfy Assumption \ref{(H)}.
Then $\widetilde Q^h_T\ConvD \eta \sim N(0,\sigma^2_h)$ as $T\to\f$.
\end{thm}
As in Remark \ref{Rem:B follow A}, the conditions $\al<1$ and $\be<1$ in Theorem \ref{th5} ensure that the
Fourier transforms of $f$ and $g$ are well defined. Observe that when $\al>0$
the process $X(t)$ may exhibit long-range dependence.
We also allow here  $\alpha+\beta$ to assume the critical value 1/2.
The assumptions $f\cdot g\in L^1(\Lambda)$, $f,g\in L^\infty(\Lambda\setminus [-a,a])$
and (\ref{m-00}) imply that $f\cdot g \in L^2(\Lambda)$, so that the variance
$\sigma^2_h$ in (\ref{1-8}) is finite.

\ssn{CLT for L\'evy-driven stationary linear models}
Now we assume that the underlying model $X(t)$ is a L\'evy-driven stationary
linear process defined by \eqref{clp}, where $a(\cdot)$ is a filter from $L^2(\mathbb{R})$, and $\xi(t)$ is a L\'evy process satisfying the conditions:
$\E \xi(t)=0$, $\E \xi^2(1)=1$ and $\E\xi^4(1)<\infty$.

The central limit theorem that follow is a tapered counterpart of Theorem
\ref{Thm:CLT} and was proved in Ginovyan and Sahakyan \cite{GS2019}.
\begin{thm}
\label{TCLT}
Assume that the filter $a(\cdot)$
and the generating kernel $\widehat g(\cdot)$ are such that
\begin{equation}
\label{eq:CLT condition-M}
a(\cdot)\in L^p(\mathbb{R})\cap L^2(\mathbb{R}),\quad \widehat g(\cdot)\in L^q(\mathbb{R}),
\q1\le p,q\le 2,\quad \frac{2}{p}+\frac{1}{q}\ge \frac{5}{2},
\end{equation}
and the taper $h$ satisfies Assumption \ref{(H)}. 
Then $\widetilde Q^h_T\ConvD \eta \sim N(0,\sigma^2_{L,h})$ as $T\to\f$, where
\beq
\label{tsigma}
\sigma^2_{L,h}=16\pi^3H_4\int_{\mathbb{R}} f^2(\lambda) g^2(\lambda) d\lambda
+\kappa_4 4\pi^2H_4\left[\int_{\mathbb{R}} f(\lambda) g(\lambda) d\lambda\right]^2,
\eeq
and $H_4$ is as in \eqref{t4}.
\end{thm}
\begin{rem}
{\rm
Notice that if the underlying process $X(t)$ is Gaussian,
then in formula (\ref{tsigma}) we have only the first term and so
$\sigma^2_{L,h}=\sigma^2_{h}$ (see \eqref{1-8}), because in this case $\kappa_4=0$.
On the other hand, the condition \eqref{eq:CLT condition-M} is more restrictive
than the conditions in Theorems \ref{th2} - \ref{th5}.
Thus, for Gaussian processes Theorems \ref{th2} - \ref{th5} improve Theorem \ref{TCLT}.}
\end{rem}

\begin{rem}
{\rm Central and non-central limit theorems for tapered quadratic forms
of a d.t.\  long memory Gaussian stationary fields have been proved
in Doukhan et al. \cite{DLS}}.
\end{rem}

\section{Applications}
\label{App}
\sn{Nonparametric estimation of spectral functionals}
\label{AppN}

Suppose we observe a realization $\mathbf{X}_T:=\{X(u)$, $0\le u\le T$
(or $u={1,\ldots,T}$ in the d.t.\  case)\}
of a centered stationary  process $X(t)$ with an {\it unknown\/} spectral
density function $f(\lambda)$, $\lambda \in \Lambda$.
We assume that $f(\lambda)$ belongs to a given (infinite-dimensional)
class $\mathbf{\Theta} \subset L^p(\Lambda)$ $(p \ge 1)$
of  spectral densities possessing some smoothness properties.
Let  $\Phi(\cdot)$ be some {\it known\/} functional,
the domain of definition of which contains $\mathbf{\Theta}$.
The problem is to estimate the value  $\Phi (f)$
of the functional $\Phi(\cdot)$  at an unknown point
$f\in \mathbf{\Theta}$ on the basis of observation $\mathbf{X}_T,$
and investigate the asymptotic (as $T \to \infty$) properties of the
suggested estimators.

This problem for linear and some nonlinear smooth functionals
for d.t.\  and c.t.\  stationary models has been extensively discussed
in the literature
(see, e.g., Dahlhaus and Wefelmeyer \cite{DW}, Ginovyan
\cite{G1988a,G1988b, G1995, G2011a, G2011b}, Ginovyan and Sahakyan \cite{GS2019, GS2019a},
Has'minskii and Ibragimov \cite{IH1}, Taniguchi \cite{Tan87},
Taniguchi and Kakizawa \cite{TK}, and references therein).

In this section we apply the results of Section \ref{CLT} to show that the statistic
$J_T^h$ given by \eqref{t7} is a consistent and asymptotically normal estimator for
the linear functional $J(f)$ given by \eqref{t1}.
We follow the papers Ginovyan and Sahakyan \cite{GS2019, GS2019a}.
To state the corresponding results, we first introduce the $L^p$-H\"older class
and set up a set of assumptions.

Given numbers $p\ge 1$, $0<\al<1$, $r\in \mathbb{N}_0:=\mathbb{N}\cup\{0\}$,
we set $\be=\al+r$  and denote by $H_p(\be)$ the $L^p$-H\"older class,
that is, the class of those functions $\psi(\la)\in L^p(\Lambda)$,
which have  $r$-th derivatives in $L^p(\Lambda)$ and with some positive
constant $C$ satisfy
$$||\psi^{(r)}(\cdot+\la)-\psi^{(r)}(\cdot)||_p\le C|\la|^\al.$$

Now we list the  assumptions.

\begin{asn}
\label{(A1)}
{\rm The filter $a(\cdot)$ and the generating kernel $\widehat g(\cdot)$
are such that
\begin{equation}
\label{eq:CLT condition-N}
a(\cdot)\in L^p(\Lambda)\cap L^2(\Lambda),
\quad \widehat g(\cdot)\in L^q(\Lambda)
\end{equation}
with
\begin{equation}\label{eq:CLT condition p,q-N}
1\le p,q\le 2,\quad \frac{2}{p}+\frac{1}{q}\ge \frac{5}{2}.
\end{equation}}
\end{asn}
\begin{asn}
\label{(A2)}
{\rm The spectral density $f$ and the generating function $g$ are such that
$f, g\in L^1(\Lambda)\cap L^2(\Lambda)$ and $g$  is of bounded variation.}
\end{asn}
\begin{asn}
\label{(A3)}
{\rm The spectral density $f$ and the generating function $g$ are such that
$f\in H_p(\be_1)$, $\be_1>0$, $p\ge1$ and $g(\la)\in H_q(\be_2)$,
$\be_2>0$, $q\ge1$ with $1/p+1/q=1$, 
and one of the conditions a)--d) is satisfied:

a) $\be_1>1/p$, $\be_2>1/q$

b) $\be_1\le1/p$, $\be_2\le1/q$ and $\be_1+\be_2>1/2$

c) $\be_1>1/p$, $1/q-1/2<\be_2\le1/q$

d) $\be_2>1/q$, $1/p-1/2<\be_1\le1/p$.}
\end{asn}

\begin{rem}
\label{RH1}
{\rm In Ginovian \cite{G1994} it was proved that if Assumption \ref{(A3)}
is satisfied, then there exist numbers
$p_1$ $(p_1>p)$ and $q_1$ $(q_1>q)$, such that $H_p(\be_1)\subset L_{p_1}$,
$H_q(\be_2)\subset L_{q_1}$ and $1/{p_1}+1/{q_1}\le1/2$.}
\end{rem}
The next theorem controls the bias $E(J_T^h)-J$ and provides sufficient conditions
assuring the proper rate of convergence of the bias to zero, which is necessary
to obtain the asymptotic normality of the estimator $J_T^h$.
Specifically, we have the following result, which was proved in Ginovyan and Sahakyan \cite{GS2019}.

\begin{thm}
\label{T-Bias}
Let the functionals $J:=J(f)$ and $J_T^h: = J(I^h_{T})$ be defined
by (\ref {t1}) and (\ref{t7}), respectively.
Then under Assumptions \ref{(A2)} (or \ref{(A3)}) and \ref{(H)}
the following asymptotic relation holds:
\bea
\label{pr1}
T^{1/2}\left[\E(J_T^h)-  J\right] \to 0 \quad {\rm as}\quad T\to\infty.
\eea
\end{thm}

The next theorem, which is an immediate consequence of Theorem \ref{TCLT},
contains sufficient conditions for functional $J_T^h$ to obey the central
limit theorem.
\begin{thm}
\label{T-CLT}
Let $J:=J(f)$ and $J_T^h: = J(I^h_{T})$ be defined by (\ref {t1}) and (\ref{t7}),
respectively. Then under Assumptions \ref{(A1)} and \ref{(H)} 
the functional $J_T^h$ obeys the central limit theorem. More precisely, we have
\bea
\label{pr2}
T^{1/2}\left[J_T^h - \E(J_T^h)\right] \ConvD \eta\quad {\rm as}\quad T\to\infty,
\eea
where the symbol $\ConvD$ stands for convergence in distribution, and
$\eta$ is a normally distributed random variable with mean zero and variance
$\si^2_h(J)$ given by
\beq
\label{tsigma-M}
\sigma^2_h(J)=4\pi e(h)\int_{\mathbb{R}} f^2(\lambda) g^2(\lambda) d\lambda
+\kappa_4 e(h)\left[\int_{\mathbb{R}} f(\lambda) g(\lambda) d\lambda\right]^2.
\eeq
Here $\kappa_4$ is the fourth cumulant of $\xi(1)$, and
\beq
\label{t9}
e(h):=\frac{H_4}{H_2^2}= \int_0^1h^4(t)dt \left(\int_0^1h^2(t)dt\right)^{-2}.
\eeq
\end{thm}

Taking into account the equality
\bea
\label{pt1}
T^{1/2}\left[J_T^h - J\right]=T^{1/2}\left[\E(J_T^h)-  J\right]
+ T^{1/2}\left[J_T^h - \E(J_T^h)\right],
\eea
as an immediate consequence of Theorems \ref{T-Bias} and \ref{T-CLT},
we obtain the next result that contains sufficient conditions for a
simple "plug-in" statistic $J(I^h_{T})$ to be an asymptotically normal
estimator for a linear spectral functional $J(f)$.

\begin{thm}
\label{TT1}
Let the functionals $J:=J(f)$ and $J_T^h: = J(I^h_{T})$ be defined
by (\ref {t1}) and (\ref{t7}), respectively.
Then under Assumptions \ref{(A1)}, \ref{(A2)} (or \ref{(A3)}) and \ref{(H)}
the statistic $J_T^h$ is an asymptotically normal estimator for functional $J$.
More precisely, we have
\bea
\label{t8}
T^{1/2}\left[J_T^h-  J\right] \ConvD \eta\quad {\rm as}\quad T\to\infty,
\eea
where $\eta$ is as in Theorem \ref{T-CLT}, that is, $\eta$ is a normally
distributed random variable with mean zero and variance $\si^2_h(J)$
given by \eqref{tsigma-M} and \eqref{t9}.
\end{thm}

In the Gaussian case we have more accurate result. The next theorem,
which was proved in Ginovyan and Sahakyan \cite{GS2019a}, states that
for Gaussian models Assumptions \ref{(A3)} and \ref{(H)} 
are sufficient for the statistic
$J_T^h$ to be an asymptotically normal estimator for functional $J$.

\begin{thm}
\label{GTT1}
Let the functionals $J:=J(f)$ and $J_T^h: = J(I^h_{T})$ be defined
by (\ref {t1}) and (\ref{t7}), respectively.
Then under Assumptions \ref{(A3)} and \ref{(H)} 
the statistic
$J_T^h$ is an asymptotically normal estimator for functional $J$.
More precisely, we have
\bea
\label{gt8}
T^{1/2}\left[J_T^h-  J\right] \ConvD \eta\quad {\rm as}\quad T\to\infty,
\eea
where $\eta$ is a normally distributed random variable with mean zero and variance
$\si^2_h(J)$ given by
\beq
\label{tsigmaG}
\sigma^2_h(J)=4\pi e(h)\int_{\mathbb{R}} f^2(\lambda) g^2(\lambda) d\lambda,
\q e(h):=H_4H_2^{-2},
\eeq
and $H_k$ is as in \eqref{t4}.
\end{thm}

\sn{Parametric estimation: The Whittle procedure}
\label{AppP}
We assume here that the spectral density $f(\la)$ belongs to a given parametric family
of spectral densities $\mathcal{F}:=\{f(\lambda,\theta): \, \theta\in\Theta\}$,
where $\theta:=(\theta_1, \ldots, \theta_p)$ is an unknown parameter and
$\Theta$ is an open subset in the Euclidean space $\mathbb{R}^p$.
The problem of interest is to estimate the unknown parameter $\theta$ on the basis
of the tapered data \eqref{t2}, 
and investigate the asymptotic (as $T \to \infty$) properties of the suggested estimators,
depending on the dependence (memory) structure of the model $X(t)$ and the
smoothness of its spectral density $f$.

There are different methods of estimation: maximum likelihood,
Whittle, minimum contrast, etc.
Here we focus on the Whittle method.

The Whittle estimation procedure, originally devised for d.t.\  short memory
stationary processes, is based on the smoothed periodogram analysis on a frequency
domain, involving approximation of the likelihood function and asymptotic properties
of empirical spectral functionals (see Whittle \cite{W}).
The Whittle estimation method since its discovery has played
a major role in the asymptotic theory of parametric estimation in the frequency domain,
and was the focus of interest of many statisticians.
Their aim was to weaken the conditions needed to guarantee the validity of the
Whittle approximation for d.t.\  short memory models,
to find analogues for long and intermediate memory models,
to find conditions under which the Whittle estimator is asymptotically equivalent
to the exact maximum likelihood estimator, and to extend the procedure to the
c.t.\  models and random fields.

For the d.t.\  case, it was shown that for Gaussian and linear
stationary models the Whittle approach leads to consistent and asymptotically
normal estimators under short, intermediate and long memory assumptions.
Moreover, it was shown that in the Gaussian case the Whittle
estimator is also asymptotically efficient in the sense of Fisher
(see, e. g., Dahlhaus \cite{Dah1}, Dzhaparidze \cite{Dz1}, Fox and Taqqu \cite{FT1},
Giraitis and Surgailis \cite{GSu},
Guyon \cite{Gu}, Heyde and Gay \cite{HG}, Taniguchi and Kakizawa \cite{TK},
Walker \cite{Wal}, and references therein).

For c.t.\  models, the Whittle estimation procedure has been
considered, for example, in Anh et al. \cite{ALS2}, Avram et al. \cite{AvLS},
Casas and Gao \cite{CG}, Dzhaparidze \cite{Dz1}, Gao \cite{Go},
Gao et al. \cite{GAHT}, Leonenko and Sakhno \cite{Leonenko:2006},
Tsai and Chan \cite{TC},
where can also be found additional references. In this case, it was proved that
the Whittle estimator is consistent and asymptotically normal.

The Whittle estimation procedure based on the d.t.\  tapered data
has been studied in Alomari et al. \cite{ALRST}, Dahlhaus \cite{D1},
Dahlhaus and K\"unsch \cite{DK}, Guyon \cite{Gu}, Lude\~na and Lavielle \cite{LL}.
In the case where the underlying model is a L\'evy-driven c.t.\
linear process with possibly unbounded or vanishing spectral density
function, consistency and asymptotic normality of the Whittle estimator
was established in Ginovyan \cite{G2020e}.

To explain the idea behind the Whittle estimation procedure, assume for simplicity
that the underlying process $X(t)$ is a d.t.\  Gaussian process,
and we want to estimate the parameter $\theta$ based on the sample
$X_T:=\{X(t), \, t=1,\ldots, T\}$. A natural approach is to find the maximum likelihood
estimator (MLE) $\widehat\theta_{T,MLE}$ of $\theta$, that is, to maximize the
log-likelihood function $L_T(\theta)$, which in this case takes the form:
$$L_T(\theta)=-\frac T2\ln2\pi-\frac12\ln\det B_T(f_\theta)-\frac12X'_T[B_T(f_\theta)]^{-1}X_T,$$
where $B_T(f_\theta)$  is the Toeplitz matrix generated by $f_\theta$.
Unfortunately, the above function is difficult to handle, and no explicit expression
for the estimator $\widehat\theta_{T,MLE}$ is known (even in the case of simple models).
An approach, suggested by P. Whittle, called the Whittle estimation procedure,
is to approximate the term $\ln\det B_T(f_\theta)$
by $\frac T2\int_{-\pi}^\pi\ln f_\theta(\la)d\la$ and the inverse matrix $[B_T(f_\theta)]^{-1}$
by the Toeplitz matrix $B_T(1/f_\theta)$.
This leads to the following approximation of the log-likelihood function $L_T(\theta)$,
introduced by P. Whittle \cite{W}, and called Whittle functional:
$$L_{T,W}(\theta)=-\frac 1{4\pi}\int_{-\pi}^\pi\left[\ln f_\theta(\la)
+\frac{I_{T}(\la)}{f_\theta(\la)}\right]\, d\la,$$
where $I_{T}(\lambda)$ is the ordinary periodogram of the process $X(t)$.

Now maximizing the Whittle functional $L_{T,W}(\theta)$ with respect to
$\theta$, we get the Whittle estimator $\widehat\theta_{T}$ for $\theta$.
It can be shown that if
$$ T^{-1/2}(L_T(\theta)- L_{T,W}(\theta)\to0 \q {\rm as}\q n\to\f \q {\rm in \,\, probability,}
$$
then the MLE $\widehat\theta_{T,MLE}$ and the Whittle estimator $\widehat\theta_{T}$
are asymptotically equivalent in the sense that $\widehat\theta_{T}$ also is
consistent, asymptotically normal and asymptotically Fisher-efficient (see, e.g., Dzhaparidze \cite{Dz1}).

In the continuous context, the Whittle procedure of estimation of a spectral
parameter $\theta$ based on the sample $X_T:=\{X(t), \, 0\leq t\leq T\}$
is to choose the estimator $\widehat\theta_{T}$ to minimize the weighted
Whittle functional:
\beq
\label{pe9a}
U_{T}(\theta): =\frac1{4\pi}\int_{\mathbb{R}}\left[\ln f(\la, \theta) +
\frac{I_{T}(\la)}{f(\la, \theta)}\right]\cd w(\la) \, d\la,
\eeq
where $I_{T}(\la)$ is the continuous periodogram of $X(t)$,
and $w(\la)$ is a weight function ($w(-\la)=w(\la)$, $w(\la)\ge0$,
$w(\la)\in L^1(\mathbb{R})$) for which the integral in (\ref{pe9a})
is well defined.
An example of common used weight function is $w(\la)=1/(1+\la^2)$.

The Whittle procedure of estimation of a spectral parameter $\theta$
based on the tapered sample \eqref{t2} is to choose the estimator
$\widehat\theta_{T,h}$ to minimize the weighted tapered Whittle functional:
\beq
\label{pe9}
U_{T,h}(\theta): =\frac1{4\pi}\int_{\Lambda}\left[\log f(\la, \theta) +
\frac{I^h_{T}(\la)}{f(\la, \theta)}\right]\cd w(\la) \, d\la,
\eeq
where $I^h_{T}(\la)$ is the tapered periodogram of $X(t)$, given by \eqref{t6},
and $w(\la)$ is a weight function for which the integral in (\ref{pe9})
is well defined. Thus,
\beq
\label{pee9}
\widehat\theta_{T,h}: =\underset{\theta\in\Theta}{\rm Arg \, min} \, U_{T,h}(\theta).
\eeq

Here we follow the paper Ginovyan \cite{G2020e}.
To state results involving properties of the Whittle estimator,
we first introduce the following set of assumptions.

\begin{asn}
\label{(B1)}
\rm{The true value $\theta_0$ of the parameter $\theta$ belongs to a compact set
$\Theta$, which is contained in an open set $S$ in the $p$-dimensional
Euclidean space $\mathbb{R}^p$,
and $f(\la,\theta_1)\neq f(\la,\theta_2)$ whenever $\theta_1\neq \theta_2$
almost everywhere in $\Lambda$ with respect to the Lebesgue measure.}
\end{asn}

\begin{asn}
\label{(B2)}
{\rm The functions $f(\la,\theta)$, $f^{-1}(\la,\theta)$ and
$\frac{\partial}{\partial\theta_k}f^{-1}(\la,\theta)$, $k=1,\ldots,p$,
are continuous in $(\la,\theta)$.}
\end{asn}

\begin{asn}
\label{(B3)}
{\rm The functions $f:=f(\la,\theta)$ and
$g:=w(\la)\frac{\partial}{\partial\theta_k}f^{-1}(\la,\theta)$
satisfy Assumptions \ref{(A2)} or \ref{(A3)} 
for all $k=1,\ldots,p$ and $\theta\in\Theta$.}
\end{asn}

\begin{asn}
\label{(B4)}
{\rm The functions $a:=a(\la,\theta)$ and $b:=\widehat g$, where $g$ is as in
Assumption \ref{(B3)}, satisfy Assumption \ref{(A1)}.} 
\end{asn}

\begin{asn}
\label{(B5)}
{\rm The functions
$\frac{\partial^2}{\partial\theta_k\partial\theta_j}f^{-1}(\la,\theta)$ and
$\frac{\partial^3}{\partial\theta_k\partial\theta_j\partial\theta_j}f^{-1}(\la,\theta)$,
$k,j,l=1,\ldots, p$, are continuous in $(\la,\theta)$ for $\la\in\Lambda$,
$\theta\in N_\de(\theta_0)$, where $N_\de(\theta_0):=\{\theta: \, |\theta-\theta_0|<\de\}$
is some neighborhood of $\theta_0$.}
\end{asn}

\begin{asn}
\label{(B6)}
{\rm The matrices
\bea
\label{W03}
W(\theta):=\|w_{ij}(\theta)\|, \,\,
A(\theta):=\|a_{ij}(\theta)\|,\,\,
B(\theta):=\|b_{ij}(\theta)\| \q i,j=1,\ldots,p
\eea
are positive definite, where
\bea
\label{W3}
w_{ij}(\theta)&=&  \frac1{4\pi}\int_{\Lambda}
\frac{\partial}{\partial\theta_i}\ln f(\la, \theta)
\frac{\partial}{\partial\theta_j}\ln f(\la, \theta)w(\la)d\la,\\
\label{W4}
a_{ij}(\theta)&=& \frac1{4\pi}\int_{\Lambda}
\frac{\partial}{\partial\theta_i}\ln f(\la, \theta)
\frac{\partial}{\partial\theta_j}\ln f(\la, \theta)w^2(\la)d\la,\\
\label{W5}
b_{ij}(\theta)&=& \frac{\kappa_4}{16\pi^2}\int_{\Lambda}
\frac{\partial}{\partial\theta_i}\ln f(\la, \theta)w(\la)d\la
\int_{\mathbb{R}} \frac{\partial}{\partial\theta_j}\ln f(\la, \theta)w(\la)d\la,
\eea
and $\kappa_4$ is the fourth cumulant of $\xi(1)$.}
\end{asn}

\n{\it Consistency of the Whittle estimator.} The next theorem contains sufficient
conditions for Whittle estimator to be consistent.
\begin{thm}
\label{CWT}
Let $\widehat\theta_{T,h}$ be the Whittle estimator defined by \eqref{pee9} and let
$\theta_0$ be the true value of parameter $\theta$. Then, under Assumptions
\ref{(B1)}--\ref{(B4)} and \ref{(H)}, 
the statistic $\widehat\theta_{T,h}$ is a consistent estimator
for $\theta$, that is, $\widehat\theta_{T,h}\to \theta_0$ in probability as $T\to\f$.
\end{thm}

\n
{\it Asymptotic normality of the Whittle estimator.}
Having established the consistency of the Whittle estimator $\widehat\theta_{T,h}$,
we can go on to obtain the limiting distribution of $T^{1/2}\left(\widehat\theta_{T,h}-\theta_0\right)$
in the usual way by applying the Taylor's formula, the mean value theorem, and Slutsky's arguments.
Specifically we have the following result, showing that under the above assumptions,
the Whittle estimator $\widehat\theta_{T,h}$ is asymptotically normal.

\begin{thm}
\label{TAN}
Suppose that Assumptions \ref{(B1)}--\ref{(B6)} and \ref{(H)}
are satisfied.
Then the Whittle estimator $\widehat\theta_{T,h}$ of an unknown spectral parameter
$\theta$ based on the tapered data \eqref{t2} is asymptotically normal. More precisely,
we have
\bea
\label{aW1}
T^{1/2}\left(\widehat\theta_{T,h}-\theta_0\right)\ConvD N_p\left(0, e(h)\G(\theta_0)\right)
\quad {\rm as}\quad T\to\infty,
\eea
where $N_p(\cdot,\cdot)$ denotes the $p$-dimensional normal law, \, $\ConvD$
stands for convergence in distribution,
\bea
\label{aW2}
\G(\theta_0) = W^{-1}(\theta_0)\left(A(\theta_0)+B(\theta_0)\right)W^{-1}(\theta_0),
\eea
where the matrices $W$, $A$ and $B$ are defined in (\ref{W03})-(\ref{W5}),
and the tapering factor $e(h)$ is given by formula \eqref{t9}.
\end{thm}

\section{Methods and tools}
\label{methods}

In this section we briefly discuss the methods and tools, used to prove
the central and noncentral limit theorems for Toeplitz type quadratic
forms and functionals stated in Sections \ref{G-CLT}--\ref{Tap1},
as well as the results stated in Section \ref{App}.

As mentioned in the introduction, the most commonly used methods to prove central
limit theorems are: 
(a) the method of characteristic functions.
(b) the method of cumulants or moments,
(c) the approximation method.

To prove the central limit theorems stated in this paper, depending on the underlying model,
was used either the method of cumulants or the approximation method.
If the underlying model is Gaussian, then the method of cumulants was applied,
otherwise the approximation method was applied, which reduces the
quadratic integral form to a single integral form.
To prove the non-central limit theorems, was used the spectral representation of the
underlying process, the properties of L\'evy-It\^o-type and Stratonovich-type multiple
stochastic integrals, and power counting theorems.
Some details of the above methods are described below.

\subsection{The characteristic functions and cumulant criteria for the CLT}
\label{cm}
The characteristic functions criterion for the CLT is based on the fact that convergence
in distribution is equivalent to the pointwise convergence of the corresponding
characteristic functions.
The general cumulant criterion for the CLT is based on the following result
(see, e.g., Giraitis et al. \cite{GKSu}, Corollary 4.2.1).
\begin{pp}
\label{gIb}
Let the random variables $\{X_T, T\in \mathbb{R}\}$ have all moments finite,
and let $\E[X_T]\to0$, $\Var[X_T]\to\si^2$, and $\chi_k(X_T) :=\Cum_k(X_T)\to0$
for all $k=3,4,\ldots$ as $T\to\f$. Then $X_T\ConvD X \sim N(0,\sigma^2)$ as $T\to\f$.
\end{pp}

The characteristic functions and cumulant criteria for the CLT for quadratic functionals is based on
the following general result (see Ibragimov \cite{I1963}).

Let $\xi$ be a Gaussian random variable with values in a separable Hilbert
space $H$. In other words, $\xi$ is a random variable with characteristic
functional:
\begin{equation}
\I(h) = {\exp}\left\{m(h) -\frac 12 (Rh,h)\right\}, \q h \in H,
 \end{equation}
where $m(h)$ is a continuous linear functional and the correlation operator
$R$ is a self-adjoint completely continuous operator with finite trace.
We assume, without loss of generality, that $m(h)\equiv0$.
Let $A$ be some linear self-adjoint bounded operator.

The proof of the following result can be found in Ibragimov \cite{I1963}.
\begin{pp}
\label{Ib}
Let the operators $R$ and $A$ be as above.
The quadratic form $(A\xi,\xi)$ has the same distribution as the sum
$\sum_{k = 1}^\f \la_k^2\xi_k^2$, where $\xi_k$ are independent $N(0,1)$
Gaussian random variables and $\la_k$ are the eigen-values of the operator
$B:= RA.$
\end{pp}
\begin{rem}
{\rm It can easily be shown that the sets of non-zero eigen-values of the operators
$RA$, $AR$ and $R^{1/2}AR^{1/2}$ coincide, where $R^{1/2}$ is the positive
definite square root of $R$.}
\end{rem}

As mentioned above, Toeplitz matrices and operators arise naturally
in the theory of stationary processes, and serve as tools, to study many
topics of the spectral and statistical analysis of d.t.\  and c.t.\
stationary processes.

Let $A_T(f)$ denote the covariance matrix (or operator) of the process
$\{X(u), \ u\in \mathbb{U}\}$, that is, $A_T(f)$ is either the $T\times T$
Toeplitz matrix ($B_T(f)$), or the $T$-truncated Toeplitz operator $W_T(f)$
generated by the spectral density $f$, and let $A_T(g)$ denote either the
$T\times T$ Toeplitz matrix, or the $T$-truncated Toeplitz operator generated
by the function $g$ (for definitions see formulas \eqref{IMT2-1} and \eqref{IMT3-1}).

As a consequence of Proposition \ref{Ib},
we have the following result, which gives a link between the distribution
of the quadratic functional $Q_T$ in (\ref{MTc-1}) and the trace problem
for Toeplitz matrices and operators (see, e.g., Ginovyan et al. \cite{GST2014},
Grenander and Szeg\H{o} \cite{GS}, Ibragimov \cite{I1963}).

\begin{itemize}
\item[1.]
The quadratic functional $Q_T$ in (\ref{MTc-1}) has the same
distribution as the sum $\sum_{k = 1}^\f \la_{k,T}^2\xi_k^2$
($\sum_{k = 1}^T \la_{k,T}^2\xi_k^2$ \ in the d.t.\  case),
where $\{\xi_k, k\ge1\}$ are independent $N(0,1)$ Gaussian random variables
and $\{\la_{k,T}, k\ge1\}$ are the eigenvalues of the operator $A_T(f)A_T(g)$.

\item[2.]
The characteristic function $\I_T (t)$ of $Q_T$ is given by
\begin{equation}
\label{MTc-055}
\I_T (t) = \prod_{k = 1}^\f|1 - 2it\la_{k,T}|^{-1/2}.
\end{equation}

\item[3.]
The $k$--th order cumulant $\chi_k(\cdot)$  of $Q_T$ is given by
\begin{eqnarray}
\label{MTc-056}
\chi_k(Q_T) = 2^{k-1}(k-1)! \sum_{j = 1}^\f \la_{j,T}^k
=2^{k-1}(k-1)!\, \tr\,[A_T(f)A_T(g)]^k,
\end{eqnarray}
where $\tr[A]$ stands for the trace of an operator $A$.
\end{itemize}
The product in \eqref{MTc-055} and the sum in \eqref{MTc-056} are over $j=1,\ldots,T$
in the d.t.\  case.

\paragraph{The tapered case.}
To study the asymptotic distribution (as $T\to\f$) of the tapered functional $Q^h_T$,
given by (\ref{tq-1}), we use the method of cumulants, the frequency-domain approach,
and the technique of truncated tapered Toeplitz matrices and operators.

Let $A_T^h(f)$ be either the $T\times T$ tapered Toeplitz matrix $B_T^h(f)$,
or the $T$-truncated tapered Toeplitz operator $W_T^h(f)$
generated by the spectral density $f$ and taper $h$, and let $A_T^h(g)$
denote either the $T\times T$ tapered Toeplitz matrix, or the $T$-truncated tapered
Toeplitz operator generated by the functions $g$ and $h$
(for definitions see formulas \eqref{1-7M} and \eqref{1-7}).
Similar to the non-tapered case, we have the following results
(cf. Ginovyan et al. \cite{GST2014}, Grenander and Szeg\H{o} \cite{GS}, Ibragimov \cite{I1963}).

\begin{itemize}
\item[1.]
The quadratic functional $Q^h_T$ in (\ref{tq-1}) has the same distribution as the sum
$\sum_{j = 1}^\f\la_{j,T}^2\xi_j^2$, where \mbox{$\{\xi_j, j\ge1\}$} are independent
$N(0,1)$ Gaussian random variables and $\{\la_{j,T}, j\ge1\}$ are the eigenvalues of
the operator $A_T^h(f)\,A_T^h(g)$.

\item[2.]
The characteristic function $\I(t)$ of $Q^h_T$ is given by formula:
\begin{equation}
\label{MTc-05}
\I(t) = \prod_{j = 1}^\f|1 - 2it\la_{j,T}|^{-1/2}.
\end{equation}

\item[3.]
The $k$--th order cumulant $\chi_k(Q^h_T)$ of $Q^h_T$ is given by formula:
\begin{eqnarray}
\label{MTc-5}
\chi_k(Q_T) = 2^{k-1}(k-1)! \sum_{j = 1}^\f\la_{j,T}^k
=2^{k-1}(k-1)!\, \tr\,[A_T^h(f)\,A_T^h(g)]^k.
\end{eqnarray}
\end{itemize}

Thus, to describe the asymptotic distribution of the quadratic functional $Q^h_T$,
we have to control the traces and eigenvalues of the products of truncated
tapered Toeplitz operators and matrices.

\subsection{Approximation of traces of products of Toeplitz matrices and operators}

The trace approximation problem for truncated Toeplitz operators and matrices
has been discussed in detail in the survey paper Ginovyan et al. \cite{GST2014}
for non-tapered case.
Here we present some important results, which were used to prove CLT for
quadratic functionals $Q^h_T$.

Let $h$ be a taper function satisfying Assumption \ref{(H)}. 
Let $A^h_T(\psi)$ be either the $T\times T$ tapered Toeplitz matrix $B^h_T(\psi)$,
or the $T$-truncated tapered Toeplitz operator $W^h_T(\psi)$
generated by a function $\psi$ (for definitions see formulas
\eqref{1-7M} and \eqref{1-7}) 

Observe that, in view of  \eqref{t4},  \eqref{t55}, \eqref{1-7M} and \eqref{1-7}, we have
\begin{equation}
\label{MT-00}
\frac1{T} \tr\left[A^h_T(\psi)\right]=\frac1T\cd \widehat \psi(0)\cd\int_0^Th_T^2(t)dt
= 2\pi H_2\int_{\Lambda} \psi(\la)d\la.
\end{equation}
What happens to the relation (\ref{MT-00}) when $A^h_T(\psi)$
is replaced by a product of Toeplitz matrices (or operators)?
Observe that the product of Toeplitz matrices (resp. operators)
is not a Toeplitz matrix (resp. operator).

The idea is to approximate the trace of the product of
Toeplitz matrices (resp. operators) by the trace of a Toeplitz matrix
(resp. operator) generated by the product of the corresponding generating functions.
More precisely, let $\{\psi_1,\psi_2,\ldots,\psi_m\}$ be a
collection of integrable real symmetric functions defined on $\Lambda$.
Let $A_T^h(\psi_i)$ be either the $T\times T$ tapered Toeplitz matrix $B_T^h(\psi_i)$,
or the $T$-truncated tapered Toeplitz operator $W_T^h(\psi_i)$
generated by a function $\psi_i$ and a taper function $h$. Define
\begin{eqnarray}
\label{n 4-4}
\nonumber
S_{A,\mathcal{H},h}(T):=\frac1T\tr\left[\prod_{i=1}^m A_T^h(\psi_i)\right
],\q M_{\Lambda,\mathcal{H},h}:=(2\pi)^{m-1}H_m\int_{\Lambda}
\left[\prod_{i=1}^m \psi_i(\la)\right]\,d\la,
\end{eqnarray}
and let
\begin{eqnarray}
\label{n 4-5}
&&\De(T):=\De_{A,{\Lambda},\mathcal{H},h}(T)=|S_{A,\mathcal{H},h}(T)-
M_{{\Lambda},\mathcal{H},h}|.
\end{eqnarray}

\begin{pp}
\label{T5}
Let $\De(T)$ be as in (\ref{n 4-5}). Each of the following conditions is sufficient for
\begin{equation}
\label{n 4-7}
\De(T)=o(1) \q{\rm as} \q T\to\f.
\end{equation}
\begin{itemize}
\item[{\bf(C1)}]
$\psi_i\in L^1(\Lambda)\cap L^{p_i}(\Lambda)$, $p_i>1$,
$i=1,2,\ldots,m$,
with $1/p_1+\cdots+1/p_m\le1$.

\item[{\bf(C2)}]
The function $\varphi({\bf u})$ defined by
%
\begin{equation}
\label{n 4-6}
\varphi({\bf u}):\,=\,
\int_\Lambda \psi_1(\la)\psi_2(\la-u_1)\psi_3(\la-u_2)\cdots \psi_m(\la-u_{m-1})\,d\la,
\end{equation}
where ${\bf u}=(u_1,u_2,\ldots,u_{m-1})\in\Lambda^{m-1}$,
belongs to $L^{m-2}(\Lambda^{m-1})$ and is continuous at
${\bf0}=(0,0,\ldots,0)\in\Lambda^{m-1}$.
\end{itemize}
\end{pp}

\begin{rem}
\label{rem4-1}
{\rm In the nontapered case, Proposition \ref{T5} was proved in Ginovyan et al. \cite{GST2014},
in the tapered case, it was proved in Ginovyan \cite{G2020g}.
Proposition \ref{T5} was used to prove parts (B) and (C) of Theorem \ref{cth1}
(in the nontapered case) and Theorems \ref{th3} and \ref{th1} (in the tapered case).
In the special case $m=4$, $\psi_1=\psi_3:=f$ and $\psi_2=\psi_4:=g$,
in Ginovyan and Sahakyan \cite{GS2005} (in the d.t.\  case) and
in Ginovyan and Sahakyan \cite{GS2007} (in the c.t.\  case)
it was proved that the conditions of Theorem \ref{cth1}(D) and Theorem \ref{cth2}
are also sufficient for (\ref{n 4-7})}.
\end{rem}
\begin{rem}
{\rm More results concerning the trace approximation problem for truncated Toeplitz
operators and matrices can be found in Ginovyan and Sahakyan \cite{GS2012,GS2013},
Ginovyan et al. \cite{GST2014} and Lieberman and Phillips \cite{LP}.}
\end{rem}

\subsection{Approximation method for the CLT}
\label{am}

The approximation method for CLT for quadratic functionals is based on
approximation of quadratic functional $Q_T$ by a univariate sum (integral)
of $m$-dependent random variables, and then using the CLT for such variables
(see, e.g., Giraitis et al. \cite{GKSu}, Section 4.5, and  Giraitis and Surgailis \cite{GSu}).

Let $X_i(u)$ ($i=1,2$) be two linear processes of the form \eqref{dlp},
subordinated to the same orthonormal sequence $\{\xi(k),\,k\in\mathbb{Z}\}$,
with square summable covariance functions $r_i(u)$ ($i=1,2$), that is,
\begin{equation}
\label{dlp-2}
X_i(u) = \sum_{k=-\f}^{\f}a_i(u-k)\xi(k), \q \q
         \sum_{k=-\f}^{\f}|a_i(k)|^2 < \f, \q i=1,2,
\end{equation}
and
\begin{equation}
\label{dlp-3}
\sum_{u=-\f}^{\f}|r_i(u)|^2 < \f, \q i=1,2.
\end{equation}
Denote
\begin{equation}
\label{dlp-4}
S_T: = \sum_{u=1}^{T}X_1(u)X_2(u).
\end{equation}
\begin{pp}
\label{GSu-L}
Let $X_i(u)$ $(i=1,2)$ and $S_T$ be as above, and let the quadratic form
$Q_T$ be as in \eqref{MTc-1}. Then the following assertions hold.
\begin{itemize}
\item[(a)]
The distribution of  $T^{-1/2}(S_T-\E[S_T])$ tends to the
centered normal distribution with variance:
\begin{equation}
\label{dlp-5}
\si^2:=\sum_{u=-\f}^{\f}r_1(u)r_2(u) +\kappa_4r_{1,2}^2 ,
\end{equation}
where $\kappa_4$ is the fourth cumulant of $\xi(0)$, and $r_{1,2}=\E[X_1(0)X_2(0)]$.\\
\item[(b)] $\Var(Q_T-S_T)=o(T)$ as $T\to\f$.
\end{itemize}
\end{pp}
\n
A similar result is true for c.t.\  linear processes of the form \eqref{clp},
where now we have (see \cite{GS2007}):
$$S_T: = \int_0^{T}X_1(u)X_2(u)du.$$

\subsection{Fej\'er-type singular integrals}

We define Fej\'er-type kernels and singular integrals, and state some of their
properties that were used to prove the limit theorems stated in Section \ref{App}.

For a number $k$ ($k=2,3,\ldots$) and a taper function $h$ satisfying Assumption \ref{(H)}
consider the following Fej\'er-type 'tapered' kernel function:
\beq \label{kerN}
\Phi^h_{k,T}(\uu):=\frac{H_T(\uu)}{(2\pi)^{k-1}H_{k,T}(0)}, \q
\uu =(u_1, \ldots, u_{k-1})\in \mathbb{R}^{k-1},
\eeq
where
\beq \label{ker1N}
H_T(\uu):={H_{1,T}(u_1)\cdots H_{1,T}(u_{k-1}) H_{1,T}\left(-\sum_{j=1}^{k-1}u_j\right)},
\eeq
and the function $H_{k,T}(\cdot)$ is defined by (\ref{t3}) with $H_{k,T}(0)=T\cdot H_k\neq0$ (see (\ref{t4})).

The proofs of propositions that follow can be found in Ginovyan and Sahakyan \cite{GS2019}.
The next result shows that, similar to the classical Fej\'er kernel,
the 'tapered' kernel $\Phi^h_{k,T}(\uu)$ is an approximation identity
(see Ginovyan and Sahakyan \cite{GS2019}, Lemma 3.4).
\begin{pp}
\label{L3}
For any $k=2,3,\ldots$ and a taper function $h$ satisfying Assumption \ref{(H)}
the kernel $\Phi^h_{k,T}(\uu)$, $\uu =(u_1, \ldots, u_{k-1})\in {\mathbb R }^{k-1}$,
possesses the following properties:
\begin{itemize}
\item [a)]
$\sup_{T>0}\int_{{\mathbb R}^{k-1}}\left |\Phi^h_{k,T}(\uu)\right|\,d\uu =C_1<\f;$
\item [b)] $\int_{{\mathbb R}^{k-1}}\Phi^h_{k,T}(\uu)\,d\uu =1;$
\item [c)] $\lim_{T\to\f}\int_{{\mathbb E}^c_\de}\left|\Phi^h_{k,T}(\uu)\right|\,d\uu=0$
 for any $\de>0;$
\item [d)] If $k>2$ for any $\de>0$ there exists a constant $M_\de>0$ such that for $T>0$
\beq\label{z0}
\left\|\Phi^h_{k,T}\right\|_{L^{p_k}({\mathbb{E}^c_\de})}\le M_\de,
\eeq
where  $p_k=\frac {k-2}{k-3}$ for $k>3$,  $p_3=\infty $ and
$$
{\mathbb E}_\de^c={\mathbb R}^{k-1}\sm{\mathbb E}_\de,\quad
{\mathbb E}_\de=\{\uu =(u_1, \ldots, u_{k-1})\in {\mathbb R }^{k-1}:
\, |u_i|\le\de, \, i=1,\ldots,k-1\}.
$$
\item [e)]
If the function $\Psi\in L^1(\mathbb{R}^{k-1})\bigcap L^{k-2}(\mathbb{R}^{k-1})$
is continuous at $\vv =(v_1, \ldots, v_{k-1})$ $\ (L^0$ is the space of measurable functions), then
\beq\label{2-3-1}
\lim_{T\to\infty} \int_{\mathbb{R}^{k-1}}\Psi(\uu+\vv)\Phi^h_{k,T}(\uu)d\uu=\Psi(\vv).
\eeq
\end{itemize}
\end{pp}

Denote
\beq\label{zz}
\Delta_{2,T}^h : = \int_{\mathbb{R}^2}f(\la)g(\la+\mu)\Phi^h_{2,T}(\mu)d\la d\mu
-\int_{\mathbb{R}}f(\la)g(\la)d\la,
\eeq
where $\Phi^h_{2,T}(\mu)$ is given by (\ref{kerN}), (\ref{ker1N}).

The next two propositions, which were used to prove Theorems \ref{T-Bias} and \ref{TT1},
give information on the rate of convergence to zero of $\Delta_{2,T}^h$ as $T\to\f$.
\begin{pp}\label{L011}
Assume that Assumptions \ref{(A2)} and \ref{(H)} are satisfied.
Then the following asymptotic relation holds:
\beq\label{z111}
\Delta_{2,T}^h=  o\left(T^{-1/2}\right)\q {\rm as} \q T\to\f.
\eeq
\end{pp}
\begin{pp}
\label{L01}
Assume that Assumptions \ref{(A3)} and \ref{(H)}are satisfied.
Then the following inequality holds:
\beq\label{mz1}
|\Delta_{2,T}^h|\leq C_h\begin{cases}
T^{-(\be_1+\be_2)},&\text{if}\  \ \be_1+\be_2<1\\
T^{-1}\ln T,&\text{if}\  \ \be_1+\be_2=1\\
T^{-1},&\text{if} \ \ \be_1+\be_2>1,
\end{cases}\qquad T>0,
\eeq
where $C_h$ is a constant depending on $h$.
\end{pp}
Notice that for non-tapered case ($h(t)={\mathbb I}_{[0,1]}(t)$), the above stated
results were proved in Ginovyan and Sahakyan \cite{GS2007}
(see also Ginovyan and Sahakyan \cite{GS2012,GS2013}).

\subsection{L\'evy-It\^o-type and Stratonovich-type multiple stochastic integrals}
\label{sec:prelim}

To prove limit theorems for quadratic functionals of L\'evy-driven
c.t.\  linear models, was used the multiple \emph{off-diagonal} (It\^o-type)
and \emph{with-diagonal} (Stratonovich-type) stochastic integrals with respect
to L\'evy noise. In this subsection we introduce these integrals, and briefly
discuss their properties (see, e.g., Bai et al. \cite{BGT2}, Farr\'e et al. \cite{FJ2010},
Peccati and Taqqu \cite{PT2011}).

Let $f$ be a function in $L^2(\mathbb{R}^k)$, then the following off-diagonal
multiple stochastic integral, called It\^o-L\'evy integral, is well-defined:
\begin{equation}\label{eq:I_k^xi}
I_k^\xi(f)=\int_{\mathbb{R}^k}' f(x_1,\ldots,x_k) \xi(dx_1) \ldots\xi (dx_k),
\end{equation}
where $\xi(t)$ is a L\'evy process with $\E\xi(t)=0$ and
$\Var[\xi(t)]=\sigma_\xi^{2}t$, and the prime $'$ indicates that we do not
integrate on the diagonals $x_i=x_j$, $i\neq j$.
The multiple integral $I_k^\xi(\cdot)$ satisfies the following inequality:
\begin{equation}\label{eq:isometry ineq}
\|I_k^\xi(f)\|_{L^2(\Omega)}^2 \le k!\sigma_\xi^{2k} \|f\|_{L^2(\mathbb{R}^k)}^2,
\end{equation}
and the inequality in (\ref{eq:isometry ineq}) becomes equality if $f$ is symmetric:
\begin{equation}\label{eq:isometry}
\|I_k^\xi(f)\|_{L^2(\Omega)}^2 = k!\sigma_\xi^{2k} \|f\|_{L^2(\mathbb{R}^k)}^2.
\end{equation}
Observe that if in \eqref{eq:I_k^xi}, $\xi(\cdot)=B(\cdot)$, where $B(\cdot)$
is the real-valued Brownian motion, then the corresponding integral:
\begin{equation}\label{eq:I_k^B}
I_k^B(f)=\int_{\mathbb{R}^k}' f(x_1,\ldots,x_k) B(dx_1) \ldots B (dx_k)
\end{equation}
is called multiple Wiener-It\^o integral (see It\^o \cite{Ito1951}).

The Wiener-It\^o integral  can also be defined with respect to the
complex-valued Brownian motion:
\begin{equation}\label{eq:I_k^W}
I_k^W(g)=\int_{\mathbb{R}^k}'' g(u_1,\ldots,u_k) W(du_1) \ldots W(du_k),
\end{equation}
where  $g\in L^2(\mathbb{R}^k)$ is a complex-valued function satisfying
$g(-u_1,\ldots,-u_k)=\overline{g(u_1,\ldots,u_k)}$, and $W(\cdot)$ is a
complex-valued Brownian motion (with real and imaginary parts being independent)
viewed as a random integrator
(see, e.g., Embrechts and Maejima \cite{EM2002}),
and the double prime $''$ indicates the exclusion of the hyper-diagonals
$u_p=\pm u_q$, $p\neq q$.

The next result, which can be deduced from  Proposition 9.3.1 of Peccati and Taqqu
\cite{PT2011} and  Proposition 4.2 of Dobrushin \cite{Dob1979} (see Bai et al. \cite{BGT2}),
gives a relationship between the integrals $I_k^B(\cdot)$ and $I_k^W(\cdot)$, defined by
(\ref{eq:I_k^B}) and  (\ref{eq:I_k^W}), respectively.
\begin{pp}\label{Pro:Time<->Spec}
Let $f_j(\cdot)$ be  real-valued functions in $L^2(\mathbb{R}^{k_j})$,
$j=1,\ldots,J$, and let
$$
\widehat{f}_j(w_1,\ldots,w_{k_j})=\int_{\mathbb{R}^{k_j}} f_j(x_1,\ldots,x_{k_j})
e^{i\left(x_1w_1+\ldots+x_{k_j}w_{k_j}\right)} dx_1\ldots dx_{k_j}
$$
be the $L^2$-Fourier transform of $f_j(\cdot)$. Then
\begin{equation*}
\Big(I_{k_1}^B(f_1),\ldots,I_{k_J}^B(f_J)\Big)\overset{d}{=}\Big((2\pi)^{-k_1/2}I_{k_1}^W\left(\widehat{f}_1 A^{\otimes_{k_1}}\right),\ldots,(2\pi)^{-k_J/2}I_{k_J}^W\left(\widehat{f}_J A^{\otimes_{k_J}}\right)\Big),
\end{equation*}
for any function $A(u):\mathbb{R}\rightarrow \mathbb{C}$ such that
$|A(u)|=1$ and $A(w)=\overline{A(-w)}$  almost everywhere,
where $A^{\otimes k}(w_1,\ldots ,w_k):=A(w_1)\cdots A(w_k)$.
\end{pp}

In the next proposition we state a stochastic Fubini's theorem
(see Bai et al. \cite{BGT2},
Lemma 3.1, or Peccati and Taqqu \cite{PT2011}, Theorem 5.12.1).
\begin{pp}\label{Lem:fubini}
Let $(S,\mu)$ be a measure space with $\mu(S)<\infty$, and let
$f(s,x_1,\ldots,x_k)$ be a function on $S\times \mathbb{R}^k$ such that
\[
\int_{S}\int_{\mathbb{R}^k} f^2(s,x_1,\ldots,x_k)dx_1\ldots dx_k \mu(ds)<\infty,
\]
then we can change the order of  the multiple stochastic integration
$I_k^\xi(\cdot)$ and the deterministic integration $\int_S f(s,\cdot) \mu(ds)$:
\[
\int_{S}I_k^\xi\big(f(s,\cdot)\big)\mu(ds)= I_k^\xi\left(\int_S f(s,\cdot)\mu(ds)\right).
\]
\end{pp}

The \emph{with-diagonal} counterpart of the L\'evy-It\^o integral $I_k^{\xi}(f)$, called
a \emph{Strato\-novich-type} stochastic integral, is defined by
\begin{equation}\label{eq:stra int}
\mathring{I}_k^{\xi}(f):=\int_{\mathbb{R}^k} f(x_1,\ldots,x_k) \xi(dx_1) \ldots\xi (dx_k),
\end{equation}
which includes all the diagonals. We refer to Farr\'e et al. \cite{FJ2010}
for a comprehensive treatment of Stratonovich-type integrals $\mathring{I}_k^{\xi}(f)$.
Observe that for the with-diagonal integral $\mathring{I}_k^{\xi}(f)$ to be well-defined,
the integrand $f$ needs also to be square-integrable on all the diagonals of
$\mathbb{R}^k$ (see Bai et al. \cite{BGT2}, Farr\'e et al. \cite{FJ2010}).

The with-diagonal integral $\mathring{I}_k^{\xi}(f)$ can be expressed
by off-diagonal integrals of lower orders using the Hu-Meyer formula
(see Farr\'e et al. \cite{FJ2010}, Theorem 5.9).
In the special case when $k=2$, we have
\begin{equation}\label{eq:off diag decomp}
\mathring{I}_2^{\xi}(f)=\int_{\mathbb{R}^2}'f(x_1,x_2)\xi(dx_1)\xi(dx_2)
+\int_{\mathbb{R}}f(x,x)\xi^{(2)}_c(dx) + \int_{\mathbb{R}}f(x,x)dx,
\end{equation}
where
\begin{equation}\label{eq:xi^2_c(t)}
\xi^{(2)}_c(t)=\xi^{(2)}(t)-\E\xi^{(2)}(t) =\xi^{(2)}(t)-|t|,
\end{equation}
and $\xi^{(2)}(t)$ is the quadratic variation of $\xi(t)$, which is
non-deterministic if $\xi(t)$ is non-Gaussian (see Farr\'e et al. \cite{FJ2010},  equation (10)).
The centered process $\xi^{(2)}_c(t)$ is called a second order \emph{Teugels martingale},
which is a L\'evy process with the same filtration as $\xi(t)$, whose quadratic variation
is deterministic:
\begin{equation*}
[\xi^{(2)}_c(t),\xi^{(2)}_c(t)]=\kappa_4 t,
\end{equation*}
where $\kappa_4$ is the fourth cumulant of $\xi(1)$.
For any $f,g\in L^2(\mathbb{R})$, one has (see Farr\'e et al. \cite{FJ2010}) 
\begin{equation}\label{eq:xi^2_c quadratic var cov}
\E \left[\int_{\mathbb{R}} g(x) \xi^{(2)}_c(dx)\int_{\mathbb{R}} h(x)
\xi^{(2)}_c(dx)\right]=\kappa_4 \int_{\mathbb{R}}f(x)g(x)dx.
\end{equation}

 The decomposition  (\ref{eq:off diag decomp}) implies that
\[\E\left[\mathring{I}_k^{\xi}(f)\right]=  \int_{\mathbb{R}}f(x,x)dx.
\]

Notice that for any $f\in L^2(\mathbb{R}^2)$ and $g\in L^2(\mathbb{R})$
the following integrals, the first of which  is an off-diagonal double
integral and the second is a single integral with respect to Teugels
martingale $\xi_c^{(2)}(t)$:
\begin{equation}\label{eq:two integrals}
\int_{\mathbb{R}^2}'f(x_1,x_2)\xi(dx_1)\xi(dx_2)\quad\text{and}
\quad \int_{\mathbb{R}}g(x)\xi^{(2)}_c(dx).
\end{equation}
are uncorrelated (see Bai et al. \cite{BGT2}).

\subsection{Power counting theorems}
\label{pct}
Power counting theorems provide convergence conditions for some classes of integrals
on $\mathbb{R}^n$ whose integrands are products of functions bounded near zero by powers
of linearly dependent affine functionals and near infinity by different powers of those
functionals.
These theorems are useful in studying asymptotic distributions
of statistics of time series models with long-range dependence.
The results stated below were used in Terrin and Taqqu \cite{TT1}
to establish non-central limit theorems for
quadratic forms $Q_T$ of d.t.\  stationary processes with long-range dependence
(see Theorem \ref{NCT}).

First we introduce some notation. Let $\mathfrak{L}:= \{L_1,\ldots, L_m\}$ be
a collection of linear functionals $L_i(\bf x)$ on $\mathbb{R}^n$, $i = 1,\ldots, m$.
For numbers $0 < a_i \leq b_i$, $c_i > 0$, and real constants $\al_i$ and $\be_i$
we define
\beq \label{pct1}
P_1({\bf x}): = f_1(L_1({\bf x}))f_2(L_2({\bf x}))\cdots f_m(L_m({\bf x})),
\eeq
where the functions $f_i$, $(i = 1,\ldots, m)$ satisfy the condition:
\beq \label{pct2}
|f_i(y)|\leq
\left \{
\begin{array}{ll}
c_i|y_i|^{\al_i} & \mbox{if} \q |y|<a_i\\
c_i|y_i|^{\be_i} & \mbox{if} \q |y|>b_i,
\end{array}
\right.
\eeq
and $|f_i(y)|$ is bounded above in the interval $(a_i, b_i)$, $i=1,...,m$.
The constants $\al_i$ and $\be_i$ are called the {\sl exponents} of $|y|$
around $0$ and $\f$, respectively. For $W\subset T$ we set $s(W):= {\rm span}(W)\cap T$,
where ${\rm span}(W)$ denotes the linear span of $W$. Define
\beq \label{pct3}
d_0(W):=r(W)+\sum_{s(W)}\al_i \q \text{and} \q
d_\f(W):=r(T)-r(W)+\sum_{T\setminus s(W)}\be_i.
\eeq
A summation over a set $E$ means summation over the set $\{i: L_i\in E\}$,
$|E|$ denotes the cardinality of $E$ and $r(E)$ is the rank of $E$, that is,
the number of linearly independent functionals in $E$. We call a set $W$
{\sl padded} if for every $L$ in $W$, $L$ is also in $s(W\setminus\{L\})$,
that is, $L$  can be obtained as a linear combination of other elements in $W$.
Observe that $W=\emptyset$ is padded, and if $W$ is linearly independent,
then $d_0(W):=|W|+\sum_{s(W)}\al_i$.

The proofs of the following results can be found in Terrin and Taqqu \cite{TT1}
(see also Fox and Taqqu \cite{FT1}, and Terrin and Taqqu \cite{TT3}).
\begin{pp}\label{pct4}
If $r(T)=n$ and (a) $d_0(W)> 0$ for every nonempty subset $W$ of $T$ with $s(W)= W$,
while (b) $d_\f(W)< 0$ for every proper subset $W$ of $T$ with $s(W) = W$,
including the empty set. Then
\beq \label{pct5}
\int_{\mathbb{R}^n}|P_1({\bf x})|d{\bf x}<\f,
\eeq
where $P_1({\bf x})$ is as in \eqref{pct1}.
\end{pp}
Proposition \ref{pct4} can be extended to the class of functionals of
the form
$L_i(\bf x)+\theta_i$, where $\theta_i$ is a constant.
For $i=1,...,m$ let $\al_i$ and $\theta_i$ be real constants,
and let $M_i(\bf x)$ be a linear functional on $\mathbb{R}^n$.
Put $L_i(\bf x):=M_i(\bf x)+\theta_i$ and set
\beq \label{pct15}
P_2({\bf x}): = |L_1({\bf x})|^{\al_1}\cdots|L_m({\bf x})|^{\al_m}.
\eeq

\begin{pp}\label{pct6} 
Let $s(W)$ and $d_0(W)$ be as above.
Suppose that $d_0(W)> 0$ for every nonempty subset $W$ of $T$ with $s(W)= W$. Then
\beq \label{pct7}
\int_{[-t,t]^n}|P_2({\bf x})|d{\bf x}<\f\q \text{for all} \q t>0,
\eeq
where $P_2({\bf x})$ is as in \eqref{pct15}.
\end{pp}
\begin{rem}
\label{pct8}
{\rm If $\al_i > -1$ and $\be_i\geq-1$ in Proposition \ref{pct4}, then it suffices
to verify the conditions (a) and (b) for sets $W$ that are also padded.
If $\al_i > -1$ in Proposition \ref{pct6}, then it suffices to consider
subsets $W$ that are also padded (see Terrin and Taqqu \cite{TT1,TT3}).}
\end{rem}


\bigskip
\small
\noindent Mamikon S. Ginovyan:\\
Department of Mathematics and Statistics, Boston University,\\
111 Cummington Mall, Boston, MA 02215, USA\\ e-mail: ginovyan@math.bu.edu.\\
\\
Murad S. Taqqu:\\ 
Department of Mathematics and Statistics, Boston University,\\ 
111 Cummington Mall, Boston, MA 02215, USA\\ e-mail: murad@bu.edu.

\end{document}